\documentclass[a4paper,11pt]{article}
\usepackage[active]{srcltx}
\usepackage{amsmath}\usepackage{amssymb}\usepackage{amscd}

\newcommand{\pos}{{\mathrm{pos}}}\newcommand{\node}{{\mathrm{node}}}
\newcommand{\tsigma}{\tilde{\sigma}}
\newcounter{pic}
\newenvironment{pic}[1][\bf Fig. \arabic{pic}]{
        
        \refstepcounter{pic}\noindent\textbf{#1.}${}$\hspace{5pt}${}$\it}{}
\newtheorem{defi}{Definition}\newtheorem{theo}[defi]{Theorem}\newtheorem{prop}[defi]{Proposition}\newtheorem{lemm}[defi]{Lemma}      
\newtheorem{coro}[defi]{Corollary}   
\renewcommand{\theequation}{{\thesection}.{\arabic{equation}}}

\newcommand{\ml}[2]{{\mathcal{X}}_{#1}^{#2}}

\newcommand{\xl}[2]{X_{#1}^{#2}}

\newcommand{\yl}[2]{Y_{#1}^{#2}}

\newcommand{\lam}{\lambda^{(m)}}

\usepackage{empheq}

\begin{document}

$\ $

\begin{center}

\bigskip\bigskip\bigskip

{\Large\bf Jucys--Murphy elements  
and representations \\[.5cm]   of cyclotomic Hecke algebras}

\vspace{1.5cm}
{\large {\bf O. V. Ogievetsky\footnote{On leave of absence from P. N. Lebedev Physical Institute, Leninsky Pr. 53,
117924 Moscow, Russia} and L. Poulain d'Andecy}}

\vskip 1cm

{\large Center of Theoretical Physics, Luminy \\[.4em] 
13288 Marseille, France}

\end{center}

\vskip 1cm
\begin{abstract} 
An inductive approach to the representation theory of cyclotomic Hecke algebras, inspired by Okounkov and Vershik \cite{OV}, is developed.
We study the common spectrum of the Jucys-Murphy elements using representations of the simplest affine Hecke algebra.
Representations are constructed with the help of a new associative algebra whose underlying vector space is the tensor product of the cyclotomic 
Hecke algebra with the free associative algebra generated by standard $m$-tableaux.
\end{abstract}
\setcounter{tocdepth}{3}
\newpage
\tableofcontents

\vskip 1cm
\section{{\hspace{-0.55cm}.\hspace{0.55cm}}Introduction\vspace{.25cm}}\label{sectionintroduction}

The A-type Hecke algebra $H_n(q)$ is a one-parameter deformation of the group ring of the symmetric group $S_n$. We shall often omit the 
reference to the deformation parameter $q$ and write simply $H_n$ (for other families of 
algebras appearing in this article the reference to the parameters of the family will also be often omitted in the notation).
The Hecke algebra of type A plays an important role in numerous subjects: we just mention the knot theory, 
the Schur--Weyl duality for the
quantum general linear group; the representation theory of the A-type Hecke algebra for $q$ a root of unity (for generic $q$ 
the Hecke algebra is isomorphic to the group ring of $S_n$) is related to the modular representation theory of the symmetric group.

\vskip .2cm
The Hecke algebras $H_n$ form, with respect to $n$, an ascending chain of algebras. 
The Hecke algebras $H_n$ possess a set of 
Jucys--Murphy elements. This is a maximal commutative set (for generic $q$ and in the classical limit) whose 
advantages are: explicit description (compared to other maximal commutative sets discovered in the study of chain models);   simple relation
to the centralizers of the members of the chain. Moreover, the inductive formula for the Jucys--Murphy elements can be lifted to the
``universal" level: there exists a chain of affine Hecke algebras $\hat{H}_n$ for which one also defines the Jucys--Murphy elements;
the subalgebra generated by the Jucys--Murphy elements is isomorphic to a free commutative algebra and there is a surjection
$\hat{H}_n\rightarrow H_n$ which sends the Jucys--Murphy elements of the affine Hecke algebra to the Jucys--Murphy elements 
of the Hecke algebra $H_n$. 

\vskip .2cm
The main relation participating in the definition of the Hecke algebras is the Artin (or braid or Yang--Baxter) relation. The main additional
relation participating in the definition of the affine Hecke algebra and governing the commutativity of the set of the Jucys--Murphy elements is the 
reflection equation (a sort of braid relation participating in the presentation of the Weyl groups of type B) . The Jucys--Murphy elements and reflection equation appear in many applications: scattering on the half-line  \cite{C}, 
theory of knots in a torus (see, \emph{e.g.}, \cite{L} and references therein), the standard complex and the BRST operators for quantum Lie algebras 
\cite{GIO,IO4}, matrix integrals \cite{ZJ}, quantum multilinear algebra \cite{IOP}, construction of the quantum Minkowski space \cite{OSWZ,AKR} 
and quantum versions of the accidental isomorphisms \cite{JO} \emph{etc}.

\vskip .2cm
In \cite{OV} Okounkov and Vershik developed an inductive approach to the representation theory of the chain of the symmetric groups. Within this 
approach the description of the irreducible representations, the Young graph and the Young orthogonal form arise from the study of the spectrum of 
the set of the Jucys--Murphy elements of the group ring $\mathbb{C}S_n$. This approach has then been successfully generalized
to the projective representations of the symmetric group \cite{SV}, to the wreath product $G\wr S_n$ of any finite group $G$ by symmetric group \cite{Pus}, to the Hecke algebras of type A \cite{IO} and to the Birman-Murakami-Wenzl algebras of type A \cite{IO3}.

\vskip .2cm
An analogue of the Hecke algebra and of the braid group exists for all complex reflection groups. The complex reflection groups generalize the 
Coxeter groups and the complete list of irreducible finite complex reflection groups consists of the series of groups denoted $G(m,p,n)$, where 
$m,p,n$ are positive integers such that $p$ divides $m$, and $34$ exceptional groups \cite{ST}.

\vskip .2cm
The Hecke algebra of $G(m,1,n)$, which we denote by $H(m,1,n)$, is our main object of study. The algebra $H(m,1,n)$ has been introduced 
in \cite{AK,Bro-M,C2} and is called the cyclotomic Hecke algebra. For $m=1$ this is the Hecke algebra of type $A$ and for $m=2$ this is the Hecke 
algebra of type $B$. The representation theory of the algebra $H(m,1,n)$ was developed in \cite{AK} 
(and in \cite{H} for the Hecke algebra of type $B$); the irreducible representations of the algebra 
$H(m,1,n)$ are labeled, as for the group $G(m,1,n)$, by $m$-tuples of Young diagrams. 

\vskip .2cm
The Hecke algebras $H(m,1,n)$ also (as the usual Hecke algebras) form, with respect to $n$, an ascending chain of algebras. Moreover, the algebra
$H(m,1,n)$ is naturally a quotient of the affine Hecke algebra $\hat{H}_n$ and therefore inherits the natural set of Jucys--Murphy elements from the affine 
Hecke algebra. Our principal aim in this paper is to reproduce the representation theory of the cyclotomic Hecke algebra by analyzing the spectrum of 
Jucys--Murphy operators. We generalize the Okounkov-Vershik approach to the representation theory of the cyclotomic Hecke algebra $H(m,1,n)$; we 
construct irreducible representations and show that the usage of this approach allows to describe all irreducible representations of $H(m,1,n)$ upon certain restrictions 
(slightly stronger than the semi-simplicity conditions) on the parameters of the algebra $H(m,1,n)$. 

\vskip .2cm
We stress that our aim is not the construction itself of the representation theory - it has been already constructed in \cite{AK} - but we want 
to re-obtain the representations directly from the analysis 
of the Jucys--Murphy operators, to encode the representation bases in terms of sets of numbers which satisfy simple 
rules and which in fact are sets of common eigenvalues of the Jucys--Murphy elements and to reinterpret the Young multi-tableaux in terms of strings 
of eigenvalues of the Jucys--Murphy operators. As a byproduct of the construction of the representation bases it follows that the set of Jucys--Murphy 
operators is \emph{maximal} commutative (this observation is present in \cite{AK}). The approach, based on the Jucys--Murphy operators, has a recursive
nature - it uses the structure of the ascending, with respect to $n$, chain of the cyclotomic algebras $H(m,1,n)$.

\vskip .2cm
The inductive approach to the representation theory of the usual, type A, Hecke algebras, that is, the algebras $H(1,1,n)$, heavily uses the 
representation theory of the affine Hecke algebra of type A. 
One could expect that, say, for the representation theory of the algebras $H(2,1,n)$, there will be a necessity to use representations of the affine Hecke algebra of type B (see \cite{EK,M} for definitions).
But -- and it is maybe surprising -- the representation 
theory of the Hecke algebras $H(m,1,n)$ in the inductive approach requires, for all $m$, the study of representations of the same affine Hecke algebra of type A. 

\vskip .2cm
Since the cyclotomic Hecke algebra is the quotient of the affine Hecke algebra, a representation of the cyclotomic Hecke algebra is also
a representation of the affine Hecke algebra. The representations of the affine Hecke algebras are usually expressed in a different language, 
see \cite{BZ,Z,R,V,AS} for the original papers and surveys of the classical and $q$-deformed situations. 

\vskip .2cm
A novelty of our construction of the representations of the cyclotomic Hecke algebra is the introduction of a new associative algebra. Namely
we equip with an algebra structure the tensor product of the algebra $H(m, 1, n)$ with a free associative algebra generated by the standard 
$m$-tableaux corresponding to a given $m$-partition of $n$. We denote the resulting algebra by ${\mathfrak{T}}$. The representations are built 
then by evaluation from the right with the help of the simplest one-dimensional representation of $H(m, 1, n)$. There is a natural ``evaluation" 
homomorphism from ${\mathfrak{T}}$ to $H(m, 1, n)$ sending the generator labeled by a standard $m$-tableau to the corresponding primitive
idempotent of the algebra $H(m, 1, n)$.  

\vskip .2cm
An interesting consequence of the existence of this ``smash" product with the free algebra is a structure of a module on the tensor product of two 
representations corresponding to two (in general, any number of) $m$-partitions. 
We determine the rules of the decomposition of these tensor products into direct sums of irreducible representations. The decomposition rules 
themselves are quite easy: $V_{\lambda^{(m)}}\hat{\otimes}V_{\lambda'^{(m)}}\cong \dim(V_{\lambda'^{(m)}})\ V_{\lambda^{(m)}}$;
however an intertwiner, establishing this isomorphism, is difficult to describe explicitly - already for small $n$ the simplest choice 
of an intertwiner is quite evolved. Also, to obtain the decomposition rules we need the completeness result saying that every $H(m,1,n)$-module is 
isomorphic to a direct sum of $H(m,1,n)$-modules corresponding to $m$-partitions of $n$; however the definition of the module structure on the tensor 
product does not use the completeness. We think that the module structure on the tensor product deserves further study. 
 
\vskip .2cm
The space of a representation obtained with the help of the smash product is equipped with a distinguished basis which is analogous to the
semi-normal basis for representations of the symmetric group. It turns out that there exist several analogues of an invariant scalar product
on the representation spaces; one is bilinear, the definitions of others involve the complex conjugation as well as the involution $\omega\colon$ $q\to q^{-1}$, $v_j\to v_j^{-1}$, 
$j=1,\dots ,m$ ($q$ and 
$v_j$, $j=1,\dots ,m$, are the parameters of the algebra $H(m, 1, n)$, see Section \ref{sec-def} for precise definitions).  We compute these products for all irreducible representations of $H(m, 1, n)$
building thereby analogues of the orthogonal representations of the symmetric group. As a consequence, there is a large class of finite-dimensional
irreducible representations of the affine Hecke algebra which are unitarizable in a certain sense;
if the parameters $q$ and $v_j$, $j=1,\dots ,m$, of the algebra $H(m, 1, n)$ take values on the unit circle (in $\mathbb{C}$) then one of these products becomes a usual Hermitian scalar product. 

\vskip .2cm
Some of results, concerning representations of the cyclotomic Hecke algebras, smash product with free algebras, 
\emph{etc.}, of this paper were announced in \cite{OPdA} without proofs. We provide all necessary proofs here.

\subsection{Organization of the paper}$\ $

In Section \ref{sec-def} we recall the definitions of various chains of groups and algebras featuring in this article, and of the Jucys--Murphy 
elements of the chain of the braid groups and chains of quotients of the braid group ring.

\vskip .2cm
In Section \ref{ygcyclo} we start the study of the representation theory of the chain, with respect to $n$, of the cyclotomic Hecke algebras $H(m,1,n)$ 
generalizing Okounkov-Vershik approach of the representation theory of the symmetric groups. An important tool here is the list of representations,
satisfying some natural properties, of the affine Hecke algebra $\hat{H}_2$. We relate the set of standard Young $m$-tableaux with the 
set ${\mathrm{Spec}}(J_1,\dots,J_n)$ of common eigenvalues of the Jucys--Murphy elements $J_1,\dots,J_n$ in a certain class of representations
(which we call $C$-representations) of $H(m,1,n)$. 
More precisely, we show that any string of numbers belonging to ${\mathrm{Spec}}(J_1,\dots,J_n)$ is contained in a set 
called ${\mathrm{Cont}}_m(n)$ which is in bijection with the set of standard Young $m$-tableaux.

\vskip .2cm 
In Section \ref{ygcyclo'} we equip, for any $m$-partition $\lambda^{(m)}$ of $n$, the space $\mathbb{C}[\mathcal{X}_{\lambda^{(m)}}]\otimes H(m,1,n)$ 
with a structure of an associative algebra. Here $\mathbb{C}[\mathcal{X}_{\lambda^{(m)}}]$ is the free associative algebra with generators 
$\mathcal{X}_{\lambda^{(m)}}$  labeled by the standard $m$-tableaux of shape $\lambda^{(m)}$. To define the algebra structure, it is convenient 
to use the Baxterized form of (only) a 
part of generators of the cyclotomic Hecke algebra. Given a one-dimensional representation of the algebra $H(m,1,n)$, we construct, with the help of 
the algebra structure on $\mathbb{C}[\mathcal{X}_{\lambda^{(m)}}]\otimes H(m,1,n)$, a representation on the space whose basis is labeled by the 
standard $m$-tableaux of shape $\lambda^{(m)}$. This construction implies  
that the set ${\mathrm{Spec}}(J_1,\dots,J_n)$ of common eigenvalues of the Jucys--Murphy elements in $C$-representations
coincides with the set ${\mathrm{Cont}}_m(n)$ and with the set of standard Young $m$-tableaux. 
At the end of Section \ref{ygcyclo'} we compute analogues of the invariant scalar product for all representations $V_{\lambda^{(m)}}$. 

\vskip .2cm
In Section \ref{sec-comp} we complete the representation theory of the cyclotomic Hecke algebras; we show that the constructed representations are 
irreducible and pairwise non-isomorphic (the proof is included for completeness; it is adopted from \cite{AK}). Using an upper bound for the dimension of the cyclotomic Hecke algebra and some results about products of Bratteli diagrams recalled in Appendix B, we conclude in 
a standard way that the class of irreducible $C$-representations exhausts the set of the irreducible representations of the cyclotomic Hecke algebra 
when the parameters of the algebra satisfy the restrictions specified in Section \ref{sec-def}.

\vskip .2cm
Further, we include in Section \ref{sec-comp} several direct consequences of the developed representation theory (valid either in the generic picture or 
under the restrictions on the parameters considered in this article): the semi-simplicity of the cyclotomic Hecke algebra $H(m,1,n)$, the simplicity of the 
branching rules for the representations of the chain of the cyclotomic Hecke algebras\footnote{The paper  \cite{GV} advocates the point of view that it is the affine 
Hecke algebra ``that is responsible for the multiplicity one phenomena"; for generic finite-dimensional representations the multiplicity one statement for the affine 
Hecke algebra of type A 
follows from  the multiplicity one statement for the cyclotomic Hecke algebra because in a generic finite-dimensional representation the 
spectrum of the Jucys--Murphy element $J_1$ is finite and thus such representation is actually a representation of a cyclotomic quotient of the affine Hecke algebra.} 
and the maximality in $H(m,1,n)$ of the commutative set formed by the Jucys--Murphy elements.  We also mention some information, implied by the 
developed representation theory,  about the structure of the centralizer of the algebra $H(m,1,n-1)$ considered as a subalgebra in $H(m,1,n)$.

\vskip .2cm
The article concludes with Appendices A, B and C.

\vskip .2cm 
In Appendix A we explain how the algebra structure on the space $\mathbb{C}[\mathcal{X}_{\lambda^{(m)}}]\otimes H(m,1,n)$  induces the 
tensor structure on the set of $C$-representations of the cyclotomic Hecke algebra $H(m,1,n)$; more generally, given an $H(m,1,n)$-module $W$, 
the algebra structure on $\mathbb{C}[\mathcal{X}_{\lambda^{(m)}}]\otimes H(m,1,n)$ leads to a structure of an $H(m,1,n)$-module on the space 
$V_{\lambda^{(m)}}\otimes W$, where $V_{\lambda^{(m)}}$ is the $C$-representation related to the $m$-partition $\lambda^{(m)}$. 
We determine the rules of the decomposition of the tensor products of $C$-representations into direct sums of irreducible representations. 
In the course of the proof we give several explicit examples of such decompositions.

\vskip .2cm
In Appendix B we recall some definitions and results concerning Bratteli diagrams and their products. We specify the information to the powers 
of the Young graph; the $m^{\textrm{th}}$ power of the Young graph is relevant in the representation theory of the chain of the cyclotomic
Hecke algebras. Appendix B
has a review character. 

\vskip .2cm
Appendix C contains several examples of the defining relations of the algebra on the tensor product of the algebra $H(m,1,n)$ with a free 
associative algebra generated by the standard $m$-tableaux; resulting 
explicit formulas for matrix elements of generators in low-dimensional  irreducible representations of the cyclotomic Hecke algebras $H(m,1,n)$
are given.

\paragraph{Notation.} $\ $

\vskip .2cm
In this article, the ground field is the field $\mathbb{C}$ of complex numbers. 

The spectrum of an operator ${\cal{T}}$ is denoted by ${\mathrm{Spec}}({\cal{T}})$.

We denote, for two integers $k,l\in\mathbb{Z}$ with $k<l$, by  $[k,l]$ the set of integers $\{k,k+1,\dots,l-1,l\}$.

The $q$-number $j_q$ is defined by $j_q:=\displaystyle{\frac{q^j-q^{-j}}{q-q^{-1}}}$. 

The diagonal matrix with entries $z_1,z_2,\dots,z_k$ (on the diagonal) is denoted by $\text{diag}(z_1,z_2,\dots,z_k)$.

\setcounter{equation}{0}
\section{{\hspace{-0.55cm}.\hspace{0.55cm}}Cyclotomic Hecke algebras and Jucys--Murphy elements}\label{sec-def}

The braid group $B_n$ 
of type A (or simply the braid group) on $n$ strands is generated by the elements $\sigma_1$, $\dots$, $\sigma_{n-1}$ with the defining relations: 
\begin{empheq}[left=\empheqlbrace]{alignat=4}
\label{def1a}&\sigma_i\sigma_{i+1}\sigma_i=\sigma_{i+1}\sigma_i\sigma_{i+1} &\hspace{1cm}& \textrm{for all $i=1,\dots,n-2$\ ,}\\[.5em]
\label{def1b}&\sigma_i\sigma_j=\sigma_j\sigma_i && \textrm{for all $i,j=1,\dots,n-1$ such that $|i-j|>1$\ .}
\end{empheq}

The braid group $\alpha B_n$ of type B (sometimes called \emph{affine} braid group)
is obtained by adding to the generators  $\sigma_1$, $\dots$, $\sigma_{n-1}$ the generator $\tau$ with the defining relations (\ref{def1a}), (\ref{def1b}) and:  
\begin{empheq}[left=\empheqlbrace]{alignat=3}
\label{def1'a}&\tau\sigma_1\tau\sigma_1=\sigma_1\tau\sigma_1\tau\ ,&\\[.5em]
\label{def1'b}&\tau\sigma_i=\sigma_i\tau & \textrm{for $i>1$\ .}
\end{empheq}

The elements $J_i$, $i=1,\dots,n$, of the braid group of type B defined inductively by the following initial condition and recursion:
\begin{equation}\label{JM}J_1=\tau\ ,\quad J_{i+1}=\sigma_iJ_i\sigma_i\ ,\end{equation}
are called Jucys--Murphy elements.
It is well known that they form a commutative set of elements. In addition,  $J_i$ commutes with all $\sigma_k$ except $\sigma_{i-1}$ and 
$\sigma_i$,
\begin{equation}\label{JMcomm}J_i\sigma_k=\sigma_k J_i\ \ \ {\textrm{if}}\ \  k>i\ \  {\textrm{or}}\ \  k<i-1\ .\end{equation} 

The affine Hecke algebra $\hat{H}_n$ is the quotient of the group algebra of the B-type braid group $\alpha B_n$ 
by:
\begin{equation}\label{def1''b}
\sigma_i^2=(q-q^{-1})\sigma_i+1\qquad  \textrm{for all $i=1,\dots,n-1$\ .}
\end{equation}
The usual Hecke algebra $H_n$ is the algebra generated by the elements $\sigma_1$, $\dots$, $\sigma_{n-1}$ with the relations
(\ref{def1a})--(\ref{def1b}) and (\ref{def1''b}).

\vskip .2cm
The cyclotomic Hecke algebra $H(m,1,n)$ is the quotient of the affine Hecke algebra $\hat{H}_n$ by
\begin{equation}\label{def1''a}
(\tau-v_1)\dots(\tau-v_m)=0\ .\end{equation}
In particular, $H(1,1,n)$ is isomorphic to the Hecke algebra of type A and $H(2,1,n)$ is 
isomorphic to the Hecke algebra of type B. 
The subalgebra of $H(m,1,n)$ generated by $\sigma_1,\dots,\sigma_{n-1}$ is isomorphic to the Hecke algebra of type A.

\vskip .2cm
The algebra $H(m,1,n)$ is a deformation of the group algebra $\mathbb{C}G(m,1,n)$ of the complex reflection group $G(m,1,n)$.
The group $G(m,1,n)$ is isomorphic to $S_n\wr C_m$, the wreath product of the cyclic group with $m$ elements by the symmetric group $S_n$.

\vskip .2cm
The deformation from $\mathbb{C}G(m,1,n)$ to $H(m,1,n)$ is flat in the sense that  $H(m,1,n)$ is a free 
 $\mathbb{C}[q,q^{-1},v_1,\dots ,v_m]$-module of dimension equal to the order of $G(m,1,n)$, that is:
\begin{equation}\label{dim}
\dim(H(m,1,n))=n!m^n\ .
\end{equation}
The flatness is proved in \cite {AK} with the help of the representation theory. In \cite{OPdA2} the flatness is proved within the theory of associative algebras. In this paper, we will only use the inequality
\begin{equation}\label{dim-ine}
\dim(H(m,1,n))\leq n!m^n\ ,
\end{equation}
which can be proved without use of the representation theory \cite{AK,Bre-M}.

\vskip .2cm
The specialization  of $H(m,1,n)$ is semi-simple if and only if the numerical values of the parameters satisfy (see \cite{Ari-sim})
\begin{equation}\label{sesi1}
1+q^2+\dots+q^{2N}\neq 0\ \ {\textrm{for all}}\ \ N:\ N<n\ \end{equation}
and
\begin{equation}\label{sesi2}
q^{2i}v_j-v_k\neq 0\ \ {\textrm{for all}}\ \ i,j,k \ \ {\textrm{such that}}\ \  j\neq k \ {\textrm{and}}\ -n<i<n\ .\end{equation}
In the sequel we work either with a generic cyclotomic Hecke algebra (that is, $v_1$, $\dots$, $v_m$ and $q$ are indeterminates) or in the semi-simple situation with an additional requirement:
\begin{equation}\label{sesi3} v_j\neq 0\ ,\ j=1,\dots,m\ .\end{equation}

As $n$ varies, the algebras $H(m,1,n)$ form an ascending chain of algebras:
\begin{equation}\label{chaine}H(m,1,0)= \mathbb{C}
\subset H(m,1,1)\subset\dots\subset H(m,1,n)\subset\dots\end{equation}
(the elements $\tau$ and $\sigma_1,\dots ,\sigma_{n-2}$ of the algebra $H(m,1,n)$ generate a subalgebra isomorphic to 
$H(m,1,n-1)$).
One has similar ascending chains of braid groups, affine braid groups and affine Hecke algebras.
Thus the reference to $n$ (as in $H_n$, $H(m,1,n)$, etc.) in the notation for the generators can be omitted.

\vskip .2cm
The representation theory of the generic algebra $H(m,1,n)$ was studied in \cite{AK}. Here we present another approach which is a generalization 
of the approach of Okounkov and Vershik to the representation theory of the symmetric group \cite{OV} and which refers to the ascending chain
(\ref{chaine}).

\vskip .2cm
We shall denote by the same symbols $J_i$ the images of the Jucys--Murphy elements in the cyclotomic Hecke algebra. As a by-product of 
the representation theory of the generic algebra $H(m,1,n)$, the set of  the Jucys--Murphy elements $\{ J_1,\dots ,J_n\}$ is maximal commutative in 
$H(m,1,n)$; more precisely, the algebra of polynomials in the Jucys--Murphy elements coincides with the algebra generated by the union of the centers 
of $H(m,1,k)$ for $k=1,\dots,n$.

\paragraph{Remark.} We use the same notation ``$\sigma_i$" for generators of different groups and algebras: these are braid and affine braid groups and 
Hecke, affine Hecke and cyclotomic Hecke algebras. The symbol $\tau$ is also used to denote a generator of several different objects.
This should not lead to any confusion, it will be clear from the context what is the algebra/group in question.

\setcounter{equation}{0}
\section{{\hspace{-0.55cm}.\hspace{0.55cm}}Spectrum of Jucys--Murphy elements and Young $m$-tableaux}\label{ygcyclo}

We begin to develop an approach, based on the Jucys--Murphy elements, to the representation theory of the chain (with respect to $n$) of the 
cyclotomic Hecke algebras $H(m,1,n)$. This is a generalization of the approach of \cite{OV}.

\paragraph{1.} The first step consists in construction of all representations of $H(m,1,n)$ verifying two conditions. First, the Jucys--Murphy elements  
$J_1,\dots,J_n$ are represented by semi-simple (diagonalizable) operators. Second, for every $i=1,\dots,n-1$ the action of the subalgebra generated 
by $J_i$, $J_{i+1}$ and $\sigma_i$ is completely reducible. We shall use the name $C$-representations ($C$ is the first letter in ``completely reducible") 
for these representations.  At the end of the construction we shall see that all irreducible representations of $H(m,1,n)$ are $C$-representations.

\vskip .2cm
Following \cite{OV} we denote by ${\mathrm{Spec}}(J_1,\dots,J_n)$ the set of strings of eigenvalues of the Jucys--Murphy elements in the set 
of $C$-representations: $\Lambda=(a^{(\Lambda)}_1,\dots,a^{(\Lambda)}_n)$ belongs to ${\mathrm{Spec}}(J_1,\dots,J_n)$ if there is a vector 
$e_{\Lambda}$ in the space of some $C$-representation such that 
$J_i(e_{\Lambda})=a^{(\Lambda)}_ie_{\Lambda}$ for all $i=1,\dots,n$. Every $C$-representation possesses a basis formed by vectors $e_{\Lambda}$
(this is a reformulation of the first condition in the definition of $C$-representations).
Since $\sigma_k$ commutes with $J_i$ for $k>i$ and $k<i-1$, the action of $\sigma_k$ on a vector  $e_{\Lambda}$, 
 $\Lambda\in {\mathrm{Spec}}(J_1,\dots,J_n)$, is ``local" in the sense that $\sigma_k(e_{\Lambda})$ is a linear combination of $e_{\Lambda'}$ such that 
$a^{(\Lambda')}_i=a^{(\Lambda)}_i$ for $i\neq k,k+1$.

\paragraph{2. Affine Hecke algebra $\hat{H}_2$.} Consider the affine Hecke algebra $\hat{H}_2$, generated by $X$, $Y$ and $\sigma$ with the relations: 
\begin{equation}\label{affHec}XY=YX,\quad Y=\sigma X\sigma,\quad \sigma^2=(q-q^{-1})\sigma+1\ .\end{equation}
For all $i=1,\dots,n-1$, the subalgebra of $H(m,1,n)$ generated by $J_i$, $J_{i+1}$ and $\sigma_i$ is a quotient of $\hat{H}_2$. 
We reproduce here the result of \cite{IO} concerning the classification of irreducible representations with diagonalizable $X$ and $Y$ of the algebra $\hat{H}_2$.

\vskip .2cm
There are one-dimensional and two-dimensional irreducible representations.
\begin{itemize}
 \item The one-dimensional irreducible representations are given by
\begin{equation}\label{mat-d1}
X\mapsto a,\quad Y\mapsto q^{2\varepsilon}a,\quad\sigma\mapsto \varepsilon q^{\varepsilon}\ ,\ \ \text{where}\ \ \varepsilon =\pm 1\ .
\end{equation}
\item The two-dimensional irreducible representations are given by
\[\sigma\mapsto\left(\begin{array}{cc}0 & 1\\ 1 & q-q^{-1}\end{array}\right),\quad X\mapsto\left(\begin{array}{cc}a & -(q-q^{-1})b\\ 0 & b\end{array}\right),\quad Y\mapsto\left(\begin{array}{cc}b & (q-q^{-1})b\\ 0 & a\end{array}\right),\]
with $b\neq a$ in order for $X$ and $Y$ to be diagonalizable and with $b\neq q^{\pm2}a$ to ensure irreducibility. By a change of basis we transform $X$ and $Y$ to a diagonal form:
\begin{equation}\label{mat-d2}
\sigma\mapsto\left(\begin{array}{cc}\frac{(q-q^{-1})b}{b-a}\ &\ 1-\frac{(q-q^{-1})^2ab}{(b-a)^2}\\ 1 \ &\ -\frac{(q-q^{-1})a}{b-a}\end{array}\right),\quad X\mapsto\left(\begin{array}{cc}a &0\\ 0 & b\end{array}\right),\quad Y\mapsto\left(\begin{array}{cc}b & 0\\ 0 & a\end{array}\right).
\end{equation}
\end{itemize}

\paragraph{3.} We return to strings of eigenvalues of the Jucys--Murphy elements.

\begin{prop}
{\hspace{-.2cm}.\hspace{.2cm}}
 \label{prop2}
Let $\Lambda=(a_1,\dots,a_i,a_{i+1},\dots,a_n)\in {\mathrm{Spec}}(J_1,\dots,J_n)$ and let  $e_{\Lambda}$ be a corresponding vector. Then 
\begin{itemize}
\item[(a)] We have $a_i\neq a_{i+1}$.
\item[(b)] If $a_{i+1}=q^{2\varepsilon}a_i$, where $\varepsilon =\pm 1$, then $\sigma_i(e_{\Lambda})=\varepsilon q^{\varepsilon}e_{\Lambda}$. 
\item[(c)] If $a_{i+1}\neq q^{\pm2}a_i$ then 
$\Lambda'=(a_1,\dots,a_{i+1},a_i,\dots,a_n)\in {\mathrm{Spec}}(J_1,\dots,J_n)$; moreover, the vector 
$\sigma_i(e_{\Lambda})-\frac{(q-q^{-1})a_{i+1}}{a_{i+1}-a_i}e_{\Lambda}$ corresponds to the string $\Lambda'$ (see (\ref{mat-d2}) with $b=a_{i+1}$ and $a=a_i$). 
\end{itemize}
\end{prop}

The proof follows directly from the representation theory of the algebra $\hat{H}_2$, described above  (cf. the Proposition 4.1 in \cite{OV} and the 
Proposition 3 in \cite{IO}). 

\paragraph{4. Content strings.}

\begin{defi}
{\hspace{-.2cm}.\hspace{.2cm}}
 \label{def-cont} A content string $(a_1,\dots,a_n)$ is a string of numbers satisfying the following conditions:
\begin{itemize}
\item[(c1)]  $a_1\in\{v_1,\dots ,v_m\}$;
\item[(c2)]  for all $j>1$: if $a_j=v_k q^{2z}$ for some $k$ and $z\neq 0$ then $\{v_k q^{2(z-1)},v_k q^{2(z+1)}\}\cap\{a_1,\dots,a_{j-1}\}\neq\varnothing$;
\item[(c3)]  if $a_i=a_j=v_kq^{2z}$ with $i<j$ for some $k$ and $z$, then $\{v_k q^{2(z-1)},v_k q^{2(z+1)}\}\subset\{a_{i+1},\dots,a_{j-1}\}$.
\end{itemize}

\vskip .2cm
The set of content strings of length $n$ we denote by ${\mathrm{Cont}}_m(n)$.
\end{defi}

Here is the ``cyclotomic" analogue of the Theorem 5.1 in \cite{OV} and the Proposition 4 in \cite{IO}. We adapt the proof paying attention
to places where the restrictions (\ref{sesi1})--(\ref{sesi3}) are essential.

\begin{prop}
{\hspace{-.2cm}.\hspace{.2cm}}
 \label{prop3}
If a string of numbers $(a_1,\dots,a_n)$ belongs to ${\mathrm{Spec}}(J_1,\dots,J_n)$ then it belongs to ${\mathrm{Cont}}_m(n)$.
\end{prop}

\emph{Proof.} Since $J_1=\tau$ the condition (c1) follows from  the characteristic equation for $\tau$. 

\vskip .2cm
Assume that (c2) is not true, that is, there is a string $(a_1,\dots,a_n)\in {\mathrm{Spec}}(J_1,\dots,J_n)$ such that for some $j>1$, some $k$ and some 
$z\neq0$ one has $a_j=v_kq^{2z}$ but $a_i\neq v_k q^{2(z-1)}$ and $a_i\neq v_k q^{2(z+1)}$ for all $i$ smaller than $j$. By 
successive applications of the statement (c) of the Proposition \ref{prop2} we obtain an element of ${\mathrm{Spec}}(J_1,\dots,J_n)$ with $v_kq^{2z}$ at the 
first position. The restrictions (\ref{sesi1})--(\ref{sesi3}) imply $v_kq^{2z}\neq v_i$ for all $i=1,\dots,m$ and this contradicts the condition (c1).

\vskip .2cm
We prove the condition (c3) by induction on $j-i$. The base of induction is the statement (a) of the Proposition \ref{prop2}. Assume that there are 
some $i$ and some $j$ such that $i<j$ and $a_i=a_j=v_kq^{2z}$ for some string $(a_1,\dots,a_n)\in {\mathrm{Spec}}(J_1,\dots,J_n)$. By induction 
we suppose that the condition (c3) is verified for all $i',j'$ such that 
$|j'-i'|<j-i$. If $\{v_k q^{2(z-1)},v_k q^{2(z+1)}\}\cap\{a_{i+1},\dots,a_{j-1}\}=\varnothing$ then by an application of the statement (c) of the
Proposition \ref{prop2} we move $a_j$ to the left to the position number $(j-1)$ (note that $(j-1)$ is still greater than $i$ by the statement (a) of the
Proposition \ref{prop2}) and obtain an element of ${\mathrm{Spec}}(J_1,\dots,J_n)$ which contradicts the induction hypothesis.
Now assume that only one element from the set $\{v_k q^{2(z-1)},v_k q^{2(z+1)}\}$ is present in $\{a_{i+1},\dots,a_{j-1}\}$. 
By the induction hypothesis, this element appears only once in $\{a_{i+1},\dots,a_{j-1}\}$. If $j-i>2$ then, by an application of the  
statement (c) of the Proposition \ref{prop2}, we obtain an element of ${\mathrm{Spec}}(J_1,\dots,J_n)$ contradicting  the induction hypothesis. Thus $j-i=2$
which is impossible because the braid relation $\sigma_i\sigma_{i+1}\sigma_i=\sigma_{i+1}\sigma_i\sigma_{i+1}$ is incompatible with assignments 
$\sigma_i\mapsto \varepsilon q^{\varepsilon}$ and $\sigma_{i+1}\mapsto -\varepsilon q^{-\varepsilon}$, where $\varepsilon =\pm 1$ (these values are implied by the statement (b) of the Proposition 
\ref{prop2}).\hfill$\square$

\paragraph{Remark.} Let $E_S$ be a vector space with a basis $\{ e_{\Lambda}\}$ whose vectors are labeled by the elements 
$\Lambda\in{\mathrm{Spec}}(J_1,\dots,J_n)$. Let $E_C$ be a vector space with a basis $\{ e_{\mu}\}$ whose vectors are labeled by the elements 
$\mu\in{\mathrm{Cont}}_m(n)$. By the Proposition \ref{prop3}, ${\mathrm{Spec}}(J_1,\dots,J_n)\subset {\mathrm{Cont}}_m(n)$ and so
$E_S$ is naturally a vector subspace of $E_C$. The space $E_S$ is equipped with the action of the algebra $H(m,1,n)$:
the operator corresponding to the generator $\tau$ is simply $J_1$;  
the precise formulas for the action of the elements $\sigma_i$, $i=1,\dots,n-1$, are given in the Proposition \ref{prop2}. 

\vskip .2cm
The  Definition \ref{def-cont} straightforwardly implies that if $(a_1,\dots,a_i,a_{i+1},\dots,a_n)\in{\mathrm{Cont}}_m(n)$ with $a_{i+1}\neq q^{\pm2}a_i$ then 
$(a_1,\dots,a_{i+1},a_i,\dots,a_n)\in{\mathrm{Cont}}_m(n)$. 
Therefore the operators corresponding to the generators $\tau$ and $\sigma_i$, $i=1,\dots,n-1$, make sense as the operators in the space $E_C$.
The aim of Subsection \ref{suserep} below is to show that these operators continue to define the action of the algebra $H(m,1,n)$ - now on the,
in general, bigger space $E_C$. One could construct representations working directly with strings but it is illuminating and convenient to interpret the 
set ${\mathrm{Cont}}_m(n)$ in the more geometric terms of Young multi-tableaux.

\vskip .2cm
At the end of the whole construction it will then follow, see Section \ref{sec-comp}, that the spaces $E_S$ and $E_C$ actually coincide.

\paragraph{5.} Using ``intertwining" operators $U_{i+1}:=\sigma_iJ_i-J_i\sigma_i$, $i=1,\dots,n-1$, it can be proved, as in \cite{IO}, that:
\begin{equation}\label{respej} {\mathrm{Spec}}(J_{i+1})\subset {\mathrm{Spec}}(J_i)\cup q^{\pm2}\cdot {\mathrm{Spec}}(J_i).\end{equation}
Since ${\mathrm{Spec}}(J_1)\subset\{v_1,\dots,v_m\}$ we arrive at the following conclusion.

\begin{prop}\label{cospe}
{\hspace{-.2cm}.\hspace{.2cm}}
{}For all $i=1,\dots,n$,
\begin{equation}\label{respej2}  {\mathrm{Spec}}(J_i)\subset\{v_k q^{2\, [1-i,i-1]},\,k=1,\dots,m\}\ .\end{equation}
\end{prop}

The Proposition \ref{cospe} follows also from the Proposition \ref{prop2} and the Proposition \ref{prop3}. Indeed assume that for a string 
$(a_1,\dots,a_n)\in {\mathrm{Spec}}(J_1,\dots,J_n)$ the Proposition \ref{cospe} does not hold. Let $i$ be the smallest integer for which 
$a_i\notin\{v_k q^{2\, [1-i,i-1]},\,k=1,\dots,m\}$. Using the statement (c) of the Proposition \ref{prop2} we move $a_i$ to the left until it reaches the first 
position in the string and obtain thereby an element of ${\mathrm{Spec}}(J_1,\dots,J_n)$ which does not verify the condition (c1). This contradicts 
the Proposition \ref{prop3}.

\paragraph{6. Young $m$-diagrams and $m$-tableaux.} A Young $m$-diagram, or $m$-partition, is an $m$-tuple of Young diagrams $\lambda^{(m)}=(\lambda_1,\dots,\lambda_m)$. The size of a Young diagram $\lambda$ 
is the number of nodes of the diagram and is denoted by $|\lambda|$. By definition the size of an $m$-tuple $\lambda^{(m)}=(\lambda_1,\dots,\lambda_m)$ is 
\begin{equation}\label{lenmdi}|\lambda^{(m)}|:=|\lambda_1|+\dots+|\lambda_m|\ .\end{equation}

We recall some standard terminology. For a usual partition $\lambda$, a node $\alpha$ is called {\it removable} if the set of nodes 
obtained from $\lambda$ by removing the node $\alpha$ is still a partition; a node $\beta$ not in $\lambda$ is called {\it  addable} if the set of nodes 
obtained from $\lambda$ by adding $\beta$ is still a partition. 

\vskip .2cm
We extend this terminology for the $m$-partitions. To this end we define the notion of an \emph{$m$-node}: an $m$-node $\alpha^{(m)}$ is a 
pair $(\alpha,p)$ 
consisting of a usual node $\alpha$ and an integer $p$ with $1\leq p\leq m$. We will refer to $\alpha$ as the {\it node of the $m$-node $\alpha^{(m)}$} 
and we will write $\node(\alpha^{(m)})=\alpha$; we will refer to the integer $p$ as the {\it position of the $m$-node $\alpha^{(m)}$} and we will note 
$\pos(\alpha^{(m)})=p$. With this definition, an $m$-partition $\lambda^{(m)}$ is a set of $m$-nodes such that, for any $p$ between $1$ and $m$, the 
subset consisting of $m$-nodes $\alpha^{(m)}$ with $\pos(\alpha^{(m)})=p$ forms a usual partition.

\vskip .2cm
Let $\lambda^{(m)}$ be an $m$-partition. An $m$-node $\alpha^{(m)}$ of $\lambda^{(m)}$ is called {\it removable} if the set of $m$-nodes obtained 
from $\lambda^{(m)}$ by removing $\alpha^{(m)}$ is still an $m$-partition. An $m$-node $\beta^{(m)}$ not in $\lambda^{(m)}$ is called {\it  addable} if the set of 
$m$-nodes obtained from $\lambda^{(m)}$ by adding $\beta^{(m)}$ is still an $m$-partition. The $m$-partition obtained from $\lambda^{(m)}$ by removing any 
removable $m$-node $\alpha^{(m)}$ will be denoted by $\lambda^{(m)}\backslash\{\alpha^{(m)}\}$. For any $m$-partition $\lambda^{(m)}$, we denote by 
${\cal{E}}_-(\lambda^{(m)})$ the set of removable $m$-nodes of $\lambda^{(m)}$ and by ${\cal{E}}_+(\lambda^{(m)})$ the set of addable $m$-nodes of 
$\lambda^{(m)}$.

\vskip .2cm
An $m$-partition whose $m$-nodes are filled with numbers is called $m$-tableau. 

\vskip .2cm
Let the size of the $m$-partition $\lambda^{(m)}$ be $n$. We place now the numbers $1,\dots,n$ in the $m$-nodes of $\lambda^{(m)}$ 
in such a way that in every diagram the numbers in the $m$-nodes are in ascending order along rows and columns in right and down directions. This is a 
\emph{standard} Young $m$-tableau of shape $\lambda^{(m)}$. 

\vskip .2cm 
We associate to each $m$-node of a Young $m$-diagram a number (the ``content") which is $v_kq^{2(s-r)}$ for the $m$-node $\alpha^{(m)}$ such that $\pos(\alpha^{(m)})=k$ and $\node(\alpha^{(m)})$ lies in the line $r$ and column $s$ (equivalently we could say that the $m$-node $\alpha^{(m)}$ lies in the line $r$ and column $s$ of the $k^{\textrm{th}}$ diagram of the $m$-diagram). 
Note that the notion of content makes sense for an arbitrary $m$-node of an arbitrary set of $m$-nodes. 

\vskip .2cm
{}For an arbitrary set of $m$-nodes, two $m$-nodes on a same diagonal of the same diagram have the same 
content which allows to speak about the ``content number'' of a diagonal. 

\vskip .2cm
Here is an example of a standard Young $m$-tableau with $m=2$ and $n=10$ (the contents of the $m$-nodes are indicated):

\begin{equation}\label{example}\left(\begin{array}{cc}
\setlength{\unitlength}{2200sp}
\begingroup\makeatletter\ifx\SetFigFont\undefined
\gdef\SetFigFont#1#2#3#4#5{\reset@font\fontsize{#1}{#2pt}  \fontfamily{#3}\fontseries{#4}\fontshape{#5}  \selectfont}\fi\endgroup
\begin{picture}(4000,2904)(2000,-4888)
\thinlines
{\put(2566,-2986){\framebox(945,945){}}}
{\put(3511,-2986){\framebox(945,945){}}}
{\put(3511,-3931){\framebox(945,945){}}}
{\put(2566,-3931){\framebox(945,945){}}}
{\put(4456,-2986){\framebox(945,945){}}}
{\put(2566,-4876){\framebox(945,945){}}
}\put(3356,-2266){\makebox(0,0)[lb]{\smash{{\SetFigFont{12}{14.4}{\rmdefault}{\mddefault}{\updefault}{\scriptsize{2}}}}}}
\put(2411,-2266){\makebox(0,0)[lb]{\smash{{\SetFigFont{12}{14.4}{\rmdefault}{\mddefault}{\updefault}{\scriptsize{1}}}}}}
\put(2411,-3211){\makebox(0,0)[lb]{\smash{{\SetFigFont{12}{14.4}{\rmdefault}{\mddefault}{\updefault}{\scriptsize{6}}}}}}
\put(3356,-3211){\makebox(0,0)[lb]{\smash{{\SetFigFont{12}{14.4}{\rmdefault}{\mddefault}{\updefault}{\scriptsize{9}}}}}}
\put(2411,-4156){\makebox(0,0)[lb]{\smash{{\SetFigFont{12}{14.4}{\rmdefault}{\mddefault}{\updefault}{\scriptsize{7}}}}}}
\put(3581,-3616){\makebox(0,0)[lb]{\smash{{\SetFigFont{12}{14.4}{\rmdefault}{\mddefault}{\updefault}{\small{$v_1$}}}}}}
\put(2521,-3616){\makebox(0,0)[lb]{\smash{{\SetFigFont{12}{14.4}{\rmdefault}{\mddefault}{\updefault}{\small{$v_1q^{-2}$}}}}}}
\put(2521,-4561){\makebox(0,0)[lb]{\smash{{\SetFigFont{12}{14.4}{\rmdefault}{\mddefault}{\updefault}{\small{$v_1q^{-4}$}}}}}}
\put(3446,-2716){\makebox(0,0)[lb]{\smash{{\SetFigFont{12}{14.4}{\rmdefault}{\mddefault}{\updefault}{\small{$v_1q^2$}}}}}}
\put(4351,-2716){\makebox(0,0)[lb]{\smash{{\SetFigFont{12}{14.4}{\rmdefault}{\mddefault}{\updefault}{\small{$v_1q^4$}}}}}}
\put(2606,-2716){\makebox(0,0)[lb]{\smash{{\SetFigFont{12}{14.4}{\rmdefault}{\mddefault}{\updefault}{\small{$v_1$}}}}}}
\put(4301,-2266){\makebox(0,0)[lb]{\smash{{\SetFigFont{12}{14.4}{\rmdefault}{\mddefault}{\updefault}{\scriptsize{4}}}}}}
\end{picture}
\text{{\Large\bf ,}} & 
\setlength{\unitlength}{2200sp}
\begingroup\makeatletter\ifx\SetFigFont\undefined
\gdef\SetFigFont#1#2#3#4#5{\reset@font\fontsize{#1}{#2pt}  \fontfamily{#3}\fontseries{#4}\fontshape{#5}
\selectfont}\fi\endgroup
\begin{picture}(4000,2904)(6000,-4888)
\thinlines
{\put(7921,-2941){\framebox(945,945){}}
}{\put(6976,-3886){\framebox(945,945){}}
}{\put(8866,-2941){\framebox(945,945){}}
}{\put(6976,-2941){\framebox(945,945){}}
}\put(6821,-2221){\makebox(0,0)[lb]{\smash{{\SetFigFont{12}{14.4}{\rmdefault}{\mddefault}{\updefault}{\scriptsize{3}}}}}}
\put(7766,-2221){\makebox(0,0)[lb]{\smash{{\SetFigFont{12}{14.4}{\rmdefault}{\mddefault}{\updefault}{\scriptsize{8}}}}}}
\put(8711,-2221){\makebox(0,0)[lb]{\smash{{\SetFigFont{12}{14.4}{\rmdefault}{\mddefault}{\updefault}{\scriptsize{10}}}}}}
\put(6821,-3166){\makebox(0,0)[lb]{\smash{{\SetFigFont{12}{14.4}{\rmdefault}{\mddefault}{\updefault}{\scriptsize{5}}}}}}
\put(6911,-3571){\makebox(0,0)[lb]{\smash{{\SetFigFont{12}{14.4}{\rmdefault}{\mddefault}{\updefault}{\small{$v_2q^{-2}$}}}}}}
\put(6956,-2671){\makebox(0,0)[lb]{\smash{{\SetFigFont{12}{14.4}{\rmdefault}{\mddefault}{\updefault}{\small{$v_2$}}}}}}
\put(7856,-2671){\makebox(0,0)[lb]{\smash{{\SetFigFont{12}{14.4}{\rmdefault}{\mddefault}{\updefault}{\small{$v_2q^2$}}}}}}
\put(8801,-2671){\makebox(0,0)[lb]{\smash{{\SetFigFont{12}{14.4}{\rmdefault}{\mddefault}{\updefault}{\small{$v_2q^4$}}}}}}
\end{picture}
\end{array}
\right)\ .
\end{equation}

\begin{prop}
{\hspace{-.2cm}.\hspace{.2cm}}
 \label{prop4}
There is a bijection between the set of standard Young $m$-tableaux of size $n$ and the set ${\mathrm{Cont}}_m(n)$.
\end{prop}

\emph{Proof.} To each standard Young $m$-tableau of size $n$ we associate a string of numbers $(a_1,\dots,a_n)$ such that $a_i$ for $i=1,\dots,n$ is the 
content of the $m$-node in which the number $i$ is placed. This string belongs to $\mathrm{Cont}_m(n)$. Indeed the condition (c1) is immediately verified. The 
condition (c2) is true: if $i$ occupies an $m$-node with a content $v_kq^{2z}$ with $z\neq0$ of a standard $m$-tableau of shape $\lambda^{(m)}$ then there is either 
an $m$-node just above in the same column or an $m$-node just to the left in the same line of the same diagram of $\lambda^{(m)}$; it carries a number less than $i$ and so its content which is 
$v_kq^{2(z+1)}$ or $v_kq^{2(z-1)}$ appears before $v_kq^{2z}$ in the string. The condition (c3) is true because if $a_i=a_j=v_kq^{2z}$ for some $i<j$ then 
due to the restrictions (\ref{sesi1})--(\ref{sesi3}) the $m$-nodes carrying $i$ and $j$ are on the same diagonal of the same tableau in the $m$-tableau. Thus the 
$m$-node just above the $m$-node carrying the number $j$ and the $m$-node just on the left of the $m$-node  carrying the number $j$ are present in $\lambda^{(m)}$; these 
$m$-nodes have contents $v_kq^{2(z+1)}$ and $v_kq^{2(z-1)}$ and are occupied by numbers $k$ and $l$ with $k,l\in\{i+1,\dots,j-1\}$ because the $m$-tableau 
is standard.

\vskip .2cm
Conversely to each string $(a_1,\dots,a_n)\in\mathrm{Cont}_m(n)$ we first associate a set of $m$-nodes of cardinality $n$. Additionally, this association distributes the numbers from $1$ to $n$ in these $m$-nodes; each $m$-node carry one number and different $m$-nodes
carry different numbers. Then we verify that the obtained $m$-tuple is a standard $m$-tableau. The construction goes as follows.

\vskip .2cm
The $m$-nodes are constructed one after another. Given a string $(a_1,\dots,a_n)\in\mathrm{Cont}_m(n)$ and assuming that $(i-1)$ $m$-nodes are already constructed, 
we add at the step number $i$ an $m$-node on the first non-occupied 
place of a diagonal whose content number is the value of $a_i$; we place the number $i$ in this $m$-node. This construction is unambiguous because the 
restrictions (\ref{sesi1})--(\ref{sesi3}) ensure that two different diagonals of such set of $m$-nodes have different content numbers. The construction of a set of $m$-nodes of total cardinality $n$ is finished.

\vskip .2cm
We shall now show that the obtained set of $m$-nodes is a standard $m$-tableau.
Assume by induction that for all $i=1,\dots,n-1$ the obtained set of $m$-nodes is a standard $m$-tableau after $i$ steps (for $i=1$ there is nothing 
to prove; the induction hypothesis is justified because for any $(a_1,\dots,a_n)\in\mathrm{Cont}_m(n)$, it is clear from the Definition \ref{def-cont} that $(a_1,\dots,a_i)\in\mathrm{Cont}_m(i)$ for all $i=1,\dots,n-1$). It is left to add an $m$-node in the position dictated by the value 
$a_n=v_kq^{2z}$, place the number $n$ in it and verify that we obtain a standard $m$-tableau.

\vskip .2cm
If, for all $i$, $0<i<n$, the number $a_i$ is different from the number $a_n$ then the $n^{\textrm{th}}$ $m$-node is added at the first place of a corresponding 
diagonal. If $z=0$ then there is nothing to prove, so assume that $z>0$ (the situation with $z<0$ is similar); the $n^{\textrm{th}}$ $m$-node is added in the first line 
and we have to prove that there is some $i$, $0<i<n$, such that $a_i=v_kq^{2(z-1)}$. Suppose that this is not the case; then by condition (c2) of the
Definition \ref{def-cont} we have some $i$, $0<i<n$, such that $a_i=v_kq^{2(z+1)}$. As $v_kq^{2z}$ is not present in the string before the $i^{\textrm{th}}$ 
position, the set of $m$-nodes constructed at the step number $i$ cannot be a standard $m$-tableau contradicting the induction hypothesis. 

\vskip .2cm
Assume that there is some $i$, $0<i<n$, such that $a_i=a_n$. We take the largest integer number $i$ satisfying this property. By construction, we add 
the $n^{\textrm{th}}$ $m$-node on the first non-occupied place of the diagonal which contains the $m$-node carrying the number $i$.
The result is a standard 
$m$-tableau only if the $m$-node just to the right of the $m$-node carrying the number $i$ and the $m$-node just below the $m$-node carrying the number $i$ are present. And 
this follows from the condition (c3) of the Definition \ref{def-cont} and the induction hypothesis.\hfill$\square$ 

\bigskip
In the example (\ref{example}) the standard Young 2-tableau is associated with the string of numbers:
\[(v_1,v_1q^2,v_2,v_1q^4,v_2q^{-2},v_1q^{-2},v_1q^{-4},v_2q^2,v_1,v_2q^4).\]

\paragraph{Remark.} The condition ``$z\neq 0$" in the part (c2) of the Definition \ref{def-cont} can be omitted for the Hecke algebras of type A; 
but this condition is necessary when $m>1$. It is transparent from the geometric point of view. For the Hecke algebra of type A, if $a_j=1$
(that is, $z=0$) for some $j>1$ then the number $j$ sits on the main diagonal but  not in the left upper corner of the standard Young tableau;
therefore (both) values $q^2$ and $q^{-2}$ are present in the string $\{ a_1,\dots,a_{j-1}\}$. However, for $m>1$, if $a_j=v_k$ for some $k$
and $j>1$, the number $j$ might occupy a  left upper corner of the standard Young $m$-tableau; in this case it clearly might happen that
none of the values $v_k q^2$ and $v_k q^{-2}$ occur in  in the string $\{ a_1,\dots,a_{j-1}\}$. 

\setcounter{equation}{0}
\section{{\hspace{-0.55cm}.\hspace{0.55cm}}Construction of representations\vspace{.25cm}}\label{ygcyclo'}

We proceed as in \cite{OP}. We first define an algebra structure on the tensor product of the algebra $H(m,1,n)$ with a free associative algebra
generated by the standard $m$-tableaux corresponding to $m$-partitions of $n$; the Baxterized elements are useful here. Then, by evaluation 
(with the help of the simplest one-dimensional representation of $H(m,1,n)$) from the right, we build representations. 

\vskip .2cm
Using the  tensor product of the algebra $H(m,1,n)$ with the free associative algebra generated by the standard $m$-tableaux we define 
and study in Appendix to this Section a 
structure of a module on the tensor product of two representations corresponding to two $m$-partitions. 

\subsection{Baxterized elements}

In the definition of the tensor product of the algebra $H(m,1,n)$ with a free associative algebra we shall frequently use the so-called
Baxterized elements.

\bigskip
Define, for any $\sigma_i$ among the generators $\sigma_1,\dots,\sigma_{n-1}$ of $H(m,1,n)$, the Baxterized elements $\sigma_i(\alpha,\beta)$ by 
\begin{equation}\label{bax-sig}
\sigma_i(\alpha,\beta):=\sigma_i+(q-q^{-1})\frac{\beta}{\alpha-\beta}\ .
\end{equation} 
The parameters $\alpha$ and $\beta$ are called spectral parameters. We recall some useful relations for the Baxterized generators $\sigma_i$. Let
\begin{equation}\label{defff}
f(\alpha,\beta)=\frac{q\alpha-q^{-1}\beta}{\alpha-\beta}\ .\end{equation}

\begin{prop}
 {\hspace{-.2cm}.\hspace{.2cm}}
 \label{prop-bax}
The following relations hold:
\begin{equation}\label{rel-bax}
\begin{array}{cc}\sigma_i(\alpha,\beta)\sigma_i(\beta,\alpha)=f(\alpha,\beta)f(\beta,\alpha),\\[0.7em]
\sigma_i(\alpha,\beta)\sigma_{i+1}(\alpha,\gamma)\sigma_i(\beta,\gamma)=\sigma_{i+1}(\beta,\gamma)\sigma_i(\alpha,\gamma)\sigma_{i+1}(\alpha,\beta),&\\[0.85em]
\sigma_i(\alpha,\beta)\sigma_j(\gamma,\delta)=\sigma_j(\gamma,\delta)\sigma_i(\alpha,\beta)&\textrm{if\ \  $|i-j|>1$.}
\end{array}
\end{equation}
\end{prop}

\emph{Proof.} It is a straightforward and well-known calculation.\hfill$\square$

\bigskip
In the construction of representations we shall often verify relations for the Baxterized elements, as in \cite{OP}. Relations will be verified
for specific values of the spectral parameters. The following Lemma shows that the original relations follow from the relations for the Baxterized elements
with fixed values of the spectral parameters.

\begin{lemm}
 {\hspace{-.2cm}.\hspace{.2cm}}
 \label{prop-bax2}
 Let $A$ and $B$ be two elements of an arbitrary associative unital algebra ${\cal{A}}$. Denote  $A(\alpha,\beta):=A+(q-q^{-1})\frac{\beta}{\alpha-\beta}$ and
 $B(\alpha,\beta):=B+(q-q^{-1})\frac{\beta}{\alpha-\beta}$ where $\alpha$ and $\beta$ are parameters.

\bigskip
(i) If 
$$A(\alpha,\beta)A(\beta,\alpha)=f(\alpha,\beta)f(\beta,\alpha)\ ,$$
where $f$ is defined in (\ref{defff}), for some (arbitrarily) fixed values of the parameters $\alpha$ and $\beta$ ($\alpha\neq\beta$) then 
$$A^2-(q-q^{-1})A-1=0\ .$$ 

\medskip
(ii) If 
$$A^2-(q-q^{-1})A-1=0\ ,\ B^2-(q-q^{-1})B-1=0$$ 
and
$$A(\alpha,\beta)B(\alpha,\gamma)A(\beta,\gamma)=B(\beta,\gamma)A(\alpha,\gamma)B(\alpha,\beta)$$
for some (arbitrarily) fixed values of the parameters $\alpha$,$\beta$ and $\gamma$ ($\alpha\neq\beta\neq\gamma\neq\alpha$) then
$$ABA=BAB\ .$$

\medskip
(iii) If 
$$A(\alpha,\beta)B(\gamma,\delta)=B(\gamma,\delta)A(\alpha,\beta)$$
for some (arbitrarily) fixed values of the parameters $\alpha$,$\beta$,$\gamma$ and $\delta$ ($\alpha\neq\beta$ and $\gamma\neq\delta$) then
$$AB=BA\ .$$
\end{lemm}

\emph{Proof.} (i) We have 
\[\begin{array}{ll}
 & A(\alpha,\beta)A(\beta,\alpha)=f(\alpha,\beta)f(\beta,\alpha) \\[1em]
\Rightarrow & A^2+(q-q^{-1})\left(\displaystyle{\frac{\beta}{\alpha-\beta}}+\displaystyle{\frac{\alpha}{\beta-\alpha}}\right) A
-(q-q^{-1})^2\displaystyle{\frac{\alpha\beta}{(\beta-\alpha)^2}}=\displaystyle{\frac{\alpha^2+\beta^2-\alpha\beta(q^2+q^{-2})}{(\beta-\alpha)^2}} \\[1.5em]
\Rightarrow & A^2-(q-q^{-1})A-1=0\ .
\end{array}\]
 
(ii) We have
\[\begin{array}{ll}
 &A(\alpha,\beta)B(\alpha,\gamma)A(\beta,\gamma)-B(\beta,\gamma)A(\alpha,\gamma)B(\alpha,\beta)=0\\[1em]
\Rightarrow & ABA-BAB+(q-q^{-1})(A^2-B^2)\displaystyle{\frac{\gamma}{\alpha-\gamma}}+\\[1em]
&+(q-q^{-1})^2(A-B)\biggl(\displaystyle{\frac{\gamma}{\alpha-\gamma}}(\displaystyle{\frac{\beta}{\alpha-\beta}}+\displaystyle{\frac{\gamma}{\beta-\gamma}})
-\displaystyle{\frac{\gamma}{\beta-\gamma}}\displaystyle{\frac{\beta}{\alpha-\beta}}\biggr)=0\\[1.5em]
\Rightarrow & ABA-BAB+(q-q^{-1})^2(A-B)\displaystyle{\frac{\gamma}{\alpha-\gamma}}\biggl(1+\displaystyle{\frac{\beta}{\alpha-\beta}}
+\displaystyle{\frac{\gamma}{\beta-\gamma}}-\displaystyle{\frac{\beta(\alpha-\gamma)}{(\beta-\gamma)(\alpha-\beta)}}\biggr)=0\\[1.5em]
\Rightarrow & ABA-BAB=0\ .
\end{array}\]

(iii) It is immediate that $A(\alpha,\beta)B(\gamma,\delta)=B(\gamma,\delta)A(\alpha,\beta)$ implies $AB=BA$.
\hfill$\square$

\subsection{Smash product}\label{subsection.smash.product}

We pass to the definition of the associative algebra structure on the product of the algebra $H(m,1,n)$ with a free associative algebra whose 
generators are indexed by the standard $m$-tableaux corresponding to $m$-partitions of $n$. The resulting algebra we shall denote by ${\mathfrak{T}}$.

\bigskip
Let $\lambda^{(m)}$ be an $m$-partition of size $n$. Consider a set of free generators labeled by standard $m$-tableaux of shape 
$\lambda^{(m)}$; for a standard $m$-tableau $X_{\lambda^{(m)}}$ we denote by $\mathcal{X}_{\lambda^{(m)}}$ the corresponding free generator 
and by $c(X_{\lambda^{(m)}}|i)$ the content  (see the preceding Section for the definition) of the $m$-node carrying the number $i$.

\vskip .2cm
In the sequel we shall use the Artin generators of the symmetric group.
Recall that the symmetric group (whose group ring is the classical limit of the A-type Hecke algebra $H_n$) in the Artin presentation is generated
by the elements $s_i$, $i=1,...,n-1$, with the defining relations
\begin{empheq}[left=\empheqlbrace]{alignat=4}
\label{def1sya}&s_i s_{i+1}s_i=s_{i+1}s_i s_{i+1} &\hspace{1cm}& \textrm{for all $i=1,\dots,n-2$\ ,}\\[.5em]
\label{def1syb}&s_i s_j=s_j s_i && \textrm{for all $i,j=1,\dots,n-1$ such that $|i-j|>1$\ ,}\\[.5em]
\label{def1syc}&s_i^2=1&\hspace{1cm}& \textrm{for all $i=1,\dots,n-1$\ .}
\end{empheq}

\vskip .2cm
Let $X_{\lambda^{(m)}}$ be an $m$-partition of $n$ whose $m$-nodes are filled with numbers from 1 to $n$; thus different $m$-nodes of the $m$-tableau 
$X_{\lambda^{(m)}}$ carry different numbers.
The $m$-tableau $X_{\lambda^{(m)}}$ is not necessarily standard. By definition, for such $m$-tableau $X_{\lambda^{(m)}}$ and any permutation 
$\pi\in S_n$, the $m$-tableau $\xl{\lam}{\pi}$ is obtained from the $m$-tableau $X_{\lambda^{(m)}}$ by applying the permutation $\pi$ to the 
numbers occupying the $m$-nodes of $X_{\lambda^{(m)}}$; for example $\xl{\lam}{s_i}$ is the $m$-tableau 
obtained from $X_{\lambda^{(m)}}$ by exchanging the positions of the numbers $i$ and $(i+1)$ in the $m$-tableau $X_{\lambda^{(m)}}$. We thus 
have by construction:
\begin{equation}\label{contents}c(\xl{\lam}{\pi}|i)=c(X_{\lambda^{(m)}}|\pi^{-1}(i))\end{equation}
for all $i=1,\dots,n$.

\vskip .2cm
{}For a standard $m$-tableau $X_{\lambda^{(m)}}$, the $m$-tableau $\xl{\lam}{\pi}$ is not necessarily standard. 
As for the generators of the free algebra, we denote the generator corresponding to the $m$-tableau $\xl{\lam}{\pi}$ by 
$\ml{\lam}{\pi}$ if the $m$-tableau $\xl{\lam}{\pi}$ is standard. And if the $m$-tableau $\xl{\lam}{\pi}$ is not standard 
then we put $\ml{\lam}{\pi}=0$. 

\begin{prop}
{\hspace{-.2cm}.\hspace{.2cm}}
 \label{prop-rel}
The relations
\begin{equation}\label{rel-a1}
\Bigl(\sigma_i+\frac{(q-q^{-1})
c(X_{\lambda^{(m)}}|i+1)}{c(X_{\lambda^{(m)}}|i)-c(X_{\lambda^{(m)}}|i+1)}\Bigr)\cdot\mathcal{X}_{\lambda^{(m)}}=\ml{\lam}{s_i}\cdot\Bigl(\sigma_i+\frac{(q-q^{-1})c(X_{\lambda^{(m)}}|i)}{c(X_{\lambda^{(m)}}|i+1)-c(X_{\lambda^{(m)}}|i)}\Bigr)\end{equation}
and
\begin{equation}\label{rel-a2}
\Bigl(\tau-c(X_{\lambda^{(m)}}|1)\Bigr)\cdot\mathcal{X}_{\lambda^{(m)}}=0
\end{equation}
are compatible with the relations for the generators $\tau,\sigma_1,\dots,\sigma_{n-1}$ of the algebra $H(m,1,n)$. 
\end{prop}

Before the proof we explain the meaning of the word ``compatible" in the formulation of the Proposition. 

\vskip .2cm
Let ${\cal{F}}$ be the free associative algebra
generated by $\tilde{\tau}$, $\tsigma_1,\dots,\tsigma_{n-1}$.  The algebra $H(m,1,n)$ is naturally the quotient of ${\cal{F}}$. 

\vskip .2cm
Let $\mathbb{C}[\mathcal{X}]$ be 
the free associative algebra whose generators  
$\mathcal{X}_{\lambda^{(m)}}$  range over all standard $m$-tableaux of shape $\lambda^{(m)}$ for all $m$-partitions $\lambda^{(m)}$ of $n$.

\vskip .2cm
Consider an algebra structure on the space $\mathbb{C}[\mathcal{X}]\otimes {\cal{F}}$ for which: (i) the map $\iota_1:x\mapsto x\otimes 1$, 
$x\in \mathbb{C}[\mathcal{X}]$, is an isomorphism of $\mathbb{C}[\mathcal{X}]$ with its image with respect to $\iota_1$; (ii) the map 
$\iota_2:\phi\mapsto 1\otimes \phi$, $\phi\in {\cal{F}}$, is an isomorphism of ${\cal{F}}$ with its image with respect to $\iota_2$; (iii) 
the formulas (\ref{rel-a1})-(\ref{rel-a2}), extended by associativity, provide the rules to rewrite elements of the form $(1\otimes \phi)(x\otimes 1)$, $x\in \mathbb{C}[\mathcal{X}]$,
$\phi\in {\cal{F}}$, as elements of $\mathbb{C}[\mathcal{X}]\otimes {\cal{F}}$. 

\vskip .2cm
The ``compatibility" means that we have an induced structure of an associative algebra, denoted by ${\mathfrak{T}}$,
on the space $\mathbb{C}[\mathcal{X}]\otimes H(m,1,n)$. 
More precisely, if we multiply any relation of the cyclotomic Hecke algebra $H(m,1,n)$ (the relation is viewed as an element of the free algebra 
${\cal{F}}$) from the right by a generator $\mathcal{X}_{\lambda^{(m)}}$ (this is a combination of the form ``a relation of $H(m,1,n)$ times 
$\mathcal{X}_{\lambda^{(m)}}$") and use the ``instructions" (\ref{rel-a1})-(\ref{rel-a2}) to move all appearing $\mathcal{X}$'s to the left (the free generator 
changes but the expression stays always linear in $\mathcal{X}$) then we obtain a linear combination of terms of the 
form ``$\ml{\lam}{\pi}$, $\pi\in S_n$, times a relation of $H(m,1,n)$". 
 
\vskip .2cm
\emph{Proof.} We rewrite the relation (\ref{rel-a1}) using the Baxterized form of the elements $\sigma_i$:
\[\sigma_i\Bigl(c(X_{\lambda^{(m)}}|i),c(X_{\lambda^{(m)}}|i+1)\Bigr)\cdot\mathcal{X}_{\lambda^{(m)}}=\ml{\lam}{s_i}\cdot\sigma_i\Bigl(c(X_{\lambda^{(m)}}|i+1),c(X_{\lambda^{(m)}}|i)\Bigr)\ .\]
For brevity we denote $c^{(k)}:=c(X_{\lambda^{(m)}}|k)$ for all $k=1,\dots,n$. 

\vskip .2cm
We shall check the compatibility of the ``instructions" (\ref{rel-a1})-(\ref{rel-a2}) with the set of defining relations (\ref{def1a})--(\ref{def1'b}) and 
(\ref{def1''b})--(\ref{def1''a}). 
We start with the relations involving the generators $\sigma_i$ only. Here we use the Baxterized form of the relations and the Lemma \ref{prop-bax2}.

\vskip .2cm
Below we shall use without mentioning the inequalities $c^{(k)}\neq c^{(k+1)}$, $c^{(k)}\neq c^{(k+2)}$, $c^{(k+1)}\neq c^{(k+2)}$ (for all $k$) 
which follow from the restrictions (\ref{sesi1})--(\ref{sesi3}). 

\vskip .2cm
(a) We consider first the relation $\sigma_i^2=(q-q^{-1})\sigma_i+1$. 
 If the $m$-tableau $\xl{\lam}{s_i}$ is standard then we analyze this
relation in its equivalent form, given in (i) in the Lemma \ref{prop-bax2}. We have
\begin{equation}\label{trpr0}
\begin{array}{ll}&\sigma_i(c^{(i+1)},c^{(i)})\ \sigma_i(c^{(i)},c^{(i+1)})\cdot\mathcal{X}_{\lambda^{(m)}}\\[1em]
=&\sigma_i(c^{(i+1)},c^{(i)})\cdot \ml{\lam}{s_i} \cdot\sigma_i(c^{(i+1)},c^{(i)})\\[1em]
=&\mathcal{X}_{\lambda^{(m)}}\cdot\sigma_i(c^{(i)},c^{(i+1)})\ \sigma_i(c^{(i+1)},c^{(i)})\ .\end{array}\end{equation}

We used (\ref{contents}) in the second equality. Therefore,
\[\begin{array}{ll}&\Bigl(\sigma_i(c^{(i+1)},c^{(i)})\ \sigma_i(c^{(i)},c^{(i+1)})-f(c^{(i+1)},c^{(i)})f(c^{(i)},c^{(i+1)})\Bigr)\cdot\mathcal{X}_{\lambda^{(m)}}\\[1em]
=&\mathcal{X}_{\lambda^{(m)}}\cdot\Bigl(\sigma_i(c^{(i)},c^{(i+1)})\ \sigma_i(c^{(i+1)},c^{(i)})-f(c^{(i+1)},c^{(i)})f(c^{(i)},c^{(i+1)})\Bigr)\end{array}
\]
and the compatibility for this relation is verified since the expression on the right of $\mathcal{X}_{\lambda^{(m)}}$ belongs, by the Proposition \ref{prop-bax}, 
to the ideal generated by relations.

\vskip .2cm
If the $m$-tableau $\xl{\lam}{s_i}$ is not standard then $(i+1)$ sits in the same tableau of the $m$-tableau $X_{\lambda^{(m)}}$ as $i$ and is 
situated directly to the right or directly down with respect to $i$. In this situation the relation (\ref{rel-a1}) can be rewritten as $\sigma_i
\mathcal{X}_{\lambda^{(m)}}=w\mathcal{X}_{\lambda^{(m)}}$, where $w$ is equal to either $q$ or $(-q^{-1})$, and the verification of the compatibility
of the relation $\sigma_i^2-(q-q^{-1})\sigma_i-1=0$ with the instructions (\ref{rel-a1})-(\ref{rel-a2}) is straightforward.

\vskip .2cm
(b) The relation $\sigma_i\sigma_{i+1}\sigma_i=\sigma_{i+1}\sigma_i\sigma_{i+1}$ we analyze in its equivalent form, given in (ii) in the 
Lemma \ref{prop-bax2}. We have
\begin{equation}\label{trpr1}
\begin{array}{ll}&\sigma_i(c^{(i+1)},c^{(i+2)})\sigma_{i+1}(c^{(i)},c^{(i+2)})\sigma_i(c^{(i)},c^{(i+1)})\cdot\mathcal{X}_{\lambda^{(m)}}\\[0.5em]
=&\ml{\lam}{s_is_{i+1}s_i}\cdot\sigma_i(c^{(i+2)},c^{(i+1)})\sigma_{i+1}(c^{(i+2)},c^{(i)})\sigma_i(c^{(i+1)},c^{(i)})\end{array}\end{equation}
and
\begin{equation}\label{trpr2}
\begin{array}{ll}&\sigma_{i+1}(c^{(i)},c^{(i+1)})\sigma_i(c^{(i)},c^{(i+2)})\sigma_{i+1}(c^{(i+1)},c^{(i+2)})\cdot\mathcal{X}_{\lambda^{(m)}}\\[0.5em]
=&\ml{\lam}{s_{i+1}s_is_{i+1}}\cdot(\sigma_{i+1}(c^{(i+1)},c^{(i)})\sigma_i(c^{(i+2)},c^{(i)})\sigma_{i+1}(c^{(i+2)},c^{(i+1)})\  .
\end{array}\end{equation}
We used several times the relations (\ref{contents}). 

\vskip .2cm
One might think that, as for the relation $\sigma_i^2-(q-q^{-1})\sigma_i-1=0$, we should consider separately the situations when 
in the process of transformation the $m$-tableau becomes non-standard. However one verifies that for an arbitrary standard $m$-tableau
$Y_{\lambda^{(m)}}$:

\begin{itemize}
\item if the $m$-tableau $\yl{\lam}{s_i}$ is not standard then the $m$-tableaux $\yl{\lam}{s_is_{i+1}}$ and $\yl{\lam}{s_is_{i+1}s_i}$ 
are not standard as well;
\item if the $m$-tableau $\yl{\lam}{s_{i+1}}$ is not standard then the $m$-tableaux $\yl{\lam}{s_{i+1}s_i}$ and 
$\yl{\lam}{s_{i+1}s_is_{i+1}}$ are not standard as well.
\end{itemize}

It then follows that

\begin{itemize}
\item if the $m$-tableau $\yl{\lam}{s_is_{i+1}s_{i}}$ is standard then the $m$-tableaux $\yl{\lam}{s_i}$ and $\yl{\lam}{s_is_{i+1}}$ are 
standard as well;
\item if the $m$-tableau $\yl{\lam}{s_{i+1}s_{i}s_{i+1}}$ is standard then the $m$-tableaux $\yl{\lam}{s_{i+1}}$ and $\yl{\lam}{s_{i+1}s_i}$ 
are standard as well.
\end{itemize}

Therefore, we cannot return to a standard $m$-tableau if an intermediate $m$-tableau was not standard. Thus the equalities (\ref{trpr1}) and (\ref{trpr2}) are 
always valid, in contrast to (\ref{trpr0}). 

\vskip .2cm
We replace $\ml{\lam}{s_{i+1}s_is_{i+1}}$ by $\ml{\lam}{s_is_{i+1}s_i}$ in the right hand side of  (\ref{trpr2}) and 
subtract  (\ref{trpr2}) from  (\ref{trpr1}). In the result, the expression on the right of $\ml{\lam}{s_is_{i+1}s_i}$ belongs, by the Proposition 
\ref{prop-bax},  to the ideal generated by relations.

\vskip .2cm
(c) The relation $\sigma_i\sigma_{j}=\sigma_{j}\sigma_i$ for $\vert i-j\vert >1$ we analyze in its equivalent form, given in (iii) in the 
Lemma \ref{prop-bax2}. We have
\begin{equation}\label{trpr3}
\begin{array}{ll}&\Bigl(\sigma_i(c^{(i)},c^{(i+1)})\sigma_{j}(c^{(j)},c^{(j+1)})-\sigma_{j}(c^{(j)},c^{(j+1)})
\sigma_i(c^{(i)},c^{(i+1)})\Bigr)\cdot\mathcal{X}_{\lambda^{(m)}}\\[.5em]
=&\ml{\lam}{s_is_j}\cdot\Bigl(\sigma_i(c^{(i+1)},c^{(i)})\sigma_{j}(c^{(j+1)},c^{(j)})-
\sigma_{j}(c^{(j+1)},c^{(j)})\sigma_i(c^{(i+1)},c^{(i)})\Bigr)\ .
\end{array}\end{equation}
Again, as for the previous relation, a direct inspection shows that (\ref{trpr3}) is always valid.

\vskip .2cm
The expression on the right of $\ml{\lam}{s_is_j}$ in the right hand side of (\ref{trpr3}) belongs again to the ideal generated by relations 
by the Proposition \ref{prop-bax}.

\vskip .2cm
(d) It is left to analyze the relations including the generator $\tau$. 

\vskip .2cm
The verification of the compatibility of the relations $(\tau-v_1)\dots(\tau-v_m)=0$ and 
$\tau\sigma_i=\sigma_i\tau$ for $i>1$ with the instructions (\ref{rel-a1})-(\ref{rel-a2}) is immediate. 

\vskip .2cm
The compatibility of the remaining relation $\tau\sigma_1\tau\sigma_1=\sigma_1\tau\sigma_1\tau$ with the instructions (\ref{rel-a1})-(\ref{rel-a2}) is a direct 
consequence of the Lemma below.\hfill$\square$

\begin{lemm}
{\hspace{-.2cm}.\hspace{.2cm}}
 \label{prop-rel2}
The relations (\ref{rel-a1})-(\ref{rel-a2}) imply the relations:
\begin{equation}\label{rel2}
\bigl(J_i-c(X_{\lambda^{(m)}}|i)\bigr)\cdot\mathcal{X}_{\lambda^{(m)}}=0\quad\textrm{for all $i=1,\dots,n$\ .}
\end{equation}
\end{lemm}

\emph{Proof.} For brevity we denote $c^{(k)}:=c(X_{\lambda^{(m)}}|k)$ for all $k=1,\dots,n$. 

\vskip .2cm
Recall that $J_1=\tau$ and $J_{i+1}=\sigma_iJ_i\sigma_i$. We use induction on $i$. For $i=1$ the relation (\ref{rel2}) is the relation (\ref{rel-a2}). 

\vskip .2cm
Assume first that the $m$-tableau $\xl{\lam}{s_i}$ is standard. Then 
\[\begin{array}{ll}
\sigma_iJ_i\sigma_i\cdot\mathcal{X}_{\lambda^{(m)}}\!\!\!\!&=\sigma_iJ_i\cdot\Bigl(-(q-q^{-1})\frac{c^{(i+1)}}{c^{(i)}-c^{(i+1)}}\mathcal{X}_{\lambda^{(m)}}+\ml{\lam}{s_i}\cdot\sigma_i(c^{(i+1)},c^{(i)})\Bigr)\\[0.8em]
&=\sigma_i\cdot\Bigl(-(q-q^{-1})c^{(i)}\frac{c^{(i+1)}}{c^{(i)}-c^{(i+1)}}\mathcal{X}_{\lambda^{(m)}}+c^{(i+1)}\ml{\lam}{s_i}\cdot\sigma_i(c^{(i+1)},c^{(i)})\Bigr)\\[0.8em]
&=(q-q^{-1})^2c^{(i)}\frac{c^{(i+1)}c^{(i+1)}}{(c^{(i+1)}-c^{(i)})^2}\mathcal{X}_{\lambda^{(m)}}-(q-q^{-1})\frac{c^{(i)}c^{(i+1)}}{c^{(i)}-c^{(i+1)}}\ml{\lam}{s_i}\cdot\sigma_i(c^{(i+1)},c^{(i)})\\[0.6em]
&\ -(q-q^{-1})\frac{c^{(i+1)}c^{(i)}}{c^{(i+1)}-c^{(i)}}\ml{\lam}{s_i}\cdot\sigma_i(c^{(i+1)},c^{(i)})+ c^{(i+1)}\mathcal{X}_{\lambda^{(m)}}\cdot\sigma_i(c^{(i)},c^{(i+1)})\sigma_i(c^{(i+1)},c^{(i)})\\[0.8em]
&=c^{(i+1)}\Bigl((q-q^{-1})^2\frac{c^{(i)}c^{(i+1)}}{(c^{(i+1)}-c^{(i)})^2}+\frac{c^{(i)}c^{(i)}+c^{(i+1)}c^{(i+1)}-c^{(i)}c^{(i+1)}(q^2+q^{-2})}{(c^{(i+1)}-c^{(i)})^2}\Bigr)\mathcal{X}_{\lambda^{(m)}}\\[0.8em]
&=c^{(i+1)}\mathcal{X}_{\lambda^{(m)}}\ .
\end{array}\]
Here we moved the elements $\sigma_i$ to the right, using the relation (\ref{rel-a1});  we used then relations (\ref{contents}), the induction hypothesis and
the first relation in the Proposition \ref{prop-bax}. 

\vskip .2cm
Then assume that the $m$-tableau $\xl{\lam}{s_i}$ is not standard. It means that the $m$-nodes  carrying numbers $i$ and $(i+1)$ are adjacent 
(neighbors in a same line or a same column of a tableau in the $m$-tableau). In this situation we have $c^{(i)}=q^{\pm2}c^{(i+1)}$ and 
\[\begin{array}{ll}
\sigma_iJ_i\sigma_i\cdot\mathcal{X}_{\lambda^{(m)}}&=\sigma_iJ_i\cdot\Bigl(-(q-q^{-1})\frac{c^{(i+1)}}{q^{\pm2}c^{(i+1)}-c^{(i+1)}}\Bigr)\mathcal{X}_{\lambda^{(m)}}\\[0.6em]
&=\sigma_i\cdot\Bigl(-(q-q^{-1})\frac{q^{\pm2}c^{(i+1)}c^{(i+1)}}{q^{\pm2}c^{(i+1)}-c^{(i+1)}}\Bigr)\mathcal{X}_{\lambda^{(m)}}\\[0.6em]
&=(q-q^{-1})^2\ \frac{q^{\pm2}c^{(i+1)}c^{(i+1)}c^{(i+1)}}{(q^{\pm2}c^{(i+1)}-c^{(i+1)})^2}\mathcal{X}_{\lambda^{(m)}}
\\[0.6em]
&=c^{(i+1)}\mathcal{X}_{\lambda^{(m)}}\ .
\end{array}\]
Here we moved the elements $\sigma_i$ to the right using the relation (\ref{rel-a1});  we used then the induction hypothesis.
\hfill$\square$

\subsection{Representations}\label{suserep}

The Proposition \ref{prop-rel} provides an effective tool for the construction of representations of the cyclotomic Hecke algebra $H(m,1,n)$.

\bigskip
Let $\vert\rangle$ be a ``vacuum" - a basic vector of a one-dimensional $H(m,1,n)$-module; for example, 
\begin{equation}\sigma_i\vert\rangle=q\vert\rangle\ \ \text{for all $i$\ and}\ \  
\tau\vert\rangle =v_1\vert\rangle\ .\label{movac}\end{equation} 
Moving, in the expressions $\phi\mathcal{X}_{\lambda^{(m)}}\vert\rangle$, $\phi\in H(m,1,n)$, the elements 
$\mathcal{X}$'s to the left and using the module structure (\ref{movac}), we build, due to the compatibility, a representation of $H(m,1,n)$ on the vector space $U_{\lambda^{(m)}}$ with the basis
$\mathcal{X}_{\lambda^{(m)}}\vert\rangle$. We shall, by a slight abuse of notation, denote the symbol $\mathcal{X}_{\lambda^{(m)}}\vert\rangle$ again
by $\mathcal{X}_{\lambda^{(m)}}$. This procedure leads to the following formulas for the action of the generators $\tau,\sigma_1,\dots,\sigma_{n-1}$ on 
the basis vectors $\mathcal{X}_{\lambda^{(m)}}$ of $U_{\lambda^{(m)}}$:
\begin{equation}\label{rep-a1}\sigma_i\,:\,\mathcal{X}_{\lambda^{(m)}}\mapsto\ 
{\displaystyle \frac{(q-q^{-1})c(X_{\lambda^{(m)}}|i+1)}{c(X_{\lambda^{(m)}}|i+1)-c(X_{\lambda^{(m)}}|i)}}\mathcal{X}_{\lambda^{(m)}}
+\ {\displaystyle \frac{qc(X_{\lambda^{(m)}}|i+1)-q^{-1}c(X_{\lambda^{(m)}}|i)}{c(X_{\lambda^{(m)}}|i+1)-
c(X_{\lambda^{(m)}}|i)}}\ml{\lam}{s_i}\end{equation}
and
\begin{equation}\label{rep-a2}
\tau\,:\,\mathcal{X}_{\lambda^{(m)}}\mapsto c(X_{\lambda^{(m)}}|1)\mathcal{X}_{\lambda^{(m)}}\ .
\end{equation}
As before, it is assumed here that $\ml{\lam}{s_i}=0$ if $\xl{\lam}{s_i}$ is not a standard $m$-tableau.
Denote this $H(m,1,n)$-module by $V_{\lambda^{(m)}}$.

\vskip .2cm
{}Assume that the $m$-tableau $\xl{\lam}{s_i}$ is standard. The two-dimensional subspace of $U_{\lambda^{(m)}}$ with the 
basis  $\{ \mathcal{X}_{\lambda^{(m)}},\ml{\lam}{s_i}\}$ is $\sigma_i$-invariant. For the future convenience we write down the 
matrix giving the action of the generator $\sigma_i$ on this two-dimensional subspace:
\begin{equation}\label{rep-sigma-two}{\displaystyle \frac{1}{c^{(i+1)}-c^{(i)}}}
\left(\begin{array}{lr} (q-q^{-1})c^{(i+1)}&q^{-1} c^{(i+1)}-qc^{(i)}\\[1em]
qc^{(i+1)}-q^{-1}c^{(i)}&-(q-q^{-1})c^{(i)}\end{array}\right)\ ,\end{equation}
where we denoted $c^{(i)}=c(X_{\lambda^{(m)}}|i)$ and $c^{(i+1)}=c(X_{\lambda^{(m)}}|i+1)$. 
 
\paragraph{Remarks.} $\ $ 

\vskip .2cm
{\bf (a)} In our construction of representations, the Baxterization of the generator $\tau$ (which can be found in \cite{IO2}) is never used while the 
Baxterized generators $\sigma_i$, $i=1,\dots,n-1$, appear. The relation (\ref{rel-a2}) says that $\tau$, placed before $\mathcal{X}_{\lambda^{(m)}}$,
can be immediately replaced by a number.  This is similar to the situation with the element $\sigma_1$ in the representation theory of the Hecke  ($m=1$)
algebras. Indeed
if $m=1$ then for any standard tableau $X_{\lambda}$, the tableau $X_{\lambda}^{s_1}$ is non-standard and so $\sigma_1$, placed before  
$\mathcal{X}_{\lambda}$, can be immediately replaced by a number; in particular, the action of $\sigma_1$ is given by a diagonal matrix
in the basis $\mathcal{X}_{\lambda}\vert\rangle$ of $U_{\lambda}$; the behavior of $\tau$ extends this phenomenon to $m>1$.

\vskip .3cm 
{\bf (b)} It follows from the preceding remark that the action of the generators in the constructed representations do not depend on the value of $\tau$ 
on the vacuum $\vert\rangle$. 
Moreover the constructed representations do not depend (up to isomorphism) on the value of $\sigma_i$, $i=1,\dots,n-1$, on the vacuum. Indeed if we take for the vacuum a basic vector $\vert\rangle'$ of the one-dimensional $H(m,1,n)$-module such that $\sigma_i\vert\rangle'=-q^{-1}\vert\rangle'$ for all $i$ and
$\tau\vert\rangle' =v_1\vert\rangle'$, the procedure described in this subsection leads to representations $\tilde{V}_{\lambda^{(m)}}$ of $H(m,1,n)$.
By construction, $V_{\lambda^{(m)}}$ and $\tilde{V}_{\lambda^{(m)}}$ have the same underlying vector space $U_{\lambda^{(m)}}$. 
It is straightforward to check that for any $m$-partition $\lambda^{(m)}$ there is an isomorphism between $H(m,1,n)$-modules 
$V_{\lambda^{(m)}}$ and $\tilde{V}_{\lambda^{(m)}}$; the operators for the representation $\tilde{V}_{\lambda^{(m)}}$ are obtained from the operators for the representation $V_{\lambda^{(m)}}$ by the following diagonal change of basis in $U_{\lambda^{(m)}}$:
\[\mathcal{X}_{\lambda^{(m)}}\mapsto \mathfrak{c}_{_{\scriptstyle\mathcal{X}_{\lambda^{(m)}}}}\mathcal{X}_{\lambda^{(m)}}\ ,\ \ \text{where}\ \  
\mathfrak{c}_{_{\scriptstyle\mathcal{X}_{\lambda^{(m)}}}}
=\prod\limits_{i\, :\, \ml{\lam}{s_i}\neq0}\ \left( q\ c(X_{\lambda^{(m)}}|i)-q^{-1} c(X_{\lambda^{(m)}}|i+1)\right)\ .\] 
By construction, $\mathfrak{c}_{\mathcal{X}_{\lambda^{(m)}}}\neq 0$.

\vskip .3cm
{\bf (c)} In the Hecke situation ($m=1$) the coefficients appearing in the action of the generators can be expressed in terms of the lengths $l_{j,j+1}$ between  
nodes (see, e.g., \cite{OP}). We do not define the length between nodes which do not belong to the same tableau of the $m$-tableau;
the form, referring to lengths, of the action is not useful any more. 

\vskip .3cm
{\bf (d)} The constructed action of the generators in the representations $V_{\lambda^{(m)}}$ coincides with the action given in \cite{AK}.

\vskip .3cm
{\bf (e)} The action of the intertwining operators $U_{i+1}=\sigma_iJ_i-J_i\sigma_i$, $i=1,\dots,n-1$, (see paragraph \textbf{5} of Section \ref{ygcyclo}) in a representation $V_{\lambda^{(m)}}$ is:
\begin{equation}\label{act-int}
U_{i+1}(\mathcal{X}_{\lambda^{(m)}})=\left(q^{-1}c^{(i)}-qc^{(i+1)}\right)\mathcal{X}_{\lambda^{(m)}}^{s_i}\ ,
\end{equation}
where $c^{(i)}=c(X_{\lambda^{(m)}}|i)$, $i=1,\dots,n$. Indeed we rewrite $U_{i+1}=\sigma_iJ_i-\sigma_i^{-1}J_{i+1}=\sigma_i(J_i-J_{i+1})+(q-q^{-1})J_{i+1}$, so, by the Lemma \ref{prop-rel2},
\[U_{i+1}(\mathcal{X}_{\lambda^{(m)}})=(c^{(i)}-c^{(i+1)})\left(\sigma_i(\mathcal{X}_{\lambda^{(m)}})+\frac{(q-q^{-1})c^{(i+1)}}{c^{(i)}-c^{(i+1)}}\mathcal{X}_{\lambda^{(m)}}\right)\ .\]
Using (\ref{rep-a1}) we obtain the formula (\ref{act-int}).

\subsection{Scalar product\vspace{.1cm}}\label{scalarproduct} 

The representations constructed on the spaces $U_{\lambda^{(m)}}$, where $\lambda^{(m)}$ is an $m$-partition of size $n$, are analogues for $H(m,1,n)$ of the seminormal representations of the symmetric group. We compute here analogues, 
for the representations spaces of $H(m,1,n)$, of invariant scalar products on representation spaces of the symmetric group.

\paragraph{$\mathsf{S}$-invariance.} Let $\mathcal{A}$ be an associative algebra, $U$ an $\mathcal{A}$-module, $\mathcal{O}$ an abelian group and 
$\langle ,\rangle :U\times U\to\mathcal{O}$ a bi-additive map. Let $\mathsf{S}$ be a set of generators of $\mathcal{A}$.

\begin{defi}{\hspace{-.2cm}.\hspace{.2cm}}
We say that $\langle ,\rangle$ is $\mathsf{S}$-invariant if, for all $u,v\in U$ and for all $\mathsf{a}\in\mathsf{S}$,
\begin{equation}\label{scaprose2}
\langle \mathsf{a}(u),\mathsf{a}(v)\rangle =\langle u,v\rangle \ .\end{equation}
We say that $\langle ,\rangle$ is $\mathsf{S}^-$-invariant if, for all $u,v\in U$ and for all $\mathsf{a}\in\mathsf{S}$,
\begin{equation}\label{scaprose}
\langle \mathsf{a}(u),\mathsf{a}^{-1}(v)\rangle =\langle u,v\rangle \ .\end{equation}
\end{defi}
{}For an $\mathsf{S}$-invariant $\langle ,\rangle$ we have, for all $u,v\in U$,
\[\langle \mathsf{a}_1\dots\mathsf{a}_k(u),\,\mathsf{a}_1\dots\mathsf{a}_k(v)\rangle =\langle u,v\rangle\ \ \ \text{for}\ \ k\geq0\ \ \text{and}\ \ \mathsf{a}_1,\dots,\mathsf{a}_k\in\{\tau,\sigma_1,\dots ,\sigma_{n-1}\}\ .\]
For an $\mathsf{S}^-$-invariant $\langle ,\rangle$ we have, for all $u,v\in U$,
\[\langle \mathsf{a}_1\dots\mathsf{a}_k(u),\,\mathsf{a}_1^{-1}\!\!\dots\mathsf{a}_k^{-1}(v)\rangle =\langle u,v\rangle\ \ \ \text{for}\ \ k\geq0\ \ \text{and}\ \ \mathsf{a}_1,\dots,\mathsf{a}_k\in\{\tau,\sigma_1,\dots ,\sigma_{n-1}\}\ .\]

\vskip .2cm
\emph{Note.} In the sequel, $\mathcal{A}=H(m,1,n)$ and we always choose $\mathsf{S}=\{\tau,\sigma_1,\dots ,\sigma_{n-1}\}$.

\paragraph{The ring $\mathfrak{R}$.} 
Let $$\mathfrak{D}:=\{ 1+q^2+\dots+q^{2N}\}_{N=1,\dots,n}\cup \{ q^{2i}v_j-v_k\}_{i,j,k\colon  j\neq k;\, -n<i<n}$$
and let $\mathfrak{R}$ be the ring $\mathbb{C}[q,q^{-1},v_1,v_1^{-1},\dots,v_m,v_m^{-1}]$ of Laurent polynomials in variables $q,v_1,\dots,v_m$ localized
with respect to the multiplicative set generated by $\mathfrak{D}$. 

\vskip .2cm
In this Subsection we shall work with the generic cyclotomic Hecke algebra which means here the cyclotomic Hecke algebra over the ring $\mathfrak{R}$.
We denote by the same symbol $U_{\lambda^{(m)}}$ the free module, now over the ring $\mathfrak{R}$, with the basis $\mathfrak{B}:=\{\mathcal{X}_{\lambda^{(m)}}\}$. As before, the generators $\tau,\sigma_1,\dots,\sigma_{n-1}$ of the algebra $H(m,1,n)$ act on  $U_{\lambda^{(m)}}$ according to the 
formulas (\ref{rep-a1})--(\ref{rep-a2}).

\vskip .2cm
We shall define a map $\mathfrak{B}\times\mathfrak{B}\to\mathfrak{R}$ which we denote by the same symbol $\langle ,\rangle$ as above
because it will be extended to several bi-additive maps $U_{\lambda^{(m)}}\times U_{\lambda^{(m)}}\to \mathfrak{R}$. 

\vskip .2cm
Let $\lambda^{(m)}$ be an $m$-partition of size $n$ and let $X_{\lambda^{(m)}}$ and $X'_{\lambda^{(m)}}$ be two different standard $m$-tableaux of shape $\lambda^{(m)}$. For brevity we set $c^{(i)}=c(X_{\lambda^{(m)}}|i)$ for all $i=1,\dots,n$. The map 
$\langle ,\rangle$ is given by 
\begin{equation}\label{H-scal-prod1}
\langle\mathcal{X}_{\lambda^{(m)}},\mathcal{X}'_{\lambda^{(m)}}\rangle =0\ ,\end{equation}

\begin{equation}\label{H-scal-prod2} 
\langle\mathcal{X}_{\lambda^{(m)}},\mathcal{X}_{\lambda^{(m)}}\rangle =\prod\limits_{j,k\colon j<k\, ,\, c^{(j)}\neq c^{(k)}\, ,\, c^{(j)}\neq c^{(k)}q^{\pm2}}\ \frac{q^{-1}c^{(j)}-q c^{(k)}}{c^{(j)}-c^{(k)}}\ .
\end{equation}

Notice that, if $\xl{\lam}{s_i}$ is a standard $m$-tableau, we have
\begin{equation}\label{H-eq-scal-prod1}
\langle\ml{\lam}{s_i},\ml{\lam}{s_i}\rangle =\frac{q c^{(i)}-q^{-1}c^{(i+1)}}{q^{-1}c^{(i)}-q c^{(i+1)}}\langle\mathcal{X}_{\lambda^{(m)}},\mathcal{X}_{\lambda^{(m)}}\rangle \ .
\end{equation}

\paragraph{New basis.} In the generic situation, analogues for $H(m,1,n)$, corresponding to the shape $\lambda^{(m)}$, of the orthogonal representation
of the symmetric group will 
be defined over an extension $\tilde{\mathfrak{R}}_{\lambda^{(m)}}$ of the ring $\mathfrak{R}$. 
Let $\mathfrak{F}_{\lambda^{(m)}}$ be the set of all standard $m$-tableaux of 
shape $\lambda^{(m)}$. Let also $\mathfrak{E}_{_{\scriptstyle X_{\lambda^{(m)}}}}:=
\{ (j,k) \vert j<k\, ,\, c^{(j)}\neq c^{(k)}\, ,\, c^{(j)}\neq c^{(k)}q^{\pm2}\}$ (this is the set over which the product in the right hand side of (\ref{H-scal-prod2}) is taken). 
Introduce, for each standard $m$-tableau $X_{\lambda^{(m)}}$, a collection of variables 
$\varsigma^{jk}_{_{\scriptstyle X_{\lambda^{(m)}}}}$ and let 
$$\tilde{\mathfrak{R}}_{\lambda^{(m)}}:=\mathfrak{R}\left[ \varsigma^{jk}_{_{\scriptstyle X_{\lambda^{(m)}}}}\right]_{X_{\lambda^{(m)}}\in 
\mathfrak{F}_{\lambda^{(m)}},(j,k)\in \mathfrak{E}_{_{\scriptstyle X_{\lambda^{(m)}}}}}\text{\huge /}\mathfrak{I}\ ,$$ 
where $\mathfrak{I}$ is the ideal generated by 
$\left\{\left(\varsigma^{jk}_{_{\scriptstyle X_{\lambda^{(m)}}}}\right)^2-
\displaystyle{\frac{c^{(j)}-c^{(k)}}{q^{-1}c^{(j)}-q c^{(k)}}}\right\}_{X_{\lambda^{(m)}}\in \mathfrak{F}_{\lambda^{(m)}},(j,k)\in 
\mathfrak{E}_{_{\scriptstyle X_{\lambda^{(m)}}}}}$.
 
\vskip .2cm 
We define a new basis $\{ \tilde{{\mathcal{X}}}_{\lambda^{(m)}} \}$ of the $\tilde{\mathfrak{R}}_{\lambda^{(m)}}$-module 
$\tilde{\mathfrak{R}}_{\lambda^{(m)}}\otimes_{\mathfrak{R}}U_{\lambda^{(m)}}$ by
the following (diagonal) change of basis:
\begin{equation}\label{anoforfo}
{\mathcal{X}}_{\lambda^{(m)}}\mapsto \tilde{{\mathcal{X}}}_{\lambda^{(m)}}:=
\mathfrak{d}_{_{\scriptstyle\mathcal{X}_{\lambda^{(m)}}}}\ \mathcal{X}_{\lambda^{(m)}}\ ,
\ \ \text{where}\ \ 
\mathfrak{d}_{_{\scriptstyle\mathcal{X}_{\lambda^{(m)}}}}=\hspace{-.6cm}\prod\limits_{\hspace{1cm} j,k\colon 
(j,k)\in \mathfrak{E}_{_{\scriptstyle X_{\lambda^{(m)}}}}}\ \varsigma^{jk}_{_{\scriptstyle X_{\lambda^{(m)}}}} 
\end{equation}
for any standard $m$-tableau $X_{\lambda^{(m)}}$ of shape $\lambda^{(m)}$. 

\paragraph{1. Bilinear scalar product.} $ $

\vskip .2cm
{\bf 1.1.} Let $U$ be a free module over $\mathfrak{R}$. We shall call 
{\it bilinear scalar product} on $U$ a bi-additive map $\langle,\rangle_{\text{bilin}} \colon U\times U\to \mathfrak{R}$ with the property
\begin{equation}\label{herm}\langle f\, u,g\, v\rangle_{\text{bilin}} =f\, g\, \langle u,v\rangle_{\text{bilin}} \ ,\ f,g\in \mathfrak{R}\,,\ \ u,v\in U\ .\end{equation}
Let $\mathfrak{B}$ be a basis of $U$. An arbitrary map $\mathfrak{B}\times\mathfrak{B}\to \mathfrak{R}$ extends, by bi-additivity and (\ref{herm}), to a bilinear scalar product on $U$.  

\vskip .2cm
Thus the formulas (\ref{H-scal-prod1})--(\ref{H-scal-prod2}) define a bilinear scalar product $\langle,\rangle_{\text{bilin}}$ 
on the $H(m,1,n)$-module $U_{\lambda^{(m)}}$. Recall that $\mathsf{S}=\{\tau,\sigma_1,\dots ,\sigma_{n-1}\}$.

\begin{prop}{\hspace{-.2cm}.\hspace{.2cm}}
\label{prop-scal-prod}
The bilinear scalar product $\langle,\rangle_{\text{bilin}}$ on $U_{\lambda^{(m)}}$ is $\mathsf{S}^-$-invariant.
\end{prop}
\emph{Proof.} 
It is immediate that (\ref{scaprose}) is valid for the generator $\tau$ of $H(m,1,n)$. The verification of 
the $\mathsf{S}^-$-invariance of (\ref{H-scal-prod1}) for the generator $\sigma_i$ of $H(m,1,n)$ is non-trivial only if $\xl{\lam}{s_i}$ is standard and $X'_{\lambda^{(m)}}=\xl{\lam}{s_i}$. This verification is done by a straightforward calculation of $\langle\sigma_i(\mathcal{X}_{\lambda^{(m)}}),\sigma_i^{-1}(\ml{\lam}{s_i})\rangle_{\text{bilin}} $ (the action of $\sigma_i^{-1}$ on the two-dimensional space spanned by 
$\mathcal{X}_{\lambda^{(m)}},\mathcal{X}^{s_i}_{\lambda^{(m)}}$
is easily obtained from (\ref{rep-sigma-two})); the product $\langle\sigma_i(\mathcal{X}_{\lambda^{(m)}}),\sigma_i^{-1}(\ml{\lam}{s_i})\rangle_{\text{bilin}}$  is equal to 
\[
-\frac{(q-q^{-1})c^{(i+1)}}{(c^{(i)}-c^{(i+1)})^2}\Bigl((q c^{(i)}-q^{-1}c^{(i+1)})\langle\mathcal{X}_{\lambda^{(m)}},\mathcal{X}_{\lambda^{(m)}}\rangle_{\text{bilin}} +(q c^{(i+1)}-q^{-1}c^{(i)})\langle\ml{\lam}{s_i},\ml{\lam}{s_i}\rangle_{\text{bilin}} \Bigr)\ ,\]
which is $0$ due to (\ref{H-eq-scal-prod1}). 

\vskip .2cm
If $\xl{\lam}{s_i}$ is not standard then $\sigma_i(\mathcal{X}_{\lambda^{(m)}})=\epsilon q^{\epsilon}\mathcal{X}_{\lambda^{(m)}}$, 
$\epsilon =\pm 1$, and the $\mathsf{S}^-$-invariance of (\ref{H-scal-prod2}) for $\sigma_i$ follows. 

\vskip .2cm
If $\xl{\lam}{s_i}$ is standard, then a direct calculation gives that $\langle\sigma_i(\mathcal{X}_{\lambda^{(m)}}),\sigma_i^{-1}(\mathcal{X}_{\lambda^{(m)}})\rangle_{\text{bilin}} $ is equal to 
\[\frac{(q-q^{-1})^2c^{(i)}c^{(i+1)}\langle\mathcal{X}_{\lambda^{(m)}},\mathcal{X}_{\lambda^{(m)}}\rangle_{\text{bilin}} +(q^{-1}c^{(i)}-q c^{(i+1)})^2\langle\ml{\lam}{s_i},\ml{\lam}{s_i}\rangle_{\text{bilin}} }{(c^{(i)}-c^{(i+1)})^2}\ .\]
Using (\ref{H-eq-scal-prod1}) one obtains 
\[\langle\sigma_i(\mathcal{X}_{\lambda^{(m)}}),\sigma_i^{-1}(\mathcal{X}_{\lambda^{(m)}})\rangle_{\text{bilin}} =\langle\mathcal{X}_{\lambda^{(m)}},\mathcal{X}_{\lambda^{(m)}}\rangle_{\text{bilin}} \ .\]
The proof of the Proposition is finished.\hfill$\square$

\paragraph{Remarks.}$\ $

\vskip .2cm
\textbf{(a)} Note that (\ref{herm}) is compatible with any specialization of the parameters $q$ and $v_j$, $j=1,\dots,m$.

\vskip .2cm
\textbf{(b)}  For $m=1$, that is, for the usual Hecke algebra $H(1,1,n)$, we have $c^{(i)}=q^{2 cc^{(i)}}$ where $cc^{(i)}$ is the classical content of the node occupied by $i$ 
and the formula (\ref{H-scal-prod2}) can be rewritten in the following form 
\begin{equation}\label{H-scal-prod-He} 
\langle\mathcal{X}_{\lambda^{(m)}},\mathcal{X}_{\lambda^{(m)}}\rangle_{\text{bilin}} =\prod\limits_{j,k\colon j<k\, ,\, c^{(j)}\neq c^{(k)}\, ,\, c^{(j)}\neq c^{(k)}q^{\pm2}}\ 
\frac{(cc^{(j)}-cc^{(k)}-1)_q}{(cc^{(j)}-cc^{(k)})_q}\ .
\end{equation}

\vskip .3cm
{\bf 1.2.} Let $\mathsf{A}$ be the matrix of the generator $\mathsf{a}\in\{\tau,\sigma_1,\dots ,\sigma_{n-1}\}$ of 
$H(m,1,n)$ in the basis $\{ \tilde{{\mathcal{X}}}_{\lambda^{(m)}}\}$, see (\ref{anoforfo}). Then 
\begin{equation}\label{mafoor}\mathsf{A}\,(\mathsf{A}^{-1})^T=\textbf{Id}\ ,\end{equation} 
where $\textbf{Id}$ is the identity matrix and, given a matrix $x$, $x^T$ means the transposed matrix.

\paragraph{2. Sesquilinear scalar products.} $ $

\vskip .2cm
{\bf 2.1.} For any Laurent polynomial $f\in\mathbb{C}[q,q^{-1},v_1,v_1^{-1},\dots,v_m,v_m^{-1}]$, denote by $\gamma(f)$ the Laurent polynomial in $q,v_1,\dots,v_m$ obtained by conjugating the coefficients of $f$. The involution $\gamma$ extends to the ring $\mathfrak{R}$. Let $U$ be a free $\mathfrak{R}$-module. We shall call 
{\it $\gamma$-sesquilinear scalar product} on $U$ a bi-additive map $\langle,\rangle_{\gamma} \colon U\times U\to \mathfrak{R}$ with the property (replacing (\ref{herm}))
\begin{equation}\label{herm2}\langle f\, u,g\, v\rangle_{\gamma} =f\, \gamma(g)\, \langle u,v\rangle_{\gamma} \ ,\ f,g\in \mathfrak{R}\,,\ \ u,v\in U\ .\end{equation}
Let $\mathfrak{B}$ be a basis of $U$. An arbitrary map $\mathfrak{B}\times\mathfrak{B}\to \mathfrak{R}$ extends, by bi-additivity and (\ref{herm2}), to a $\gamma$-sesquilinear scalar product on $U$.  

\vskip .2cm
Thus the formulas (\ref{H-scal-prod1})--(\ref{H-scal-prod2}) define a $\gamma$-sesquilinear scalar product on $U_{\lambda^{(m)}}$ which is $\mathsf{S}^-$-invariant, where $\mathsf{S}=\{\tau,\sigma_1,\dots,\sigma_{n-1}\}$; the verification repeats the proof of the Proposition \ref{prop-scal-prod}.

\vskip .2cm
The $\mathbb{C}$-specializations compatible with the definition of a $\gamma$-sesquilinear scalar product are those for which $q$ and $v_j$, $j=1,\dots,m$, are real numbers; for such specializations, 
the $\gamma$-sesquilinear scalar product becomes a usual Hermitian scalar product on a complex vector space.

\vskip .3cm
\textbf{2.2.} Denote by $\omega$ the involution on 
$\mathbb{C}[q,q^{-1},v_1,v_1^{-1},\dots,v_m,v_m^{-1}]$ which sends $q$ to $q^{-1}$ and $v_j$ to $v_j^{-1}$, $j=1,\dots ,m$. The involution $\omega$
is compatible with the localization with respect to the multiplicative set generated by $\mathfrak{D}$ and thus extends to $\mathfrak{R}$. Let $U$ be as before a free module over $\mathfrak{R}$. We shall call 
{\it $\omega$-sesquilinear scalar product} a bi-additive map $\langle,\rangle_{\omega} \colon U\times U\to \mathfrak{R}$ with the property
\begin{equation}\label{herm-w}\langle f\, u,g\, v\rangle_{\omega} =f\, \omega(g)\, \langle u,v\rangle_{\omega} \ ,\ f,g\in \mathfrak{R}
\,,\ \ u,v\in U\ .\end{equation}
Let $\mathfrak{B}$ be a basis of $U$. An arbitrary map $\mathfrak{B}\times\mathfrak{B}\to \mathfrak{R}$ extends, by bi-additivity and (\ref{herm-w}), to an $\omega$-sesquilinear scalar product on $U$.  

\vskip .2cm
Thus the formulas (\ref{H-scal-prod1})--(\ref{H-scal-prod2}) define an $\omega$-sesquilinear scalar product on $U_{\lambda^{(m)}}$ which is $\mathsf{S}$-invariant, where $\mathsf{S}=\{\tau,\sigma_1,\dots,\sigma_{n-1}\}$; the verification is similar to the proof of the Proposition \ref{prop-scal-prod}.

\vskip .2cm
The only $\mathbb{C}$-specializations compatible with the definition of an $\omega$-sesquilinear scalar product are those for which $q$ and $v_j$, $j=1,\dots,m$, belong to $\{-1,1\}$.

\vskip .3cm
\textbf{2.3.} The involutions $\gamma$ and $\omega$ commute. Let $\varpi$ be the involution on $\mathfrak{R}$ defined as the composition of $\gamma$ and $\omega$, $\varpi(f):=\omega\bigl(\gamma(f)\bigr)$ for $f\in\mathfrak{R}$.
We shall call {\it $\varpi$-sesquilinear scalar product} a bi-additive map $\langle,\rangle_{\varpi} \colon U\times U\to \mathfrak{R}$ with the property
\begin{equation}\label{herm-w2}\langle f\, u,g\, v\rangle_{\varpi} =f\, \varpi(g)\, \langle u,v\rangle_{\varpi} \ ,\ f,g\in \mathfrak{R}
\,,\ \ u,v\in U\ .\end{equation} 
Let $\mathfrak{B}$ be a basis of $U$. An arbitrary map $\mathfrak{B}\times\mathfrak{B}\to \mathfrak{R}$ extends, by bi-additivity and (\ref{herm-w2}), to a $\varpi$-sesquilinear scalar product on $U$.  

\vskip .2cm
Thus the formulas (\ref{H-scal-prod1})--(\ref{H-scal-prod2}) define a $\varpi$-sesquilinear scalar product on $U_{\lambda^{(m)}}$ which is $\mathsf{S}$-invariant, where $\mathsf{S}=\{\tau,\sigma_1,\dots,\sigma_{n-1}\}$; the verification is similar to the proof of the Proposition \ref{prop-scal-prod}.

\vskip .2cm
The $\mathbb{C}$-specializations compatible with the definition of a $\varpi$-sesquilinear scalar product are those for which $q$ and $v_j$, $j=1,\dots,m$, are complex numbers on the unit circle; for such specializations, the $\varpi$-sesquilinear scalar product becomes a usual Hermitian scalar product on a complex vector space.

\vskip .3cm
\textbf{2.4.} Reintroduce the deformation parameters $q,v_1,\dots ,v_m$ in the notation for the cyclotomic Hecke algebra: 
$H_{q,v_1,\dots ,v_m}(m,1,n)$. There is another way to interpret the formulas (\ref{H-scal-prod1})-(\ref{H-scal-prod2}) in the $\omega$-sesquilinear and the $\varpi$-sesquilinear situations. Namely these formulas define a pairing $\langle,\rangle$ between the representation spaces $U_{\lambda^{(m)}}$ and $U'_{\lambda^{(m)}}$ where the first space 
$U_{\lambda^{(m)}}$ carries the representation of the algebra $H_{q,v_1,\dots ,v_m}(m,1,n)^{\phantom{A}}\!\!\!$ and the second space 
$U'_{\lambda^{(m)}}$ carries the representation of the algebra $H_{q^{-1},v_1^{-1},\dots ,v_m^{-1}}(m,1,n)^{\phantom{A}}\!\!\!$. In the $\omega$-sesquilinear situation, the pairing is bilinear in the usual sense, $\langle fu,gv\rangle^{\phantom{A}}\!\!\! =fg\langle u,v\rangle$, $f,g\in \mathfrak{R}$. In the $\varpi$-sesquilinear situation, the pairing is $\gamma$-sesquilinear, $\langle fu,gv\rangle=f\gamma(g)\langle u,v\rangle $, $f,g\in \mathfrak{R}^{^{\phantom{A}}}\!\!\!$.

\vskip .2cm
Both $U_{\lambda^{(m)}}$ and $U'_{\lambda^{(m)}}$ have the same basis $\mathfrak{B}=\{\mathcal{X}_{\lambda^{(m)}}\}$. The formula (\ref{scaprose2}), stating the $\mathsf{S}$-invariance of the
pairing, holds;  
now in the formula (\ref{scaprose2}), $x(u)^{^{\phantom{A}}}\!\!\!$ is to be understood as the result of the action of the generator $x\in H_{q,v_1,\dots ,v_m}(m,1,n)$ on the vector $u\in U_{\lambda^{(m)}}$ while $x(v)^{^{\phantom{A}}}\!\!\!$ is the result of the action of the generator $x\in H_{q^{-1},v_1^{-1},\dots ,v_m^{-1}}(m,1,n)^{^{\phantom{A}}}\!\!\!$ on the vector $v\in U'_{\lambda^{(m)}}$.

\vskip .3cm
\textbf{2.5.} Since each factor $\displaystyle{\frac{c^{(j)}-c^{(k)}}{q^{-1}c^{(j)}-q c^{(k)}}}$ in the product in the right hand side of (\ref{H-scal-prod2})
is stable with respect to the involutions $\gamma$ and $\omega$, one can extend the involutions $\gamma$ and $\omega$ to the 
ring $\tilde{\mathfrak{R}}_{\lambda^{(m)}}$ by the rules $\gamma \left(\varsigma^{jk}_{_{\scriptstyle X_{\lambda^{(m)}}}}\right) =
\varsigma^{jk}_{_{\scriptstyle X_{\lambda^{(m)}}}}$ and $\omega \left(\varsigma^{jk}_{_{\scriptstyle X_{\lambda^{(m)}}}}\right) =
\varsigma^{jk}_{_{\scriptstyle X_{\lambda^{(m)}}}}$ for all $X_{\lambda^{(m)}}\in \mathfrak{F}_{\lambda^{(m)}}$ and $(j,k)\in 
\mathfrak{E}_{_{\scriptstyle X_{\lambda^{(m)}}}}$. 
 
\vskip .2cm 
Let $\mathsf{A}$ be the matrix of the generator $\mathsf{a}\in\{\tau,\sigma_1,\dots ,\sigma_{n-1}\}$ of 
$H(m,1,n)$ in the basis $\{ \tilde{{\mathcal{X}}}_{\lambda^{(m)}}\}$, see (\ref{anoforfo}). Then 
\begin{equation}\label{mafoor2}\mathsf{A}\, \gamma(\mathsf{A}^{-1})^T=\textbf{Id}\ ,\ \ \  \mathsf{A}\, \omega(\mathsf{A})^T=\textbf{Id}\ \ \ \text{and}\ \ \ \mathsf{A}\, \varpi(\mathsf{A})^T=\textbf{Id}\ .\end{equation}

For a specialization in which $q$ and $v_j$, $j=1,\dots,m$, are complex numbers on the unit circle, the matrices of the generators of $H(m,1,n)$ in the basis $\{ \tilde{{\mathcal{X}}}_{\lambda^{(m)}}\}$ are unitary in the usual sense.

\paragraph{3.} Let $\rho\colon \hat{H}_n\to\text{End}({\cal{V}})$ be an irreducible representation of the affine Hecke algebra $\hat{H}_n$ on a 
complex vector space ${\cal{V}}$ of finite dimension $L$. Assume that the operator $\rho (\tau )$ is diagonalizable and the spectrum of
$\rho (\tau )$ is $\{ (v_1)_{\mathfrak{l}_1},(v_2)_{\mathfrak{l}_2},\dots ,(v_m)_{\mathfrak{l}_m}\}$; here the numbers $\{ v_1,v_2,\dots ,v_m\}$ are pairwise
different, ${\mathfrak{l}_j}$ is the multiplicity of the eigenvalue $v_j$, $j=1,\dots ,m$. Then the representation $\rho$ passes through the cyclotomic quotient 
$H_{q,v_1,\dots ,v_m}(m,1,n)$ of the affine Hecke algebra. Assume that the parameters $q,v_1,\dots ,v_m$ satisfy the
restrictions (\ref{sesi1})-(\ref{sesi3}). 
By the completeness result from Section \ref{sec-comp}, the representation $\rho$ is isomorphic, as the representation of $H_{q,v_1,\dots ,v_m}(m,1,n)$, 
to $V_{\lambda^{(m)}}$ for a certain $m$-partition $\lambda^{(m)}$. Note that, given the knowledge of the values of $q$ and $v_j$, $j=1,\dots ,m$,
the basis $\{ {\mathcal{X}}_{\lambda^{(m)}}\}$ is determined uniquely up to a global rescaling -- if there were two bases then the operator, transforming 
one into another, would contradict to the irreducibility of the representation $\rho$. 
The diagonal change of basis does not pose any problem: in the product in (\ref{anoforfo}), 
$\varsigma^{jk}_{_{\scriptstyle X_{\lambda^{(m)}}}}$ can be chosen to be an arbitrary square root of $\displaystyle{\frac{c^{(j)}-c^{(k)}}{q^{-1}c^{(j)}-q c^{(k)}}}$, $(j,k)\in\mathfrak{E}_{_{\scriptstyle X_{\lambda^{(m)}}}}$.

Let $\mathsf{A}$ be the matrix of the generator $\mathsf{a}\in\{\tau,\sigma_1,\dots ,\sigma_{n-1}\}$ of 
$\hat{H}_n$ in the basis $\{ \tilde{{\mathcal{X}}}_{\lambda^{(m)}}\}$. Then
\[\mathsf{A}\, (\mathsf{A}^{-1})^T=\textbf{Id}_{{\cal{V}}}\ ,\]
where $\textbf{Id}_{{\cal{V}}}$ is the identity operator on the space ${\cal{V}}$.

\vskip .2cm
If the parameters $q$ and $v_j$, $j=1,\dots,m$ are complex numbers on the unit circle, the matrices of the generators, in the basis $\{ \tilde{{\mathcal{X}}}_{\lambda^{(m)}}\}$, are unitary in the usual sense, $\mathsf{A}\, \mathsf{A}^{\dagger}=\textbf{Id}_{{\cal{V}}}$, where $A^{\dagger}$ is the transposed conjugated matrix.
 
\setcounter{equation}{0}
\section{Completeness}\label{sec-comp}

\paragraph{1.} In the preceding Section we constructed, for every $m$-partition $\lambda^{(m)}$, a representation of $H(m,1,n)$. The spectrum of the 
Jucys--Murphy elements $J_1,\dots,J_n$ in this representation is the set of strings corresponding to the standard $m$-tableaux of shape $\lambda^{(m)}$, 
see the Lemma \ref{prop-rel2}. This construction provides an inclusion of the set of standard Young $m$-tableaux of size $n$ into 
${\mathrm{Spec}}(J_1,\dots,J_n)$. 
On the other hand, the Proposition \ref{prop3} and the Proposition \ref{prop4} provide an inclusion of ${\mathrm{Spec}}(J_1,\dots,J_n)$ into 
the set of standard Young $m$-tableaux of size $n$. These operations, by construction, are inverse to each other. We sum up the results.

\vskip .2cm
We underline that the restrictions (\ref{sesi1})--(\ref{sesi3}) are essential for the statements below. 

\begin{prop}
{\hspace{-.2cm}.\hspace{.2cm}}
 \label{prop-cont-spec}
The set ${\mathrm{Spec}}(J_1,\dots,J_n)$, the set ${\mathrm{Cont}}_m(n)$ and the set of standard $m$-tableaux are in bijection.
\end{prop}

\begin{coro}
{\hspace{-.2cm}.\hspace{.2cm}}
 \label{lem-fin}
The spectrum of the Jucys--Murphy elements is simple in the representations $V_{\lambda^{(m)}}$ (labeled by the $m$-partitions).
\end{coro}
It means that for two different standard $m$-tableaux (not necessarily of the same shape) the elements of ${\mathrm{Spec}}(J_1,\dots,J_n)$ associated to them by the Proposition \ref{prop-cont-spec} are different (two strings $(a_1,\dots,a_n)$ and $(a'_1,\dots,a'_n)$ are different if there is some $i$ such that $a_i\neq a'_i$).

\paragraph{2.} It remains to verify that we obtain within this approach all irreducible representations of the algebra $H(m,1,n)$. 

\vskip .2cm
According to Appendix B, see (\ref{sumdim2}),
the sum of the squares of the dimensions of the constructed representations
is not less than
the dimension of $H(m,1,n)$. It is therefore sufficient to prove
that these representations are irreducible and pairwise non-isomorphic. It is done in the Theorem \ref{prop-fin}. 

\vskip .2cm
Similarly to \cite{AK}, we obtain as a by-product that, under the restrictions (\ref{sesi1})--(\ref{sesi3}), $$\dim(H(m,1,n))=n!m^n\ ,$$ and the algebra $H(m,1,n)$ is semi-simple.

\begin{theo} 
{\hspace{-.2cm}.\hspace{.2cm}}
 \label{prop-fin}
The representations $V_{\lambda^{(m)}}$ (labeled by the $m$-partitions) of the algebra $H(m,1,n)$ constructed in Section \ref{ygcyclo'} are 
irreducible and pairwise non-isomorphic.
\end{theo}

\emph{Proof.} The proof can be found in \cite{AK} (as well as the semi-simplicity result). We briefly repeat the argument for completeness.

\vskip .2cm
In the proof we use induction on $n$. It is justified since the restrictions (\ref{sesi1})--(\ref{sesi3}) for $H(m,1,n)$ imply the 
restrictions (\ref{sesi1})--(\ref{sesi3}), in which $n$ is replaced by $n'$, for $H(m,1,n')$ with arbitrary $n'$, $0<n'<n$. 

\vskip .2cm
The Corollary \ref{lem-fin} directly implies that the representations $V_{\lambda^{(m)}}$ and $V_{\lambda'^{(m)}}$ are non-isomorphic if 
$\lambda^{(m)}\neq\lambda'^{(m)}$.

\vskip .2cm
Suppose by induction that the representations $V_{\mu^{(m)}}$, for all $m$-partitions $\mu^{(m)}$ of $n-1$, are irreducible representations of 
$H(m,1,n-1)$. The base of induction is $n=1$; there is nothing to prove here.

\vskip .2cm
Fix an $m$-partition $\lambda^{(m)}$ with $|\lambda^{(m)}|=n$. Let $\{ \mu_i^{(m)}\}$, $i=1,\dots,l$, be the set of all $m$-sub-partitions of $\lambda^{(m)}$ with $|\mu_i^{(m)}|=n-1$.

\vskip .2cm
{}For each $i$, the representation $V_{\mu_i^{(m)}}$ of $H(m,1,n-1)$ is a sub-representation of the restriction of the representation $V_{\lambda^{(m)}}$ 
to $H(m,1,n-1)$. The dimension of $V_{\lambda^{(m)}}$ (the number of standard $m$-tableaux of the shape $\lambda^{(m)}$) is the sum (over $i$) of 
dimensions of $V_{\mu_i^{(m)}}$. Therefore, the representation $V_{\lambda^{(m)}}$ of $H(m,1,n)$ decomposes with respect to $H(m,1,n-1)$ into a direct 
sum of representations $V_{\mu_i^{(m)}}$. 

\vskip .2cm
The $m$-sub-partitions $\mu_i^{(m)}$ are different and correspond thus to non-isomorphic irreducible representations 
of $H(m,1,n-1)$. It follows that the positions of $V_{\mu_i^{(m)}}$ as subspaces in $V_{\lambda^{(m)}}$ are well-defined. Therefore, if $V_{\lambda^{(m)}}$ 
has a non-trivial invariant subspace $U$ then $U$ must contain at least one of the $V_{\mu_i^{(m)}}$, say $V_{\mu_1^{(m)}}$.

\vskip .2cm
It is sufficient to show that starting from elements of $V_{\mu_1^{(m)}}$ one can obtain an 
element of $V_{\mu_j^{(m)}}$ for any $j\neq1$ by the action of operators from $H(m,1,n)$. A basis vector of $V_{\lambda^{(m)}}$ labeled by a 
standard $m$-tableau $X_{\lambda^{(m)}}$ of shape $\lambda^{(m)}$ belongs to the subspace $V_{\mu_j^{(m)}}$ where $\mu_j^{(m)}$ is the 
$m$-sub-partition of 
size $(n-1)$ formed by the $m$-nodes with $1,\dots,n-1$ of $X_{\lambda^{(m)}}$.
For any $j\neq1$ the $m$-sub-partition $\mu_j^{(m)}$ is obtained from $\mu_1^{(m)}$ by 
removing one $m$-node and adding some other $m$-node, different from the removed one; it is easy to see that the two $m$-nodes involved are non-adjacent 
and, even more, are not on neighboring diagonals. Take the 
standard $m$-tableau of shape $\lambda^{(m)}$ for which the numbers $1,\dots,n-1$ are placed in the $m$-sub-partition $\mu_1^{(m)}$ of 
$\lambda^{(m)}$ and moreover the number $n-1$ is in the only $m$-node of $\mu_1^{(m)}$ which is not in the $m$-sub-partition $\mu_j^{(m)}$. 
The vector $\vec{v}$ of $V_{\lambda^{(m)}}$ labeled by this $m$-tableau 
belongs to the subspace $V_{\mu_1^{(m)}}$ and is sent by $\sigma_{n-1}$ to a combination of the vector $\vec{v}$ and a vector belonging to 
$V_{\mu_j^{(m)}}$. The formula (\ref{rep-a1}) shows that this vector of $V_{\mu_j^{(m)}}$ is non-zero. \hfill$\square$

\bigskip
Let ${\mathcal{B}}$ be an associative subalgebra of an associative algebra ${\mathcal{A}}$. An indecomposable (irreducible if the algebra 
${\mathcal{A}}$ is semi-simple) representation of the algebra ${\mathcal{A}}$ ``branches" with respect to the algebra ${\mathcal{B}}$, that is, decomposes 
into a direct sum of indecomposable (irreducible if the algebra ${\mathcal{B}}$ is semi-simple) 
representations of ${\mathcal{B}}$. The information about branchings of all representations of the algebra ${\mathcal{A}}$ with respect to the
subalgebra ${\mathcal{B}}$ is called \emph{branching rules} for the pair $({\mathcal{A}},{\mathcal{B}})$.

\vskip .2cm
As a corollary of the whole construction we obtain under the restrictions (\ref{sesi1})--(\ref{sesi3}) the branching rules
for the pair $\bigl( H(m,1,n),H(m,1,n-1)\bigr)$; the representation of the algebra $H(m,1,n)$ labeled by an $m$-partition $\lambda^{(m)}$ of $n$
decomposes into the direct sum of the representations of the algebra $H(m,1,n)$ labeled by the $m$-sub-partitions of $\lambda^{(m)}$ of size $n-1$.
In particular we obtain the following Corollary.

\begin{coro}
{\hspace{-.2cm}.\hspace{.2cm}}
 \label{coro-sim}
Under the restrictions (\ref{sesi1})--(\ref{sesi3}) the branching rules for the chain, with respect to $n$, of the algebras $H(m,1,n)$ are multiplicity-free.
\end{coro}

It means that under the restrictions (\ref{sesi1})--(\ref{sesi3}) in the decomposition of an irreducible representation of the algebra $H(m,1,n)$ each 
irreducible representation of the algebra $H(m,1,n-1)$ appears with the multiplicity equal either to $0$ or to $1$.

\vskip .2cm
By the general arguments it follows that under the restrictions (\ref{sesi1})--(\ref{sesi3}) the centralizer of the sub-algebra $H(m,1,n-1)$ in the 
algebra $H(m,1,n)$ is commutative for each $n=1,2,3,\dots$ 

\vskip .2cm
It also follows from the constructed representation theory that under the restrictions (\ref{sesi1})--(\ref{sesi3}) 
\begin{itemize}
\item the centralizer of the subalgebra $H(m,1,n-1)$ in the algebra $H(m,1,n)$ is 
generated by the center of $H(m,1,n-1)$ and the Jucys--Murphy element $J_n$; 
\item the subalgebra generated by the Jucys--Murphy elements $J_1,\dots,J_n$ of the algebra $H(m,1,n)$ is maximal commutative.
\end{itemize} 

\paragraph{Remarks.}$\ $

\vskip .2cm
 \textbf{(a)} For every standard $m$-tableau $X_{\lambda^{(m)}}$ define the element $\mathsf{E}_{_{\scriptstyle{X_{\lambda^{(m)}}}}}$ 
of the cyclotomic Hecke algebra $H(m,1,n)$ by the following recursion. The initial condition is $\mathsf{E}_{\varnothing}=1$. Let $\alpha^{(m)}$ be the
$m$-node occupied by the number $n$ in $X_{\lambda^{(m)}}$; define $\mu^{(m)}:=\lambda^{(m)}\backslash\{ \alpha^{(m)}\}$. Denote by 
$X_{\mu^{(m)}}$ the standard $m$-tableau with the numbers $1,\dots ,n-1$ at the same $m$-nodes as in $X_{\lambda^{(m)}}$. Then the recursion is given by
\begin{equation}\label{nopr} \mathsf{E}_{_{\scriptstyle{X_{\lambda^{(m)}}}}}:= \mathsf{E}_{_{\scriptstyle{X_{\mu^{(m)}}}}}
\prod_{\beta^{(m)}\colon \beta^{(m)}\in {\cal{E}}_+(\mu^{(m)})\, ,\, \beta^{(m)}\neq \alpha^{(m)}      }
\frac{J_n-c(\beta^{(m)})}{c(\alpha^{(m)})-c(\beta^{(m)})}\end{equation}
where $c(\beta^{(m)})$ is the content of the $m$-node $\beta^{(m)}$. Due to the completeness results of this Section, the elements $\mathsf{E}_{_{\scriptstyle{X_{\lambda^{(m)}}}}}$ form a complete set of pairwise orthogonal primitive idempotents of the algebra $H(m,1,n)$.

 \vskip .2cm
We shall prove that, moreover, we have a well-defined homomorphism $\varrho\,\colon\ {\mathfrak{T}}\to H(m,1,n)$ which is identical on the generators $\tau,\sigma_1,\dots,\sigma_{n-1}$ and sends $\mathcal{X}_{\lambda^{(m)}}$ to $\mathsf{E}_{_{\scriptstyle{X_{\lambda^{(m)}}}}}$ for all standard $m$-tableaux $X_{\lambda^{(m)}}$. Using the completeness, the only non-trivial verification one has to do is to check that, for any standard $m$-tableau $X_{\lambda^{(m)}}$ such that $X_{\lambda^{(m)}}^{s_i}$ is standard, the defining relation (\ref{rel-a1}) of the algebra ${\mathfrak{T}}$ is satisfied by the images of $\sigma_i$, $\mathcal{X}_{\lambda^{(m)}}$ and $\mathcal{X}_{\lambda^{(m)}}^{s_i}$ through the homomorphism $\varrho$. The verification reduces to the following equality for matrices (see (\ref{rep-sigma-two})): 
\begin{equation}\label{offdia}\left(\begin{array}{cc}0 & A\\
B & -C\end{array}\right)\left(\begin{array}{cc}1&0\\0&0\end{array}\right)=
\left(\begin{array}{cc}0&0\\B&0\end{array}\right)=
\left(\begin{array}{cc}0&0\\0&1\end{array}\right)\left(\begin{array}{cc}C & A\\
B & 0\end{array}\right)\ ,\end{equation}
where $A={\displaystyle \frac{qc^{(i)}-q^{-1}c^{(i+1)}}{c^{(i)}-
c^{(i+1)}}}$, $B={\displaystyle \frac{qc^{(i+1)}-q^{-1}c^{(i)}}{c^{(i+1)}-
c^{(i)}}}$, $C={\displaystyle \frac{(q-q^{-1})(c^{(i+1)}+c^{(i)})}{c^{(i+1)}-
c^{(i)}}}$ and $c^{(i)}:=c(X_{\lambda^{(m)}}|i)$ for all $i=1,\dots,n$. The elements $\varrho (\mathcal{X}_{\lambda^{(m)}})$ are the diagonal matrix units;
the elements 
\begin{equation}\label{rel-a1rho}
\varrho\Bigl(\sigma_i+\frac{(q-q^{-1})
c(X_{\lambda^{(m)}}|i+1)}{c(X_{\lambda^{(m)}}|i)-c(X_{\lambda^{(m)}}|i+1)}\Bigr)\ \varrho(\mathcal{X}_{\lambda^{(m)}})=
\varrho (\ml{\lam}{s_i})\varrho\Bigl(\sigma_i+\frac{(q-q^{-1})c(X_{\lambda^{(m)}}|i)}{c(X_{\lambda^{(m)}}|i+1)-c(X_{\lambda^{(m)}}|i)}\Bigr)\end{equation}
form a part of off-diagonal (non-normalized) matrix units - the calculation (\ref{offdia}) shows that these elements are non-zero. For the usual Hecke 
algebra $H(1,1,n)$ the equality 
(\ref{rel-a1rho}) was established in \cite{Mu}. The complete set of off-diagonal matrix units was constructed in \cite{OP} for the usual Hecke algebra. The construction for the cyclotomic Hecke algebra is similar and we leave details to the reader. 

\vskip .3cm
\textbf{(b)} For a subset $\{ v_{i_1},\dots ,v_{i_l}\}$ with $l<m$ let $\mathfrak{z}:=(\tau -v_{i_1})...(\tau -v_{i_l})$. Taking a quotient of $H(m,1,n)$
by the ideal generated by $\mathfrak{z}$ we get a homomorphism $\mathfrak{p}\colon H(m,1,n)\to H(l,1,n)$ where $H(l,1,n)$ is the cyclotomic
Hecke algebra with the parameters $q$, $v_{i_1},\dots ,v_{i_l}$ (note that the restrictions (\ref{sesi1})--(\ref{sesi3}) hold for this choice of the parameters 
if they hold for  $q$, $v_1,\dots ,v_n$). 
The representations of $H(m,1,n)$
for which the diagonal entries of the (diagonal) matrix (\ref{rep-a2}) belong to $\{ v_{i_1},\dots ,v_{i_l}\}$ (these representations are labeled by 
$m$-partitions with empty partitions on the places corresponding to $v_j$ which are omitted in $\{ v_{i_1},\dots ,v_{i_l}\}$) 
pass through the image 
$\mathfrak{p}(H(m,1,n))$ in $H(m,1,n)$. The sum of squares of dimensions of these representations is not less than  
the dimension of $ H(l,1,n)$. It follows that 
$\mathfrak{p}$ is surjective. 

\renewcommand{\theequation}{{A}.{\arabic{equation}}}
\setcounter{equation}{0}
\section*{Appendix A.  \hspace{.2cm}  Module structure on tensor products}
\addcontentsline{toc}{section}{Appendix A.$\ $ Module structure on tensor products\vspace{.1cm}}

The algebra ${\mathfrak{T}}$, defined in Subsection \ref{subsection.smash.product}, was used in Subsection \ref{suserep} to construct
modules over the cyclotomic Hecke algebra $H(m,1,n)$. An extension of this construction equips the tensor products of the underlying spaces
of the $H(m,1,n)$-modules, corresponding to $m$-partitions of $n$, with a structure of an  $H(m,1,n)$-module. In this Appendix we give precise 
definitions and investigate the appearing tensor product, denoted by $\hat{\otimes}$, of representations.

\vskip .2cm
In Section \ref{sec-comp} we prove that under the restrictions (\ref{sesi1})--(\ref{sesi3}) on the parameters of $H(m,1,n)$, the irreducible
representations of $H(m,1,n)$ are exhausted by the representations corresponding to $m$-partitions of $n$. Our method of studying the tensor product
$\hat{\otimes}$ is inductive and is heavily based on the completeness result from Section \ref{sec-comp}. A priori, we do not know the nature of representations
appearing in the decomposition of the tensor product $V_{\lambda^{(m)}}\hat{\otimes}V_{\lambda'^{(m)}}$ of two representations corresponding to
$m$-partitions ${\lambda^{(m)}}$ and ${\lambda'^{(m)}}$. It is here that we need the heavy completeness result saying that every $H(m,1,n)$-module is 
isomorphic to a direct sum of $H(m,1,n)$-modules corresponding to $m$-partitions of $n$. Note that in its turn the completeness result is established in  
Section \ref{sec-comp} indirectly, by counting dimensions. 

\vskip .2cm
The decomposition rules for the tensor product $\hat{\otimes}$ are given in the Proposition \ref{tens-prod} in this Appendix. Qualitatively, the result is
formulated very easily: the tensor product $V_{\lambda^{(m)}}\hat{\otimes}V_{\lambda'^{(m)}}$ is isomorphic to the direct sum of dim$(V_{\lambda'^{(m)}})$ 
copies of the representation $V_{\lambda^{(m)}}$. For simplest choices of ${\lambda^{(m)}}$ and ${\lambda'^{(m)}}$, the decomposition of the tensor
product $V_{\lambda^{(m)}}\hat{\otimes}V_{\lambda'^{(m)}}$ can be done by a direct calculation. However, in spite of the easiness of the formulation of the result 
of the Proposition \ref{tens-prod},  we could not find a way to perform an explicit calculation for two arbitrary $m$-partitions.
Besides, rectangular partitions play a distinguished role in our way of proof of the Proposition \ref{tens-prod} but not in the final formula for the decomposition rules. 
 It would be interesting to find a more 
explicit way to establish the Proposition \ref{tens-prod}, without using the completeness assertion from Section \ref{sec-comp}. 

\vskip .2cm
By construction from Subsection \ref{suserep}, the representation $V_{\lambda^{(m)}}$ corresponding to an $m$-partition ${\lambda^{(m)}}$ is equipped 
with the natural basis $\mathcal{X}_{\lambda^{(m)}}$ indexed by standard $m$-tableaux of shape ${\lambda^{(m)}}$. The explicit form of an
isomorphism $V_{\lambda^{(m)}}\hat{\otimes}V_{\lambda'^{(m)}}\cong \dim(V_{\lambda'^{(m)}})\ V_{\lambda^{(m)}}$ 
is quite evolved, showing again that the tensor product $\hat{\otimes}$ requires further understanding. We give several examples.

\vskip .2cm
Certain statements below are valid in a more general situation, without the completeness assertion from Section \ref{sec-comp}. To accurately formulate these 
statements, we shall say that a representation of the cyclotomic Hecke algebra belongs to a class ${\cal{S}}$ if it is isomorphic to a direct sum of representations, 
corresponding to multi-partitions.

\subsection*{A.1  \hspace{.2cm} Definition of the tensor product and simplest examples}
\addcontentsline{toc}{subsection}{A.1$\ $ Module structure on tensor products}

\paragraph{1. Definition of the tensor product $\hat{\otimes}$.} Given an $m$-partition $\lambda^{(m)}$ of size $n$, recall that $U_{\lambda^{(m)}}$ is 
the vector space with the chosen basis $\{\mathcal{X}_{\lambda^{(m)}}\}$. 
We stress that $U_{\lambda^{(m)}}$ is understood only as a vector space, without any $H(m,1,n)$-module structure specified, whereas 
$V_{\lambda^{(m)}}$ is understood as the $H(m,1,n)$-module given by the formulas (\ref{rep-a1})-(\ref{rep-a2}) with underlying vector space 
$U_{\lambda^{(m)}}$. 
In particular, a representation of $H(m,1,n)$ is of the class ${\cal{S}}$ if it is isomorphic to a direct sum of representations of the form $V_{\lambda^{(m)}}$.

\vskip .2cm
Let ${\lambda^{(m)}}$ and ${\lambda'^{(m)}}$ be two $m$-partitions of size $n$. The instructions from
the Proposition \ref{prop-rel} are homogeneous in the generators $\mathcal{X}$. A basis of the tensor product of
$U_{\lambda^{(m)}}$ and $U_{\lambda'^{(m)}}$ is naturally indexed
 by the products $\mathcal{X}_{\lambda^{(m)}}\mathcal{X}_{\lambda'^{(m)}}$,
where $\mathcal{X}_{\lambda^{(m)}}$ is the generator labeled by the standard $m$-tableau $X_{\lambda^{(m)}}$ (of the shape ${\lambda^{(m)}}$) 
and $\mathcal{X}_{\lambda'^{(m)}}$ is the generator labeled by the standard $m$-tableau $X_{\lambda'^{(m)}}$ (of the shape ${\lambda'^{(m)}}$). 

\vskip .2cm
Moving now (following the instructions from the Proposition \ref{prop-rel}) in the expressions 
$\phi\mathcal{X}_{\lambda^{(m)}}\mathcal{X}_{\lambda'^{(m)}}\vert\rangle$, where $\phi\in H(m,1,n)$, the elements 
$\mathcal{X}$'s to the left and evaluating, with the help of (\ref{movac}), the elements of the cyclotomic algebra on the vacuum, we define the $H(m,1,n)$-module structure on 
the tensor product $U_{\lambda^{(m)}}\otimes U_{\lambda'^{(m)}}$ of the vector spaces underlying the representations $V_{\lambda^{(m)}}$ and $V_{\lambda'^{(m)}}$.
We denote the resulting representation of $H(m,1,n)$ by $V_{\lambda^{(m)}}\hat{\otimes}V_{\lambda'^{(m)}}$.

\vskip .2cm
In principle, the tensor product $\hat{\otimes}$ is defined for given $m$ and $n$ and should be rather denoted $\hat{\otimes}_{m,n}$. For brevity, we omit $m$ in 
the notation for the product, the value of $m$ is fixed in our treatment. As for $n$, below we introduce the operation of restriction from $n$ to $(n-1)$ and explain 
that the tensor product $\hat{\otimes}$ is compatible with the restriction; due to the compatibility, we omit $n$ in the notation for the product as well.

\vskip .2cm
Notice that for any product $\mathcal{X}_{\lambda^{(m)}}\mathcal{X}_{\lambda'^{(m)}}$, the generator $\tau$ never passes through $\mathcal{X}_{\lambda^{(m)}}$ to the right (due to the particular form (\ref{rel-a2}) of the instructions for the generator $\tau$). Only the generators $\sigma_1,\dots,\sigma_{n-1}$ pass through $\mathcal{X}_{\lambda^{(m)}}$ and then act on $\mathcal{X}_{\lambda'^{(m)}}|\rangle$. The subalgebra of $H(m,1,n)$ generated by $\sigma_1,\dots,\sigma_{n-1}$ is isomorphic to the Hecke algebra $H(1,1,n)$. Thus, it makes sense to consider tensor products  of the form $V_{\lambda^{(m)}}\hat{\otimes}V$, where $V$ is a representation of $H(1,1,n)$, as a representation of $H(m,1,n)$. Moreover, by construction, the representation $V_{\lambda^{(m)}}\hat{\otimes}V_{\lambda^{'(m)}}$ is naturally isomorphic to the representation $V_{\lambda^{(m)}}\hat{\otimes}W$, where $W$ is the restriction of the representation $V_{\lambda^{'(m)}}$ to the subalgebra generated by $\sigma_1,\dots,\sigma_{n-1}$. 

\paragraph{2. Product $\hat{\otimes}\colon$ simplest examples.}
Let $\varpi^{(m)}$ be the $m$-partition $(\lambda_1,\dots,\lambda_m)$ 
of size $n$ such that $\lambda_1$ is the one-row partition $(n)$
and $\lambda_2,\dots,\lambda_m$ are empty partitions. There is only one standard $m$-tableau of 
shape $\varpi^{(m)}$ which we denote by $X_{\varpi^{(m)}}$. For this particular $m$-partition, the formulas (\ref{rel-a1})--(\ref{rel-a2}) become:
\begin{equation}\label{acon0} (\sigma_i-q){\mathcal{X}}_{\varpi^{(m)}}=0\ \ \ \textrm{for all}\ \ i=1,\dots,n-1\ \ \textrm{and}\quad
(\tau-v_1){\mathcal{X}}_{\varpi^{(m)}}=0\ .\end{equation}
Thus the representation $V_{\varpi^{(m)}}$ is isomorphic to the one-dimensional representation of the algebra $H(m,1,n)$ spanned by the vacuum 
$|\rangle$. By construction, the following properties are verified (the isomorphisms are to be understood as isomorphisms of $H(m,1,n)$-modules):
for any $m$-partition $\lambda^{(m)}$ of size $n$,
\begin{equation}\label{teby1right}V_{\lambda^{(m)}}\hat{\otimes} V_{\varpi^{(m)}}\cong V_{\lambda^{(m)}}\ ,\end{equation}
and
\begin{equation}\label{teby1left}V_{\varpi^{(m)}}\hat{\otimes} V_{\lambda^{(m)}}\cong V_{\varpi^{(m)}}\oplus\dots\oplus 
V_{\varpi^{(m)}}\cong \dim(V_{\lambda^{(m)}})\ V_{\varpi^{(m)}}\ . \end{equation}
Actually, in the formulas (\ref{teby1right})--(\ref{teby1left}) one can replace $\varpi^{(m)}$ by any $m$-partition 
$\varpi'^{(m)}$ such that $V_{\varpi'^{(m)}}$ is one-dimensional; these are the $m$-partitions $(\lambda_1,...,\lambda_m)$ with only one 
non-empty $\lambda_j$ which equals $(n)$ or $(1^n)$. For the validity of  (\ref{teby1right}) for $\varpi'^{(m)}$ see, for example, the 
Remark {\bf (b)} after formulas (\ref{rep-a1})--(\ref{rep-a2}); the validity of (\ref{teby1left})  for $\varpi'^{(m)}$ is immediate.

\vskip .2cm 
The formulas  (\ref{teby1right})--(\ref{teby1left}) are obtained in a straightforward manner. The Proposition \ref{tens-prod} below describes the product 
$V_{\lambda^{(m)}}\hat{\otimes}V_{\lambda'^{(m)}}$ of representations corresponding to two arbitrary $m$-partitions of the same size. 
However the proof of the general formula (\ref{tens-prod-form}) is not direct and relies on the completeness assertion which in turn requires the restrictions 
(\ref{sesi1})--(\ref{sesi3}).

\subsection*{A.2  \hspace{.2cm} Decomposition rules for the tensor product $\hat{\otimes}$}
\addcontentsline{toc}{subsection}{A.2$\ $ Decomposition rules for the tensor product $\hat{\otimes}$}

\paragraph{1. Restriction.} An essential role in our study of the product $\hat{\otimes}$ is played by the operation of restriction which allows to use the 
induction arguments.

\vskip .2cm   
For any representation ${\cal{W}}$ of the algebra $H(m,1,n)$ we denote by ${\mathrm{Res}}^n_{n-1}({\cal{W}})$ the 
restriction of ${\cal{W}}$ to the subalgebra of $H(m,1,n)$ generated by $\tau,\sigma_1,\dots,\sigma_{n-2}$; this subalgebra is isomorphic to $H(m,1,n-1)$. This justifies the notation ${\mathrm{Res}}^n_{n-1}$.

\vskip .2cm
The class ${\cal{S}}$ of representations is stable with respect to the restriction 
${\mathrm{Res}}^n_{n-1}$. Indeed, the formulas (\ref{rel-a1})-(\ref{rel-a2}) imply that for any $m$-partition $\lambda^{(m)}$ of size $n$,
\begin{equation}\label{eqres}{\mathrm{Res}}^n_{n-1}(V_{\lambda^{(m)}})\cong\bigoplus_{\alpha^{(m)}\colon\,\alpha^{(m)}\in {\cal{E}}_-(\lambda^{(m)})}\ \  V_{\lambda^{(m)}\backslash\{\alpha^{(m)}\}}\ ,\end{equation}
where, we recall, ${\cal{E}}_-(\lambda^{(m)})$ is the set of removable $m$-nodes of $\lambda^{(m)}$. The stability follows from (\ref{eqres}).

\vskip .2cm
Geometrically, it is clear that one can reconstruct the $m$-partition $\lambda^{(m)}$ from the set ${\cal{E}}_-(\lambda^{(m)})$ of its removable $m$-nodes.
Therefore, by (\ref{eqres}), for an irreducible $H(m,1,n)$-representation $V$, belonging to the class ${\cal{S}}$, $V\cong V_{\lambda^{(m)}}$, its restriction
${\mathrm{Res}}^n_{n-1}(V)$ characterizes the representation $V$ of $H(m,1,n)$ uniquely up to isomorphism.

\vskip .2cm
Moreover, a direct inspection shows that the operation $\hat{\otimes}$ on representations belonging to the class ${\cal{S}}$ is compatible with 
the operation of restriction in the following sense: for any two $H(m,1,n)$-representations ${\cal{W}}$ and ${\cal{W}}'$ belonging to ${\cal{S}}$, 
we have
\begin{equation}\label{eqres2}{\mathrm{Res}}^n_{n-1}({\cal{W}}\hat{\otimes}{\cal{W}}')\,\cong\,\bigl({\mathrm{Res}}^n_{n-1}({\cal{W}})\bigr)\hat{\otimes}\bigl({\mathrm{Res}}^n_{n-1}({\cal{W}}')\bigr)\ .\end{equation}
Here the symbol $\hat{\otimes}$ in the left hand side is the product for $H(m,1,n)$; in the right hand side it is the product for $H(m,1,n-1)$. The formula
(\ref{eqres2}) justifies the usage of the symbol $\hat{\otimes}$ instead of more rigorous $\hat{\otimes}_n$.

\vskip .2cm
Note that (\ref{eqres}) and (\ref{eqres2}) are valid whenever the formula (\ref{rel-a1}) makes sense for the participating representations (that is, the 
denominators in (\ref{rel-a1}) do not vanish); we do not need here the completeness result from Section  \ref{sec-comp}.

\paragraph{2. Decomposition rules.}  
Under the constraints (\ref{sesi1})--(\ref{sesi3}) on the parameters of $H(m,1,n)$, the product $\hat{\otimes}$ of two representations from the class
${\cal{S}}$ belongs again to the class ${\cal{S}}$ due to the completeness result from Section  \ref{sec-comp}. The following Proposition gives the 
decomposition rules for the tensor product $\hat{\otimes}$ of the representations from the class ${\cal{S}}$.

\begin{prop}\label{tens-prod}
{\hspace{-.2cm}.\hspace{.2cm}} Let $\lambda^{(m)}$ and $\lambda'^{(m)}$ be two arbitrary $m$-partitions of size $n$. Assume that the 
conditions (\ref{sesi1})--(\ref{sesi3}) on the parameters of $H(m,1,n)$ hold. Then the representation 
$V_{\lambda^{(m)}}\hat{\otimes}V_{\lambda'^{(m)}}$ of $H(m,1,n)$ is isomorphic to the direct sum of $\dim(V_{\lambda'^{(m)}})$ copies of 
$V_{\lambda^{(m)}}$ :
\begin{equation}\label{tens-prod-form}
V_{\lambda^{(m)}}\hat{\otimes}V_{\lambda'^{(m)}}\cong \dim(V_{\lambda'^{(m)}})\ V_{\lambda^{(m)}}.
\end{equation}
\end{prop}
\emph{Proof.} We shall use here that two representations $V_{\pi_1^{(m)}}$ and $V_{\pi_2^{(m)}}$ of $H(m,1,n)$ corresponding to two different 
$m$-partitions $\pi_1^{(m)}$ and $\pi_2^{(m)}$ of $n$ are non-isomorphic and that any representation of $H(m,1,n)$ belongs to the class 
${\cal{S}}$ 
(this is proved, under the constraints (\ref{sesi1})--(\ref{sesi3}), in Section  \ref{sec-comp}).
 
\vskip .2cm
We shall need the following Lemma. 

\begin{lemm}\label{lemm-tens-prod}
{\hspace{-.2cm}.\hspace{.2cm}}
(i) Let $\lambda^{(m)}$ be an $m$-partition of size $n$ satisfying the following conditions:
\begin{itemize}
\item $\lambda^{(m)}$ is different from $m$-partitions of the form $(\varnothing,\dots,\varnothing,\lambda,\varnothing,\dots,\varnothing)$ where $\lambda$ is a partition such 
that 
$|\lambda|\leq 2$ or $\lambda=(2,1)$ or $\lambda=(2,1,1)$ or $\lambda=(3,1)$; 
\item $\lambda^{(m)}$ is different from $m$-partitions of the form $(\varnothing,\dots,\varnothing,\Box,\varnothing,\dots,\varnothing,\Box,\varnothing,\dots,\varnothing)$.
\end{itemize}
(ii) Let $\mathfrak{L}=\{\lambda_1^{(m)},\dots,\lambda_l^{(m)}\}$ be any $l$-tuple
of $m$-partitions of size $n$ different from $\lambda^{(m)}$, $\lambda_j^{(m)}\neq\lambda^{(m)}$ for all $j=1,...,l$. 

Then the two following sets of $m$-partitions 
$$\lambda^{(m)-}:=\bigl\{ \lambda^{(m)}\backslash\{\alpha^{(m)}\}\bigr\}_{\alpha^{(m)}\colon \alpha^{(m)}\in{\cal{E}}_-(\lambda^{(m)})}\ \ \text{and}\ \ 
\mathfrak{L}^-:=\bigl\{ \lambda_j^{(m)}\backslash\{\alpha_j^{(m)}\}\bigr\}_{j\colon j=1,\dots ,l\, ;\ \alpha_j^{(m)}\colon  
\alpha_j^{(m)}\in{\cal{E}}_-(\lambda_j^{(m)})}$$ 
do not coincide.
\end{lemm}

The Lemma \ref{lemm-tens-prod} is a purely combinatorial statement. Its proof will be combinatorial as well.
We shall prove that, if $\lambda^{(m)-}$ is contained in $\mathfrak{L}^-$
then there exists a sub-$m$-partition of size $n-1$ of one of $\lambda_j^{(m)}\in\mathfrak{L}$ which is not a sub-$m$-partition of $\lambda^{(m)}$. The representation-theoretic translation of the last sentence is the following. 

\begin{coro}\label{cor-lemm-tens-prod}
{\hspace{-.2cm}.\hspace{.2cm}}
Under the conditions (i) and (ii), if ${\mathrm{Res}}^n_{n-1}(V_{\lambda^{(m)}})$ is isomorphic to a sub-represen\-ta\-tion of ${\mathrm{Res}}^n_{n-1} (V_{\lambda_1^{(m)}}\oplus\dots\oplus V_{\lambda_l^{(m)}})$, then there exists an $m$-partition $\nu^{(m)}$ of size $n-1$ such that $V_{\nu^{(m)}}$ is  isomorphic to a sub-represen\-ta\-tion of 
${\mathrm{Res}}^n_{n-1}(V_{\lambda_1^{(m)}}\oplus\dots\oplus V_{\lambda_l^{(m)}})$
but not to a sub-representation of ${\mathrm{Res}}^n_{n-1}(V_{\lambda^{(m)}})$. 
\end{coro}

\emph{Proof of the Lemma.} Assume that the $l$-tuple $\mathfrak{L}=\{\lambda_1^{(m)},\dots,\lambda_l^{(m)}\}$, formed by $m$-partitions different from 
$\lambda^{(m)}$, is such that $\lambda^{(m)-}$ is contained in $\mathfrak{L}^-$. 
Then for each $\alpha^{(m)}\in {\cal{E}}_-(\lambda^{(m)})$ 
the $m$-partition $\lambda^{(m)}\backslash\{\alpha^{(m)}\}$ of $n-1$ is a sub-$m$-partition of some 
$m$-partition from $\mathfrak{L}$; 
for a given $\alpha^{(m)}\in {\cal{E}}_-(\lambda^{(m)})$,
there may be several $m$-partitions from $\mathfrak{L}$ with this property; choose one of them and denote it by 
$\mu^{(m)}$. By condition $(ii)$, 
the $m$-partition $\mu^{(m)}$ is obtained from $\lambda^{(m)}\backslash\{\alpha^{(m)}\}$ by adding an 
$m$-node $\beta^{(m)}\in{\cal{E}}_+(\lambda^{(m)}\backslash\{\alpha^{(m)}\})$ different from $\alpha^{(m)}$. 

\vskip .2cm
For a removable $m$-node $\gamma^{(m)}$ of 
$\mu^{(m)}$, 
it is geometrically clear that 
the $m$-partition $\mu^{(m)}\backslash\{\gamma^{(m)}\}$ is not a sub-$m$-partition of $\lambda^{(m)}$ 
if and only if $\gamma^{(m)}$ is different from $\beta^{(m)}$; 
indeed, 
as for ordinary partitions, the inclusion graph (see Appendix B for the definition) of $m$-partitions is a lattice; in particular, it cannot contain a
subgraph of the form 

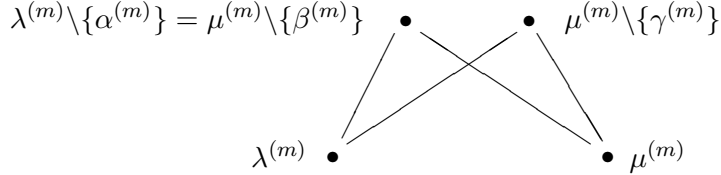
\setlength{\unitlength}{1600sp}
\begin{picture}(0,3500)(-4150,-6600)
\thinlines
{\put(3900,-4086){\line(-1,-2){810}}}{\put(6081,-4086){\line( 3,-5){1000}}}{\put(3170,-5800){\line(3,2){2600}}}{\put(6900,-5800){\line(-3,2){2600}}}
\put(4070,-3900){\circle*{144}}\put(5960,-3900){\circle*{144}}\put(2950,-6000){\circle*{144}}\put(7170,-6000){\circle*{144}}
\put(-2000,-4000){$\lambda^{(m)}\backslash\{\alpha^{(m)}\}=\mu^{(m)}\backslash\{\beta^{(m)}\}$}
\put(6500,-4000){$\mu^{(m)}\backslash\{\gamma^{(m)}\}$}
\put(1700,-6150){$\lambda^{(m)}$}\put(7500,-6150){$\mu^{(m)}$}
\end{picture}

So, if  $\mu^{(m)}$ has a removable $m$-node $\gamma^{(m)}$ different from $\beta^{(m)}$ then 
$\mu^{(m)}\backslash\{\gamma^{(m)}\}$ is a sub-$m$-partition of one of the $\lambda_j^{(m)}\in\mathfrak{L}$ and is not a sub-$m$-partition of $\lambda^{(m)}$. 
Therefore it suffices to show that
there is an $m$-node $\alpha^{(m)}\in {\cal{E}}_-(\lambda^{(m)})$ such that, for any 
$m$-partition $\mu^{(m)}$ obtained by adding to $\lambda^{(m)}\backslash\{\alpha^{(m)}\}$ an $m$-node 
$\beta^{(m)}\in  {\cal{E}}_+(\lambda^{(m)}\backslash\{\alpha^{(m)}\})$, different from $\alpha^{(m)}$, 
there is a removable $m$-node $\gamma^{(m)}\in{\cal{E}}_-(\mu^{(m)})$, different from $\beta^{(m)}$.

\vskip .2cm
Such $\gamma^{(m)}$ does not exist if and only if the $m$-partition $\mu^{(m)}$ has only one removable $m$-node, that is, 
\begin{equation}\label{condilamu}\begin{array}{l}
\text{only one partition in the $m$-tuple $\mu^{(m)}$ is non-empty}\\[.1em] 
\text{and this partition is of rectangular shape.}\end{array}\end{equation} 

Thus, an $m$-partition $\lambda^{(m)}$ contradicting to the assertion of the  Lemma \ref{lemm-tens-prod} should verify the property:
for every $\alpha^{(m)}\in  {\cal{E}}_-(\lambda^{(m)})$ there exists $\beta^{(m)}\in  {\cal{E}}_+(\lambda^{(m)}\backslash\{\alpha^{(m)}\})$, different from 
$\alpha^{(m)}$,  such that the $m$-partition $\mu^{(m)}$, obtained by adding to $\lambda^{(m)}\backslash\{\alpha^{(m)}\}$ 
the $m$-node 
$\beta^{(m)}$, is described by (\ref{condilamu}). A direct inspection shows that such $m$-partitions are exactly those which are excluded by the 
part $(i)$ of the formulation of the Lemma \ref{lemm-tens-prod}. \hfill$\square$

\vskip .2cm
We now return to the proof of the Proposition \ref{tens-prod}. 

\vskip .2cm
We proceed by induction on $n$; the formula (\ref{tens-prod-form}) is trivial for $n=0$, that is, for $\lambda^{(m)}=(\varnothing,\dots,\varnothing)$. 
As the proof of the induction step below shows,
we shall need to verify separately several cases to complete the proof of the Proposition \ref{tens-prod}. 

\paragraph{Induction step.} Let $\lambda^{(m)}$ and $\lambda'^{(m)}$ be two $m$-partitions of size $n$ 
such that $\lambda^{(m)}$ satisfies the conditions from the part (i) of the Lemma \ref{lemm-tens-prod}. 
Due to the formulas (\ref{eqres})--(\ref{eqres2}),  
\begin{equation}\label{vsvsvs}{\mathrm{Res}}^n_{n-1}(V_{\lambda^{(m)}}\hat{\otimes}V_{\lambda^{'(m)}})\,\cong\,\bigoplus_{\alpha^{(m)},\,\alpha^{'(m)}\ 
\colon\, \textrm{
{$\!\!\begin{array}{l}\\[-1.25em] \alpha^{(m)}\in{\cal{E}}_-(\lambda^{(m)})\\[0.01em]\alpha^{'(m)}\in{\cal{E}}_-(\lambda^{'(m)})\end{array}$}}}\ 
V_{\lambda^{(m)}\backslash\{\alpha^{(m)}\}}\hat{\otimes}V_{\lambda^{'(m)}\backslash\{\alpha^{'(m)}\}}.\end{equation}
Our induction hypothesis is: the formula (\ref{tens-prod-form}) is valid for the products in the right hand side of (\ref{vsvsvs}).
Assuming the induction hypothesis, we transform the right hand side of (\ref{vsvsvs}),
\begin{equation}\label{eqres3}{\mathrm{Res}}^n_{n-1}(V_{\lambda^{(m)}}\hat{\otimes}V_{\lambda^{'(m)}})
\cong \dim(V_{\lambda^{'(m)}})\ {\mathrm{Res}}^n_{n-1}(V_{\lambda^{(m)}}).\end{equation}

Now we shall employ several times the completeness result from Section \ref{sec-comp}. First of all, the representation 
$V_{\lambda^{(m)}}\hat{\otimes}V_{\lambda^{'(m)}}$ belongs to ${\cal{S}}$
so we can write $V_{\lambda^{(m)}}\hat{\otimes}V_{\lambda^{'(m)}}=c\ V_{\lambda^{(m)}}\oplus W$ where 
\begin{itemize}
\item{} $W$ belongs to ${\cal{S}}$ and has no irreducible constituents isomorphic 
to $V_{\lambda^{(m)}}$, 
\item{} $c$ is a non-negative integer; $c\leq \dim(V_{\lambda^{'(m)}})$ by the dimension argument. 
\end{itemize}

If $c<\dim(V_{\lambda^{'(m)}})$ we 
use again results of Section \ref{sec-comp} to obtain a contradiction. Due to the semi-simplicity at level $n-1$ (note that the conditions (\ref{sesi1})--(\ref{sesi3}) 
at level $n-1$ are implied by the conditions (\ref{sesi1})--(\ref{sesi3}) at level $n$), the representation monoid of $H(m,1,n-1)$ is cancellative; we can thus 
simplify (\ref{eqres3}) by $c\, {\mathrm{Res}}^n_{n-1}(V_{\lambda^{(m)}})$ on both sides.
We obtain: \[{\mathrm{Res}}^n_{n-1}(W)\cong \Bigl(\dim(V_{\lambda^{'(m)}})-c\Bigr)\ {\mathrm{Res}}^n_{n-1}(V_{\lambda^{(m)}})\ .\] 
Since ${\mathrm{Res}}^n_{n-1}(V_{\lambda^{(m)}})$ is isomorphic to a sub-representation of 
${\mathrm{Res}}^n_{n-1}(W)$, the Corollary \ref{cor-lemm-tens-prod}
implies the existence of an $m$-partition $\nu^{(m)}$ of $n-1$ such that $\nu^{(m)}$ is not a 
sub-$m$-partition of $\lambda^{(m)}$ but $V_{\nu^{(m)}}$ is isomorphic to a sub-representation of ${\mathrm{Res}}^n_{n-1}(W)$. But now the representation 
${\mathrm{Res}}^n_{n-1}(W)$ is isomorphic to a direct sum of several copies of ${\mathrm{Res}}^n_{n-1}(V_{\lambda^{(m)}})$ implying that such 
$\nu^{(m)}$ cannot exist, a contradiction. Thus $c=\dim(V_{\lambda^{'(m)}})$ and
\[V_{\lambda^{(m)}}\hat{\otimes}V_{\lambda^{'(m)}}\cong \dim(V_{\lambda^{'(m)}})\ V_{\lambda^{(m)}}.\]
The proof of the induction step is finished.

\paragraph{End of the proof of the Proposition \ref{tens-prod}.}  If the formula (\ref{tens-prod-form}) is verified for an $m$-partition $\lambda^{(m)}$ and an arbitrary $m$-partition $\lambda^{'(m)}$ 
of the same size 
as $\lambda^{(m)}$, we shall simply say that  (\ref{tens-prod-form}) is verified for $\lambda^{(m)}$.
Our way of proof of the induction step referred to the Lemma \ref{lemm-tens-prod}. So if the formula (\ref{tens-prod-form}) is established for all
$m$-partitions of size $n$ then it is established for all $m$-partitions $\lambda^{(m)}$ of size $(n+1)$ unless  $\lambda^{(m)}$ belongs 
to the set of $m$-partitions 
listed in part (i) of the Lemma \ref{lemm-tens-prod}. For the $m$-partitions  listed in part (i) of the Lemma \ref{lemm-tens-prod} an independent
proof is needed. Besides, by (\ref{teby1left}),
the formula (\ref{tens-prod-form}) is already 
established for $m$-partitions $\lambda^{(m)}$ of 1, that is, for  $\lambda^{(m)}=(\varnothing,\dots,\varnothing,\Box,\varnothing,\dots,\varnothing)$ and for $m$-partitions 
of the form
\begin{equation}\label{inba0}(\varnothing,\dots,\varnothing,\lambda,\varnothing,\dots,\varnothing)\, ,
\ \, \text{where}\ \lambda\ \text{is}\,\; (2)\, \ \text{or}\ (1,1)\ .\end{equation}
Below we shall separately verify that (\ref{tens-prod-form}) holds 
for other $m$-partitions $\lambda^{(m)}$ listed in the part $(i)$ of the Lemma \ref{lemm-tens-prod}, that is,  $m$-partitions $\lambda^{(m)}$ of the form
\begin{equation}\label{inba1} (\varnothing ,\dots ,\varnothing ,\lambda ,\varnothing ,\dots ,\varnothing )\, ,
\ \, \text{where}\ \lambda\ \text{is}\,\; (2,1)\, ,\ (2,1,1)\,\ \text{or}\ (3,1)\end{equation}
or of the form
\begin{equation}\label{inba2}(\varnothing,\dots,\varnothing,\Box,\varnothing,\dots,\varnothing,\Box,\varnothing,\dots,\varnothing)\ .\end{equation}

\paragraph{Proof of (\ref{tens-prod-form})
for the $m$-partitions of the forms (\ref{inba1}) and (\ref{inba2}).}$\ $

\vskip .2cm
We recall here that, for any $\lambda^{(m)}$, $\lambda^{'(m)}$, the representation $V_{\lambda^{(m)}}\hat{\otimes}V_{\lambda^{'(m)}}$ is naturally isomorphic to a representation $V_{\lambda^{(m)}}\hat{\otimes}W$ where $W$ is a representation of the Hecke algebra $H(1,1,n)$ (see the end of the  paragraph {\bf 1} of Subsection A.1). Due to the completeness result of Section \ref{sec-comp}, it is enough to consider the cases $W\cong V_{\lambda}$ for all partitions $\lambda$.

\vskip .2cm
For the $m$-partitions $\lambda^{(m)}$ of the form (\ref{inba1}), the generator $\tau$ acts by a constant in $V_{\lambda^{(m)}}$ and it is thus sufficient to work with the Hecke algebra $H(1,1,n)$.

\vskip .2cm
Reintroduce for this paragraph the deformation parameter $q$ in the notation for the Hecke algebra: $H_q(1,1,n)$. Due to the relations (\ref{def1a})--(\ref{def1b}) and (\ref{def1''b}), we have an isomorphism $\theta\colon H_q(1,1,n)\to H_{-q^{-1}}(1,1,n)$ of algebras, defined on 
generators by $H_q(1,1,n)\ni\sigma_i\mapsto\sigma_i\in H_{-q^{-1}}(1,1,n)$. The composition  with $\theta$  of a representation of  $H_{-q^{-1}}(1,1,n)$, 
corresponding to a partition $\lambda$,
sends the representations $V_{(3,1)}$, $V_{(2,2)}$ and $V_{(2,1,1)}$ of $H_{-q^{-1}}(1,1,n)$ to, respectively, the representations $V_{(2,1,1)}$, 
$V_{(2,2)}$ and $V_{(3,1)}$ of $H_q(1,1,n)$. Thus the formula (\ref{tens-prod-form}) for $\lambda =(3,1)$ follows from 
the formula (\ref{tens-prod-form}) for $\lambda =(2,1,1)$.

\vskip .2cm
We remind that (\ref{tens-prod-form}) has already been proved for any $m$-partitions $\lambda'^{(m)}$ such that $V_{\lambda'^{(m)}}$ is one-dimensional, see (\ref{teby1left}).

\vskip .2cm
For the $m$-partitions $\lambda^{(m)}$ of the form (\ref{inba2}), the proof of (\ref{tens-prod-form}) is reduced to the situation where $V_{\lambda'^{(m)}}$ is replaced by $V_{\lambda'}$ where $\lambda'$ is $(2)$ or $(1,1)$ in which case the representation $V_{\lambda'}$ is one-dimensional and the formula (\ref{tens-prod-form}) follows.

\vskip .3cm
{\bf 1.} For the $m$-partition $\lambda^{(m)}$ of the form (\ref{inba1}) with $\lambda =(2,1)$ we have reduced the proof to the situation $m=1$, and it remains to establish the result (\ref{tens-prod-form}) only for 
$V_{\lambda}\hat{\otimes}V_{\lambda}$.

\vskip .2cm
The underlying vector space of $V_{\lambda}$ has a basis $\{{\cal{X}}_1,{\cal{X}}_2\}$ where 
${\cal{X}}_1:={\cal{X}}_{\textrm{\tiny{$\!\!
\begin{array}{c}\fbox{\scriptsize{$1$}}\\[-0.10em] 
\fbox{\scriptsize{$3$}}\end{array}\hspace{-.348cm} \begin{array}{c}\fbox{\scriptsize{$2$}}\\[-0.10em] 
\phantom{\fbox{\scriptsize{$4$}}}
\end{array}$}}}$ and ${\cal{X}}_2:={\cal{X}}_{\textrm{\tiny{$\!\!
\begin{array}{c}\fbox{\scriptsize{$1$}}\\[-0.10em] 
\fbox{\scriptsize{$2$}}\end{array}\hspace{-.348cm} \begin{array}{c}\fbox{\scriptsize{$3$}}\\[-0.10em] 
\phantom{\fbox{\scriptsize{$4$}}}
\end{array}$}}}$. 
In this basis, the generators $\sigma_1$ and $\sigma_2$ are realized as follows:
\begin{equation}\label{rep21}\sigma_1\mapsto\text{diag}(q,-q^{-1})\ ,\ \sigma_2\mapsto\frac{1}{2_q}\left(\begin{array}{cc}
-q^{-2}&3_q\\1&q^2\end{array}\right)\ .\end{equation}

We order the basis ${\cal{X}}_i{\cal{X}}_j$ of the underlying vector space of $V_{\lambda}\hat{\otimes}V_{\lambda}$ lexicographically; that is, we choose 
the order $\{ {\cal{X}}_1{\cal{X}}_1,{\cal{X}}_1{\cal{X}}_2,{\cal{X}}_2{\cal{X}}_1,{\cal{X}}_2{\cal{X}}_2\}$. In this basis, the matrices of the generators $\sigma_1$ and $\sigma_2$ are:
$$\sigma_1\mapsto\text{diag}(q,q,-q^{-1},-q^{-1})\ ,\ \sigma_2\mapsto\frac{1}{2_q}\left(\begin{array}{cccc}
-q^{-2}&0&0&3_q\\ 0&-q^{-2}&1&q^2+q^{-2}\\ -q^2-q^{-2}&3_q&q^2&0\\ 1&0&0&q^2\end{array}\right)\ .$$
The two subspaces with the bases 
\begin{equation}\label{tsuspa}\begin{array}{c}\{ {\cal{X}}_1{\cal{X}}_2,3_q {\cal{X}}_2{\cal{X}}_1\}\ ,\\[1em]
\{ {\cal{X}}_1{\cal{X}}_1+{\cal{X}}_1{\cal{X}}_2, {\cal{X}}_2{\cal{X}}_1+{\cal{X}}_2{\cal{X}}_2\}\end{array}\end{equation}
carry the representation (\ref{rep21}) which implies (\ref{tens-prod-form}) in this case. Since $V_{\lambda}\hat{\otimes}V_{\lambda}$
decomposes into a direct sum of two isomorphic representations, the 
choice (\ref{tsuspa}) of subspaces is not unique. We just make (here and for the other cases below) a simple choice.

\vskip .3cm
{\bf 2.} For the $m$-partition $\lambda^{(m)}$ of the form (\ref{inba1}) with $\lambda =(2,1,1)$ we have reduced the proof to the situation $m=1$, and it remains to establish the result (\ref{tens-prod-form}) only for 
$V_{\lambda}\hat{\otimes}V_{\lambda'}$ with $\lambda'=(2,2), (2,1,1)$ and $(3,1)$.

\vskip .2cm
The underlying vector space of $V_{\lambda}$ has a basis $\{{\cal{X}}_1,{\cal{X}}_2,{\cal{X}}_3\}$ where
${\cal{X}}_1:={\cal{X}}_{\textrm{\tiny{$\!\!
\begin{array}{c}\fbox{\scriptsize{$1$}}\\[-0.10em] 
\fbox{\scriptsize{$3$}}\\[-0.10em]
\fbox{\scriptsize{$4$}}
\end{array}\hspace{-.348cm} \begin{array}{c}\fbox{\scriptsize{$2$}}\\[-0.10em] 
\phantom{\fbox{\scriptsize{$4$}}}\\[-0.10em]\phantom{\fbox{\scriptsize{$4$}}}
\end{array}$}}}$, 
${\cal{X}}_2:={\cal{X}}_{\textrm{\tiny{$\!\!
\begin{array}{c}\fbox{\scriptsize{$1$}}\\[-0.10em] 
\fbox{\scriptsize{$2$}}\\[-0.10em]
\fbox{\scriptsize{$4$}}
\end{array}\hspace{-.348cm} \begin{array}{c}\fbox{\scriptsize{$3$}}\\[-0.10em] 
\phantom{\fbox{\scriptsize{$4$}}}\\[-0.10em]\phantom{\fbox{\scriptsize{$4$}}}
\end{array}$}}}$ and  
${\cal{X}}_3:={\cal{X}}_{\textrm{\tiny{$\!\!
\begin{array}{c}\fbox{\scriptsize{$1$}}\\[-0.10em] 
\fbox{\scriptsize{$2$}}\\[-0.10em]
\fbox{\scriptsize{$3$}}
\end{array}\hspace{-.348cm} \begin{array}{c}\fbox{\scriptsize{$4$}}\\[-0.10em] 
\phantom{\fbox{\scriptsize{$4$}}}\\[-0.10em]\phantom{\fbox{\scriptsize{$4$}}}
\end{array}$}}}$. In this basis, the generators $\sigma_1$, $\sigma_2$ and $\sigma_3$ are realized as follows: 
\begin{equation}\label{rep211}\begin{array}{c}
\sigma_1\mapsto\text{diag}(q,-q^{-1},-q^{-1})\ ,\\[1.5em]
\sigma_2\mapsto\frac{1}{2_q}\left(\begin{array}{ccc}
-q^{-2}&3_q&0\\1&q^2&0\\0&0&-q^{-1}2_q
\end{array}\right)\ ,\ \sigma_3\mapsto\frac{1}{3_q}\left(\begin{array}{ccc}-q^{-1}3_q&0&0\\0&-q^{-3}&4_q\\0&2_q&q^3
\end{array}\right)\ .
\end{array}\end{equation}

\vskip .2cm
{\bf 2a.} $\lambda'=(2,2)$. 

\vskip .2cm
The underlying vector space of $V_{\lambda'}$ has a basis $\{{\cal{Y}}_1,{\cal{Y}}_2\}$, where  ${\cal{Y}}_1:={\cal{X}}_{\textrm{\tiny{$\!\!
\begin{array}{c}\fbox{\scriptsize{$1$}}\\[-0.10em] 
\fbox{\scriptsize{$3$}}\end{array}\hspace{-.348cm} \begin{array}{c}\fbox{\scriptsize{$2$}}\\[-0.10em] 
\fbox{\scriptsize{$4$}}\end{array}$}}}$ and
 ${\cal{Y}}_2:={\cal{X}}_{\textrm{\tiny{$\!\!
\begin{array}{c}\fbox{\scriptsize{$1$}}\\[-0.10em] 
\fbox{\scriptsize{$2$}}\end{array}\hspace{-.348cm} \begin{array}{c}\fbox{\scriptsize{$3$}}\\[-0.10em] 
\fbox{\scriptsize{$4$}}\end{array}$}}}$.
We order the basis ${\cal{X}}_i{\cal{Y}}_j$ of the underlying vector space of $V_{\lambda}\hat{\otimes}V_{\lambda'}$ lexicographically. In this basis the generators $\sigma_1$, $\sigma_2$ and $\sigma_3$ are realized as follows:
$$\sigma_1\mapsto\text{diag}(q,q,-q^{-1},-q^{-1},-q^{-1},-q^{-1})\ ,$$
\vskip .2cm
$$ \sigma_2\mapsto\frac{1}{2_q}\left(\begin{array}{cccccc}
-q^{-2}&0&0&3_q&0&0\\ 0&-q^{-2}&1&q^2+q^{-2}&0&0\\ -q^2-q^{-2}&3_q&q^2&0&0&0\\ 1&0&0&q^2&0&0\\
0&0&0&0&-q^{-1}2_q&0\\0&0&0&0&0&-q^{-1}2_q\end{array}\right)\ ,$$
\vskip .4cm
$$ \sigma_3\mapsto\frac{1}{3_q}\left(\begin{array}{cccccc}
-q^{-1}3_q&0&0&0&0&0\\ 0&-q^{-1}3_q&0&0&0&0\\ 0&0&-q^{-3}&0&4_q&0\\ 0&0&0&-q^{-3}&0&-2_q\\
0&0&2_q&0&q^{3}&0\\0&0&0&-4_q&0&q^{3}\end{array}\right)\ ,$$
The two subspaces with the bases 
$$\{ {\cal{X}}_1{\cal{Y}}_2, 3_q{\cal{X}}_2{\cal{Y}}_1,3_q{\cal{X}}_3{\cal{Y}}_1\}\ ,$$ 
$$\{ {\cal{X}}_1{\cal{Y}}_1+{\cal{X}}_1{\cal{Y}}_2, {\cal{X}}_2{\cal{Y}}_1+{\cal{X}}_2{\cal{Y}}_2,{\cal{X}}_3{\cal{Y}}_1-
(q^2+q^{-2}){\cal{X}}_3{\cal{Y}}_2\}$$ 
carry the representation (\ref{rep211}), which implies (\ref{tens-prod-form}) in this case.

\vskip .3cm
{\bf 2b.} $\lambda'=\lambda=(2,1,1)$. 

\vskip .2cm
We order the basis ${\cal{X}}_i{\cal{X}}_j$ of the underlying vector space of $V_{\lambda}\hat{\otimes}V_{\lambda}$ lexicographically. In this basis the generators $\sigma_1$, $\sigma_2$ and $\sigma_3$ are realized as follows:
$$\sigma_1\mapsto\text{diag}(q,q,q,-q^{-1},-q^{-1},-q^{-1},-q^{-1},-q^{-1},-q^{-1})\ ,$$
\vskip .2cm
$$ \sigma_2\mapsto\frac{1}{2_q}\left(\begin{array}{ccccccccc}
-q^{-2}&0&0&0&3_q&0&0&0&0\\ 0&-q^{-2}&0&1&q^2+q^{-2}&0&0&0&0\\0&0&-q^{-2}&0&0&-1&0&0&0 \\
-q^2-q^{-2}&3_q&0&q^2&0&0&0&0&0\\ 1&0&0&0&q^2&0&0&0&0\\0&0&-3_q&0&0&q^2&0&0&0
\\0&0&0&0&0&0&-q^{-1}2_q&0&0\\0&0&0&0&0&0&0&-q^{-1}2_q&0\\0&0&0&0&0&0&0&0&-q^{-1}2_q\end{array}\right)\ ,$$
\vskip .4cm
$$ \sigma_3\mapsto\frac{1}{3_q}\left(\begin{array}{ccccccccc}
-q^{-1}3_q&0&0&0&0&0&0&0&0\\ 0&-q^{-1}3_q&0&0&0&0&0&0&0\\0&0&-q^{-1}3_q&0&0&0&0&0&0
\\ 0&0&0&-q^{-3}&0&0&-2_q&0&0\\ 0&0&0&0&-q^{-3}&0&0&0&4_q\\
0&0&0&0&0&-q^{-3}&0&2_q&q^3+q^{-3}\\0&0&0&-4_q&0&0&q^{3}&0&0\\0&0&0&0&-q^3-q^{-3}&4_q&0&q^3&0\\
0&0&0&0&2_q&0&0&0&q^3
\end{array}\right)\ ,$$
The three subspaces with the bases 
$$\{ {\cal{X}}_1{\cal{Y}}_2, 3_q{\cal{X}}_2{\cal{Y}}_1,-3_q(q^2+q^{-2}){\cal{X}}_3{\cal{Y}}_1\}\ ,$$  
$$\{ {\cal{X}}_1{\cal{Y}}_1+{\cal{X}}_1{\cal{Y}}_2, {\cal{X}}_2{\cal{Y}}_1+{\cal{X}}_2{\cal{Y}}_2,-(q^2+q^{-2})
{\cal{X}}_3{\cal{Y}}_1-
\displaystyle{\frac{q^3+q^{-3}}{2_q}}{\cal{X}}_3{\cal{Y}}_2+{\cal{X}}_3{\cal{Y}}_3\}\ ,$$ 
$$\{ {\cal{X}}_1{\cal{Y}}_3,-3_q{\cal{X}}_2{\cal{Y}}_3,-3_q(q^2+q^{-2}){\cal{X}}_3{\cal{Y}}_2\}$$ 
carry the representation (\ref{rep211}), which implies (\ref{tens-prod-form}) in this case.

\vskip .3cm
{\bf 2c.} $\lambda'=(3,1)$. 

\vskip .2cm
The underlying vector space of $V_{\lambda'}$ has a basis $\{{\cal{Y}}_1,{\cal{Y}}_2,{\cal{Y}}_3\}$, where ${\cal{Y}}_1:={\cal{X}}_
{\textrm{\tiny{$\!\!
\begin{array}{c}\fbox{\scriptsize{$1$}}\\[-0.10em] 
\fbox{\scriptsize{$2$}}\end{array}
\hspace{-.348cm} \begin{array}{c}\fbox{\scriptsize{$3$}}\\[-0.10em] 
\phantom{\fbox{\scriptsize{$4$}}}
\end{array}
\hspace{-.348cm} \begin{array}{c}\fbox{\scriptsize{$4$}}\\[-0.10em] 
\phantom{\fbox{\scriptsize{$4$}}}
\end{array}$}}}$,
${\cal{Y}}_2:={\cal{X}}_
{\textrm{\tiny{$\!\!
\begin{array}{c}\fbox{\scriptsize{$1$}}\\[-0.10em] 
\fbox{\scriptsize{$3$}}\end{array}
\hspace{-.348cm} \begin{array}{c}\fbox{\scriptsize{$2$}}\\[-0.10em] 
\phantom{\fbox{\scriptsize{$4$}}}
\end{array}
\hspace{-.348cm} \begin{array}{c}\fbox{\scriptsize{$4$}}\\[-0.10em] 
\phantom{\fbox{\scriptsize{$4$}}}
\end{array}$}}}$ and
${\cal{Y}}_3:={\cal{X}}_
{\textrm{\tiny{$\!\!
\begin{array}{c}\fbox{\scriptsize{$1$}}\\[-0.10em] 
\fbox{\scriptsize{$4$}}\end{array}
\hspace{-.348cm} \begin{array}{c}\fbox{\scriptsize{$2$}}\\[-0.10em] 
\phantom{\fbox{\scriptsize{$4$}}}
\end{array}
\hspace{-.348cm} \begin{array}{c}\fbox{\scriptsize{$3$}}\\[-0.10em] 
\phantom{\fbox{\scriptsize{$4$}}}
\end{array}$}}}$.
We order the basis ${\cal{X}}_i{\cal{Y}}_j$ of the underlying vector space of $V_{\lambda}\hat{\otimes}V_{\lambda'}$ lexicographically. In this basis the generators $\sigma_1$, $\sigma_2$ and $\sigma_3$ are realized as follows:
$$\sigma_1\mapsto\text{diag}(q,q,q,-q^{-1},-q^{-1},-q^{-1},-q^{-1},-q^{-1},-q^{-1})\ ,$$
\vskip .2cm
$$ \sigma_2\mapsto\frac{1}{2_q}\left(\begin{array}{ccccccccc}
-q^{-2}&0&0&q^2+q^{-2}&1&0&0&0&0\\ 0&-q^{-2}&0&3_q&0&0&0&0&0\\0&0&-q^{-2}&0&0&3_q&0&0&0 \\
0&1&0&q^2&0&0&0&0&0\\ 3_q&-q^2-q^{-2}&0&0&q^2&0&0&0&0\\0&0&1&0&0&q^2&0&0&0
\\0&0&0&0&0&0&-q^{-1}2_q&0&0\\0&0&0&0&0&0&0&-q^{-1}2_q&0\\0&0&0&0&0&0&0&0&-q^{-1}2_q\end{array}\right)\ ,$$
\vskip .4cm
$$ \sigma_3\mapsto\frac{1}{3_q}\left(\begin{array}{ccccccccc}
-q^{-1}3_q&0&0&0&0&0&0&0&0\\ 0&-q^{-1}3_q&0&0&0&0&0&0&0\\0&0&-q^{-1}3_q&0&0&0&0&0&0
\\ 0&0&0&-q^{-3}&0&0&4_q&0&0\\ 0&0&0&0&-q^{-3}&0&0&q^3+q^{-3}&2_q\\
0&0&0&0&0&-q^{-3}&0&4_q&0\\0&0&0&2_q&0&0&q^{3}&0&0\\0&0&0&0&0&2_q&0&q^3&0\\
0&0&0&0&4_q&-q^3-q^{-3}&0&0&q^3
\end{array}\right)\ ,$$
The three subspaces with the bases 
$$\{ {\cal{X}}_1{\cal{Y}}_1, 3_q{\cal{X}}_2{\cal{Y}}_2,3_q(q^2+q^{-2}){\cal{X}}_3{\cal{Y}}_3\}\ ,$$  
$$\{ {\cal{X}}_1{\cal{Y}}_1+{\cal{X}}_1{\cal{Y}}_2, {\cal{X}}_2{\cal{Y}}_1+{\cal{X}}_2{\cal{Y}}_2,
{\cal{X}}_3{\cal{Y}}_1+
(q^2+q^{-2}){\cal{X}}_3{\cal{Y}}_3\}\ ,$$ 
$$\{ {\cal{X}}_1{\cal{Y}}_3,{\cal{X}}_2{\cal{Y}}_3,{\cal{X}}_3{\cal{Y}}_2-\frac{q^3+q^{-3}}{2_q}{\cal{X}}_3{\cal{Y}}_3\}$$ 
carry the representation (\ref{rep211}), which implies (\ref{tens-prod-form}) in this case.

\vskip .2cm
The proof of the Proposition \ref{tens-prod} is completed.\hfill$\square$

\paragraph{Remarks.} $\ $ 

\vskip .2cm
 {\bf (a)} By the Proposition \ref{tens-prod}, the  tensor product $\hat{\otimes}$ is obviously associative. We remark that if we had an independent proof 
 of the associativity of $\hat{\otimes}$ it would immediately imply the Proposition \ref{tens-prod}:  
 \[ V_{\lambda^{(m)}}\hat{\otimes}V_{\lambda^{'(m)}}\cong (V_{\lambda^{(m)}}\hat{\otimes}V_{\varpi^{(m)}})\hat{\otimes}V_{\lambda^{'(m)}}
 \cong V_{\lambda^{(m)}}\hat{\otimes}(V_{\varpi^{(m)}}\hat{\otimes}V_{\lambda^{'(m)}})\cong
 \dim(V_{\lambda^{'(m)}})\ V_{\lambda^{(m)}}\ ;\]
 we used (\ref{teby1right}) for the first isomorphism and (\ref{teby1right})-(\ref{teby1left}) for the last one. 

\vskip .3cm
{\bf (b)} Let $\lambda^{(m)}$ be an $m$-partition of size $n$ and $\rho\colon H(m,1,n)\to\text{End}(\mathfrak{V})$ a representation of $H(m,1,n)$.
One can construct a representation of $H(m,1,n)$ on the space $U_{\lambda^{(m)}}\otimes \mathfrak{V}$ by moving the  elements of $H(m,1,n)$ through
the basis elements  $\mathcal{X}_{\lambda^{(m)}}$ of $U_{\lambda^{(m)}}$ with the help of instructions from the Proposition \ref{prop-rel} and then 
applying the representation $\rho$. Denote this representation by $V_{\lambda^{(m)}}\boxtimes \mathfrak{V}$. If $\rho$ is a representation from the class $S$, this construction is equivalent to the tensor product $\hat{\otimes}$, 
\begin{equation}\label{boxhatt}V_{\lambda^{(m)}} \boxtimes\mathfrak{V}\cong
V_{\lambda^{(m)}} \hat{\otimes}\mathfrak{V}\ .\end{equation}

Similarly to the construction of a $H(m,1,n)$-module structure on the space $U_{\lambda^{(m)}}\otimes U_{\lambda^{'(m)}}$, one can construct a $H(m,1,n)$-module structure on the tensor product of spaces corresponding to $l$ arbitrary $m$-partitions, with $l\in {\mathbb Z}_{\geq 0}$. One has to replace in the construction the quadratic combinations $\mathcal{X}_{\lambda^{(m)}}\mathcal{X}_{\lambda'^{(m)}}$ by combinations 
$\mathcal{X}_{\lambda_l^{(m)}} \mathcal{X}_{\lambda_{l-1}^{(m)}}\dots\mathcal{X}_{\lambda_1^{(m)}}$
of degree $l$, move the elements of $H(m,1,n)$ through these combinations using the homogeneous (in $\mathcal{X}$) relations (\ref{rel-a1})--(\ref{rel-a2}) and then evaluating them on the vacuum $|\rangle$. Denote this representation by $\mathfrak{W}_{\lambda_l^{(m)},\lambda_{l-1}^{(m)},\dots ,\lambda_{1}^{(m)}}$.
We claim that $\mathfrak{W}_{\lambda_l^{(m)},\lambda_{l-1}^{(m)},\dots ,\lambda_{1}^{(m)}}$ is isomorphic to the direct sum of $\dim(V_{\lambda_1^{(m)}})\times\dots\times\dim(V_{\lambda_{l-1}^{(m)}})$ copies of $V_{\lambda_l^{(m)}}$. Indeed, 
by (\ref{boxhatt}), $\mathfrak{W}_{\lambda_l^{(m)},\lambda_{l-1}^{(m)},\dots ,\lambda_{1}^{(m)}}$ is equivalent to the representation $V_{\lambda_l^{(m)}}\boxtimes \mathfrak{W}_{\lambda_{l-1}^{(m)},\dots ,\lambda_{1}^{(m)}}$. By induction (the induction base is the formula (\ref{tens-prod-form})) the representation $\mathfrak{W}_{\lambda_{l-1}^{(m)},\dots ,\lambda_{1}^{(m)}}$
is isomorphic to  the direct sum of $\dim(V_{\lambda_1^{(m)}})\times\dots\times\dim(V_{\lambda_{l-2}^{(m)}})$ copies of 
$V_{\lambda_{l-1}^{(m)}}$. By (\ref{tens-prod-form}) and (\ref{boxhatt}), the claim follows.

\vskip .3cm
{\bf (c)} The partition $(2,1,1)$ appears in the list from the part $(i)$ of the Lemma \ref{lemm-tens-prod} because 
$${\mathrm{Res}}^4_3(V_{(2,1,1)})\cong {\mathrm{Res}}^4_3(V_{(2,2)})\oplus {\mathrm{Res}}^4_3(V_{(1,1,1,1)})\ .$$ 
For the representation $V_{(2,1,1)}$ the matrix of the operator $\sigma_1\sigma_3$ is
$$\frac{1}{3_q}\left(\begin{array}{ccc}-3_q&0&0\\0&q^{-4}&-q^{-1}4_q\\0&-q^{-1}2_q&-q^2\end{array}\right)\ .$$
Thus, $\text{tr}_{V_{(2,1,1)}}(\sigma_1\sigma_3)=-2+q^{-2}$. In the representation $V_{(1,1,1,1)}$ we have $\sigma_1,\sigma_3\mapsto (-q^{-1})$ while 
in the representation $V_{(2,2)}$ we have $\sigma_1,\sigma_3\mapsto \text{diag}(q,-q^{-1})$, so 
$\text{tr}_{V_{(1,1,1,1)}}(\sigma_1\sigma_3)+\text{tr}_{V_{(2,2)}}(\sigma_1\sigma_3)=q^2+2q^{-2}$.
This differs from $\text{tr}_{V_{(2,1,1)}}(\sigma_1\sigma_3)=-2+q^{-2}$ if and only if $(q+q^{-1})^2\neq 0$ that is, $q+q^{-1}\neq 0$. Therefore to 
establish the formula 
(\ref{tens-prod-form}) for  the representation $V_{(2,1,1)}$ it is enough to calculate the trace of $\sigma_1\sigma_3$ in the representations
$V_{(2,1,1)}\hat{\otimes}V_{\lambda}$ with $\lambda=(2,2), (2,1,1)$ and $(3,1)$. Note that this argument works, in particular, in the classical limit 
$q\rightarrow 1$.

\renewcommand{\theequation}{{B}.{\arabic{equation}}}
\setcounter{equation}{0}
\section*{Appendix B.  \hspace{0.2cm} Bratteli diagrams and their products\vspace{.25cm}}\addcontentsline{toc}{section}{Appendix B.$\ $ Bratteli diagrams 
and their products\vspace{.25cm}}

Here we recall several facts about Bratteli diagrams (see, e.g., \cite{GHJ}) and their graded products. Then we recall the information, needed
in the main body of the text, about the dimensions of the vertices in the powers of the Young graph.

\subsection*{B.1 \hspace{.2cm}  Bratteli diagrams}\addcontentsline{toc}{subsection}{B.1 $\ $ Bratteli diagrams}

A Bratteli diagram is a graded graph; this means that there is a function, called {\it degree}, from the set of vertices to the set of non-negative 
integers. There is only one  degree $0$ vertex which is denoted by $\varnothing$. The edges of the graph can only connect two vertices with neighboring 
degrees (``neighboring" means that the absolute value of the difference of the degrees is $1$). When one draws a Bratteli diagram it is convenient to place on 
the level $a$ all the vertices of degree $a$ (the level $a$ stands for the value of the ordinate $(-a)$). Then a vertex on the level $a$ has ``incoming" edges from 
the level $a-1$ and ``outgoing" edges to the level $a+1$. The dimension of a vertex $x$ is the number of paths which go down from the vertex 
$\varnothing$ to $x$.

\vskip .2cm
In representation theory, one associates a Bratteli diagram to an ascending chain 
\begin{equation}\label{asccoa}{\cal{A}}_0={\mathbb{C}}\subset {\cal{A}}_1\subset\dots\subset {\cal{A}}_n\subset\dots\end{equation}
of associative algebras: vertices of degree $k$ correspond to representations (depending on circumstances, indecomposable, irreducible \emph{etc.}) 
of the algebra ${\cal{A}}_k$ and the Bratteli diagram visualizes the branching rules for the pairs $({\cal{A}}_{k+1},{\cal{A}}_k)$, $k\in {\mathbb{Z}}$.
In this situation, the dimension of a vertex is simply equal to the dimension of the corresponding representation. 

\vskip .2cm
Let $\mathfrak{G}^{(1)}$ and $\mathfrak{G}^{(2)}$  be two Bratteli diagrams. The vertices of the product $\mathfrak{G}^{(1)}\times\mathfrak{G}^{(2)}$ are by definition couples $(x,y)$ where $x$ is a vertex of $\mathfrak{G}^{(1)}$ and $y$ is a vertex of $\mathfrak{G}^{(2)}$. The degree of $(x,y)$ is the sum of the degree of $x$ and the degree of $y$. The top vertex $\mathfrak{G}^{(1)}\times\mathfrak{G}^{(2)}$ which is denoted again by $\varnothing$ is therefore $(\varnothing,\varnothing)$.
If there is an edge between $x$ and $x'$ in $\mathfrak{G}^{(1)}$ one draws an edge between $(x,y)$ and $(x',y)$ for all $y$; we say that these edges are of type $1$. If there is an edge between $y$ and $y'$ in $\mathfrak{G}^{(2)}$ one draws an edge between $(x,y)$ and $(x,y')$ for all $x$; we say that these edges are of type $2$. By definition these are all edges: each edge is either of type $1$ or of type $2$. 

\vskip .2cm
Iterating, we define the product of an arbitrary number $m$ of Bratteli diagrams.

\subsection*{B.2 \hspace{.2cm} Dimensions of vertices of the product}\addcontentsline{toc}{subsection}{B.2 $\ $ Dimensions of vertices of the product}

Let $\mathfrak{G}$ be a Bratteli diagram. Denote by $\mathfrak{P}(\mathfrak{G})$ the set of paths which begin at the top vertex
of $\mathfrak{G}$ and go down. For $p\in \mathfrak{P}(\mathfrak{G})$ denote by ${\mathcal{E}}(p)$ the collection of edges of $p$ and by $\text{end} (p)$
the end point of $p$; if $z=\text{end} (p)$ then $\deg (z)$ equals the length of $p$, the cardinality of ${\mathcal{E}}(p)$. The set  ${\mathcal{E}}(p)$ is 
naturally ordered:  the steps in the path follow one after another. 

\vskip .2cm
Let $(x,y)$ be a vertex of $\mathfrak{G}^{(1)}\times\mathfrak{G}^{(2)}$. Let $p$ be a path from $\mathfrak{P}(\mathfrak{G}^{(1)}\times\mathfrak{G}^{(2)})$ 
with $\text{end} (p)=(x,y)$. The set  ${\mathcal{E}}(p)$ is the disjoint union of two subsets,  ${\mathcal{E}}_1(p)$ and  
${\mathcal{E}}_2(p)$, where  ${\mathcal{E}}_1(p)$ (respectively,  ${\mathcal{E}}_2(p)$) is the subset of  ${\mathcal{E}}(p)$ consisting of edges of type 1 
(respectively, type 2). Each edge from ${\mathcal{E}}_1(p)$ naturally defines an edge in the graph $\mathfrak{G}^{(1)}$ and the set of edges thus defined
form a path $p_1$ in the graph $\mathfrak{G}^{(1)}$ going down from $\varnothing$ of  $\mathfrak{G}^{(1)}$ to $x$, $p_1\in 
\mathfrak{P}(\mathfrak{G}^{(1)})$; similarly, each edge 
 from ${\mathcal{E}}_2(p)$ naturally defines an edge in the graph $\mathfrak{G}^{(2)}$ and the set of edges thus defined form a path $p_2$ in the graph 
$\mathfrak{G}^{(2)}$ going down from $\varnothing$ of $\mathfrak{G}^{(2)}$ to $y$, $p_2\in \mathfrak{P}(\mathfrak{G}^{(2)})$. We have therefore a map from 
$\mathfrak{P}(\mathfrak{G}^{(1)}\times\mathfrak{G}^{(2)})$ to the product $\mathfrak{P}(\mathfrak{G}^{(1)})\times\mathfrak{P}(\mathfrak{G}^{(2)})$,
defined by 
\begin{equation}\label{prftsofp} \pi :\, p\mapsto (p_1,p_2)\ .\end{equation} 
One cannot reconstruct uniquely the path $p$ knowing the paths $p_1$ and $p_2$. Let $a$ be the degree of $x$ and $b$ the degree of $y$. 
Any order $\succ$ on the union ${\mathcal{E}}(p_1)\cup{\mathcal{E}}(p_2)$ which is compatible with the natural orders on 
${\mathcal{E}}_1(p)$ and ${\mathcal{E}}_2(p)$ (in the sense that if a step $\gamma$ is after a step $\gamma'$ in ${\mathcal{E}}_1(p)$ or ${\mathcal{E}}_2(p)$
then  $\gamma$ is after $\gamma'$ with respect to the order $\succ$) determines a path from $\mathfrak{P}(\mathfrak{G}^{(1)}\times\mathfrak{G}^{(2)})$
In other words, in the sequence of $a+b$ edges of a path of length $a+b$ from $\mathfrak{P}(\mathfrak{G}^{(1)}\times\mathfrak{G}^{(2)})$ one can assign 
the type $1$ to an arbitrarily chosen subset of $a$ edges so the cardinality of the preimage of the element $ (p_1,p_2)$ with respect to the map $\pi$ given 
by (\ref{prftsofp}) is $\left(\begin{array}{c}a+b\\b\end{array}\right)$; this cardinality depends only on the end points $x$ and $y$ of the paths $p_1$ and $p_2$
so we have
\begin{equation}\label{dim-xy}
\dim \bigl((x,y)\bigr)=\left(\begin{array}{c}a+b\\b\end{array}\right)\dim(x)\dim(y)\ .
\end{equation}
 
{}For a Bratteli diagram $\mathfrak{G}$ define $D(\mathfrak{G})_a$ by
\begin{equation}\label{dim-Ga}
D(\mathfrak{G})_a:=\sum\limits_{x:\deg(x)=a}\bigl(\dim(x)\bigr)^2\ .
\end{equation}

When the Bratteli diagram is associated to an ascending chain of finite-dimensional semi-simple associative algebras, like (\ref{asccoa}), 
and the vertices correspond to irreducible representations, the number $D(\mathfrak{G})_a$ is the dimension of the algebra ${\cal{A}}_a$. 

\vskip .2cm
By (\ref{dim-xy}), we have for the product
\begin{equation}\label{dim-prod}\begin{array}{lcl}
D(\mathfrak{G}^{(1)}\times\mathfrak{G}^{(2)})_c&=&\!\!\!\!\!\!\!
{\displaystyle\sum_{\textrm{\footnotesize
$\begin{array}{c}a,b:a+b=c\\ x:deg(x)=a\\ y:deg(y)=b\end{array}$}} } {\textstyle (\dim\bigl((x,y)\bigr)^2}
=\!\!\!\!\!\!\!{\displaystyle\sum_{\textrm{\footnotesize$\begin{array}{c}a,b:a+b=c\\ x:deg(x)=a\\ y:deg(y)=b\end{array}$}}}\!\!\!
{\textstyle \left(\begin{array}{c}a+b\\b\end{array}\right)^2\bigl(\dim(x)\bigr)^2\bigl(\dim(y)\bigr)^2}\\[5em]
&=&{\displaystyle\sum_{a=0}^c}\left(\begin{array}{c}c\\a\end{array}\right)^2
D(\mathfrak{G}^{(1)})_a\ \  D(\mathfrak{G}^{(2)})_{c-a}\ .\end{array}
\end{equation}

\subsection*{B.3 \hspace{.2cm} Powers of the Young graph}\addcontentsline{toc}{subsection}{B.3 $\ $ Powers of the Young graph}

As we have seen in Section \ref{ygcyclo}, the vertices of the Bratteli diagram for the chain (with respect to $n$) of the 
algebras $H(m,1,n)$
naturally correspond to $m$-partitions, the level $a$ consists of all $m$-partitions of $a$;
the edges outgoing from the level $a$ correspond to inclusions of $m$-partitions of $a$ into $m$-partitions of $a+1$. Thus the Bratteli diagram for the
chain $H(m,1,n)$ is the $m^{\textrm{th}}$ power of the Young graph. 

\paragraph{1. Dimensions.}
We need to determine the dimensions of the vertices of the powers of the Young graph. We recall the definition of the hook length and the formula for 
the dimensions of the vertices of the Young graph. 
For a node $\alpha$ of a Young diagram the hook of $\alpha$ is the set of nodes containing $\alpha$ and the nodes 
which lie either under $\alpha$ in the same column or to the right of $\alpha$ in the same row. The hook length $h_{\alpha}$ of $\alpha$ is the number of 
nodes in the hook of $\alpha$. The dimension of a representation (of a symmetric group) corresponding to a partition $\lambda$ of $n$ is given by the 
classical hook formula,
\begin{equation}\dim(V_{\lambda})=\frac{n!}{\prod\limits_{\alpha\in\lambda}h_{\alpha}}\ ,\end{equation}
where the product $\prod\limits_{\alpha\in \lambda}h_{\alpha}$ means the product of the hook lengths of all nodes $\alpha$ of the Young diagram of 
shape $\lambda$.

\vskip .2cm
Consider an $m$-partition $\lambda^{(m)}:=(\lambda_1,\dots,\lambda_m)$ such that $|\lambda^{(m)}|=n$ (we remind that 
$|\lambda^{(m)}|=|\lambda_1|+\dots+|\lambda_m|$). We denote by $V_{\lambda^{(m)}}$ the irreducible representation of $H(m,1,n)$ associated with 
$\lambda^{(m)}$. By the generalization of (\ref{dim-xy}) to the product of $m$
graded graphs, the dimension of $V_{\lambda^{(m)}}$ is 
\begin{equation}\label{dimV}
\dim(V_{\lambda^{(m)}})=\frac{n!}{\vert\lambda_1\vert !\dots \vert\lambda_m\vert !}\
\frac{\vert\lambda_1\vert !}{\prod\limits_{\alpha\in\lambda_1}h_{\alpha}}\dots \frac{\vert\lambda_m\vert !}{\prod\limits_{\alpha\in\lambda_m}h_{\alpha}}
=\frac{n!}{\prod\limits_{i=1}^m\prod\limits_{\alpha\in\lambda_i}h_{\alpha}}\ ,
\end{equation}

\bigskip
\begin{lemm}
{\hspace{-.2cm}.\hspace{.2cm}}
 \label{sum-dim}
We have
\begin{equation}\label{sumdim}
\sum\limits_{\lambda^{(m)}}(\dim(V_{\lambda^{(m)}}))^2=n!m^n\ ,
\end{equation}
where the sum is over all $m$-partitions $\lambda^{(m)}=(\lambda_1,\dots,\lambda_m)$ such that $|\lambda^{(m)}|=n$.
\end{lemm}

\emph{Proof.}
For $m=1$, we know that the representations $V_{\lambda}$ where $\lambda$ is a partition of $n$ are all the irreducible representations of 
the symmetric group 
$S_n$ and so:
\begin{equation}\label{sucodis}
\sum\limits_{\lambda}(\dim(V_{\lambda}))^2=\sum\limits_{\lambda}\left(\frac{n!}{\prod\limits_{\alpha\in\lambda}h_{\alpha}}\right)^2=n!\ .\end{equation}
The proof of (\ref{sumdim}) is by induction on $m$. We have:
\[\begin{array}{ll}\sum\limits_{\lambda^{(m)}:\vert\lambda^{(m)}\vert=n}(\dim(V_{\lambda^{(m)}}))^2 &=
{\displaystyle{\sum_{j=0}^n}}
\left(\frac{n!}{(n-j)!j!}\right)^2\hspace{-.1cm} 
\sum\limits_{\lambda^{(1)}:\vert\lambda^{(1)}\vert=j}
\bigl(\dim(V_{\lambda^{(1)}})\bigr)^2\hspace{-.1cm}
\sum\limits_{\lambda^{(m-1)}:|\lambda^{(m-1)}|=n-j}\hspace{-.5cm}
\bigl(\dim(V_{\lambda^{(m-1)}})\bigr)^2 
\\[2em]
& = \sum\limits_{j=0}^n\Biggl(\left(\frac{n!}{(n-j)!j!}\right)^2\cdot j!\cdot (n-j)!\cdot (m-1)^{n-j}\Biggr)\\[2.5em]
& =n!\cdot\sum\limits_{j=0}^n \frac{n!}{(n-j)!j!}(m-1)^{n-j}=n!\, m^n\ ;
\end{array}\]
here $\lambda^{(1)}$ is a usual partition and $\lambda^{(m-1)}$ is an $(m-1)$-partition.
In the first equality we used (\ref{dim-prod}); in the second equality we used (\ref{sucodis}) and the induction hypothesis; we simplified the result in the 
third equality and used the binomial theorem in the fourth equality.\hfill$\square$

\paragraph{Remark.}  The $m^{\textrm{th}}$ power of the Young graph is an $m$-differential poset; the formula (\ref{sumdim}) holds for arbitrary $m$-differential posets 
(see \cite{St} for definitions and details). 

\paragraph{2. Standard $m$-tableaux and dimensions.} It is well known that the dimension of a representation of a symmetric group corresponding to
some Young diagram $\lambda$ equals the dimension of the corresponding vertex in the Young graph and equals the number of standard tableaux of  
the shape $\lambda$. It is straightforward to generalize these equalities to the cyclotomic case: the dimension of a representation of the algebra $H(m,1,n)$
corresponding to some Young $m$-diagram $\lambda^{(m)}$ equals the dimension of the corresponding vertex in the $m^{\textrm{th}}$ power of the
Young graph and equals the number of standard $m$-tableaux of the shape $\lambda^{(m)}$.

\vskip .2cm
With the help of the Lemma \ref{sum-dim} we check that the sum of the squares of the dimensions of the representations constructed in Subsection \ref{suserep} 
is not less than
the dimension of the algebra $H(m,1,n)$
(cf. (\ref{dim-ine})):
\begin{equation}\label{sumdim2}
\sum\limits_{\lambda^{(m)}}(\dim(V_{\lambda^{(m)}}))^2\geq \dim(H(m,1,n))\ .
\end{equation}

\paragraph{3. Example: square of the Young graph.} Below the beginning of the Bratteli diagram for the chain of algebras $H(2,1,n)$, the square of the 
Young graph, is drawn. The labels on the edges correspond to the eigenvalues of the Jucys--Murphy elements of  $H(2,1,n)$ (the edges going down 
from the level $i$ to the level $i+1$ are labeled by the eigenvalues of the element $J_{i+1}$; the top vertex is situated at level 0).

\vspace{-1cm}

\setlength{\unitlength}{2600sp}
\begingroup\makeatletter\ifx\SetFigFontNFSS\undefined\gdef\SetFigFontNFSS#1#2#3#4#5{
  \reset@font\fontsize{#1}{#2pt}  \fontfamily{#3}\fontseries{#4}\fontshape{#5}  \selectfont}
\fi\endgroup
\begin{picture}(506,5718)(576,-5800)
\thinlines{\put(5471,-1591){\line(-1,-1){810}}}{\put(5561,-1591){\line( 1,-1){855}}}{\put(4666,-2716){\line(-5,-2){2025}}}{\put(4701,-2761){\line(-1,-1){765}}}{\put(4751,-2761){\line( 4,-3){1080}}}{\put(6353,-2662){\line(-1,-2){450}}}{\put(6439,-2700){\line( 4,-3){1080}}}{\put(6555,-2657){\line( 3,-1){2430}}}{\put(2511,-4041){\line(-1,-1){1665}}}{\put(2601,-4041){\line(-1,-2){810}}}{\put(3941,-4041){\line(-4,-3){2160}}}{\put(4036,-4086){\line(-1,-2){810}}}{\put(6081,-4086){\line( 3,-5){945}}}{\put(6086,-4041){\line( 4,-3){2160}}}{\put(9026,-3996){\line(-1,-2){810}}}{\put(2602,-4073){\line( 6,-5){1890}}}{\put(5918,-4059){\line(-5,-6){1350}}}{\put(5956,-4048){\line(-1,-5){315}}}{\put(4058,-4063){\line( 1,-1){1620}}}{\put(7699,-4044){\line(-2,-5){630}}}{\put(7771,-4041){\line( 1,-1){1575}}}{\put(9120,-3989){\line( 5,-3){2700}}}{\put(9093,-4063){\line( 1,-1){1530}}}{\put(7836,-4029){\line( 5,-3){2700}}
}\put(5446,-3946){\makebox(0,0)[lb]{\smash{{\SetFigFontNFSS{12}{14.4}{\rmdefault}{\mddefault}{\updefault}{$\scriptstyle{\left(\Box\,,\,\Box\right)}$}}}}}\put(3376,-3901){\makebox(0,0)[lb]{\smash{{\SetFigFontNFSS{12}{14.4}{\rmdefault}{\mddefault}{\updefault}{$\scriptstyle{\left(\Box\!\Box\,,\,\varnothing\right)}$}}}}}\put(1981,-3901){\makebox(0,0)[lb]{\smash{{\SetFigFontNFSS{12}{14.4}{\rmdefault}{\mddefault}{\updefault}{$\scriptstyle{\left(\!\!\!\begin{array}{l}\scriptstyle{\Box}\\[-0.75em]\scriptstyle{\Box}\end{array}\!\!\!,\,
\varnothing\right)}$}}}}}\put(4096,-2626){\makebox(0,0)[lb]{\smash{{\SetFigFontNFSS{12}{14.4}{\rmdefault}{\mddefault}{\updefault}{$\scriptstyle{\left(\Box\,,\,\varnothing\right)}$}}}}}\put(5896,-2626){\makebox(0,0)[lb]{\smash{{\SetFigFontNFSS{12}{14.4}{\rmdefault}{\mddefault}{\updefault}{$\scriptstyle{\left(\varnothing\,,\,\Box\right)}$}}}}}\put(4800,-1450){\makebox(0,0)[lb]{\smash{{\SetFigFontNFSS{12}{14.4}{\rmdefault}{\mddefault}{\updefault}{$\scriptstyle{
\varnothing =\left(\varnothing\,,\,\varnothing\right)}$}}}}}\put(8506,-3856){\makebox(0,0)[lb]{\smash{{\SetFigFontNFSS{12}{14.4}{\rmdefault}{\mddefault}{\updefault}{$\scriptstyle{\left(\varnothing\,,\!\!\!\begin{array}{l}\scriptstyle{\Box}\\[-0.75em]\scriptstyle{\Box}\end{array}\!\!\!\right)}$}}}}}\put(6976,-3906){\makebox(0,0)[lb]{\smash{{\SetFigFontNFSS{12}{14.4}{\rmdefault}{\mddefault}{\updefault}{$\scriptstyle{\left(\varnothing\,,\,\Box\!\Box\right)}$}}}}}\put(1081,-5981){\makebox(0,0)[lb]{\smash{{\SetFigFontNFSS{12}{14.4}{\rmdefault}{\mddefault}{\updefault}{$\scriptstyle{\left(\!\!\!\begin{array}{l}\scriptstyle{\Box}\!\scriptstyle{\Box}\\[-0.75em]\scriptstyle{\Box}\end{array}\!\!\!,\,\varnothing\right)}$}}}}}\put(91,-5981){\makebox(0,0)[lb]{\smash{{\SetFigFontNFSS{12}{14.4}{\rmdefault}{\mddefault}{\updefault}{$\scriptstyle{\left(\!\!\!\begin{array}{l}\scriptstyle{\Box}\\[-0.75em]\scriptstyle{\Box}\\[-0.75em]\scriptstyle{\Box}\end{array}\!\!\!,\,\varnothing\right)}$}}}}}\put(2476,-5981){\makebox(0,0)[lb]{\smash{{\SetFigFontNFSS{12}{14.4}{\rmdefault}{\mddefault}{\updefault}{$\scriptstyle{\left(\Box\!\Box\!\Box\,,\,\varnothing\right)}$}}}}}\put(3961,-5981){\makebox(0,0)[lb]{\smash{{\SetFigFontNFSS{12}{14.4}{\rmdefault}{\mddefault}{\updefault}{$\scriptstyle{\left(\!\!\!\begin{array}{l}\scriptstyle{\Box}\\[-0.75em]\scriptstyle{\Box}\end{array}\!\!\!,\,\Box\right)}$}}}}}\put(4996,-5981){\makebox(0,0)[lb]{\smash{{\SetFigFontNFSS{12}{14.4}{\rmdefault}{\mddefault}{\updefault}{$\scriptstyle{\left(\Box\!\Box\,,\,\Box\right)}$}}}}}\put(6346,-5936){\makebox(0,0)[lb]{\smash{{\SetFigFontNFSS{12}{14.4}{\rmdefault}{\mddefault}{\updefault}{$\scriptstyle{\left(\Box\,,\,\Box\!\Box\right)}$}}}}}\put(8641,-5891){\makebox(0,0)[lb]{\smash{{\SetFigFontNFSS{12}{14.4}{\rmdefault}{\mddefault}{\updefault}{$\scriptstyle{\left(\varnothing\,,\,\Box\!\Box\!\Box\right)}$}}}}}\put(7651,-5936){\makebox(0,0)[lb]{\smash{{\SetFigFontNFSS{12}{14.4}{\rmdefault}{\mddefault}{\updefault}{$\scriptstyle{\left(\Box\,,\!\!\!\begin{array}{l}\scriptstyle{\Box}\\[-0.75em]\scriptstyle{\Box}\end{array}\!\!\!\right)}$}}}}}\put(9991,-5891){\makebox(0,0)[lb]{\smash{{\SetFigFontNFSS{12}{14.4}{\rmdefault}{\mddefault}{\updefault}{$\scriptstyle{\left(\varnothing\,,\!\!\!\begin{array}{l}\scriptstyle{\Box}\!\scriptstyle{\Box}\\[-0.75em] \scriptstyle{\Box}\end{array}\!\!\!\right)}$}}}}}\put(11341,-5891){\makebox(0,0)[lb]{\smash{{\SetFigFontNFSS{12}{14.4}{\rmdefault}{\mddefault}{\updefault}{$\scriptstyle{\left(\varnothing\,,\!\!\!\begin{array}{l}\scriptstyle{\Box}\\[-0.75em]\scriptstyle{\Box}\\[-0.75em]\scriptstyle{\Box}\end{array}\!\!\!\right)}$}}}}}
\put(4550,-1951){\makebox(0,0)[lb]{\smash{{\SetFigFontNFSS{12}{14.4}{\rmdefault}{\mddefault}{\updefault}{$\scriptstyle{v_1}$}}}}}\put(5700,-1951){\makebox(0,0)[lb]{\smash{{\SetFigFontNFSS{12}{14.4}{\rmdefault}{\mddefault}{\updefault}{$\scriptstyle{v_2}$}}}}}\put(3016,-3076){\makebox(0,0)[lb]{\smash{{\SetFigFontNFSS{12}{14.4}{\rmdefault}{\mddefault}{\updefault}{$\scriptstyle{v_1q^{-2}}$}}}}}\put(4006,-3301){\makebox(0,0)[lb]{\smash{{\SetFigFontNFSS{12}{14.4}{\rmdefault}{\mddefault}{\updefault}{$\scriptstyle{v_1q^2}$}}}}}\put(4951,-3076){\makebox(0,0)[lb]{\smash{{\SetFigFontNFSS{12}{14.4}{\rmdefault}{\mddefault}{\updefault}{$\scriptstyle{v_2}$}}}}}\put(5941,-3211){\makebox(0,0)[lb]{\smash{{\SetFigFontNFSS{12}{14.4}{\rmdefault}{\mddefault}{\updefault}{$\scriptstyle{v_1}$}}}}}\put(6796,-3166){\makebox(0,0)[lb]{\smash{{\SetFigFontNFSS{12}{14.4}{\rmdefault}{\mddefault}{\updefault}{$\scriptstyle{v_2q^2}$}}}}}\put(7606,-3076){\makebox(0,0)[lb]{\smash{{\SetFigFontNFSS{12}{14.4}{\rmdefault}{\mddefault}{\updefault}{$\scriptstyle{v_2q^{-2}}$}}}}}\put(651,-5070){\makebox(0,0)[lb]{\smash{{\SetFigFontNFSS{12}{14.4}{\rmdefault}{\mddefault}{\updefault}{$\scriptstyle{v_1q^{-4}}$}}}}}\put(1446,-5000){\makebox(0,0)[lb]{\smash{{\SetFigFontNFSS{12}{14.4}{\rmdefault}{\mddefault}{\updefault}{$\scriptstyle{v_1q^2}$}}}}}\put(4700,-4600){\makebox(0,0)[lb]{\smash{{\SetFigFontNFSS{12}{14.4}{\rmdefault}{\mddefault}{\updefault}{$\scriptstyle{v_1q^{-2}}$}}}}}\put(7540,-4671){\makebox(0,0)[lb]{\smash{{\SetFigFontNFSS{12}{14.4}{\rmdefault}{\mddefault}{\updefault}{$\scriptstyle{v_2q^4}$}}}}}\put(2106,-5351){\makebox(0,0)[lb]{\smash{{\SetFigFontNFSS{12}{14.4}{\rmdefault}{\mddefault}{\updefault}{$\scriptstyle{v_1q^{-2}}$}}}}}\put(2800,-5256){\makebox(0,0)[lb]{\smash{{\SetFigFontNFSS{12}{14.4}{\rmdefault}{\mddefault}{\updefault}{$\scriptstyle{v_1q^4}$}}}}}\put(3600,-5100){\makebox(0,0)[lb]{\smash{{\SetFigFontNFSS{12}{14.4}{\rmdefault}{\mddefault}{\updefault}{$\scriptstyle{v_2}$}}}}}\put(4250,-4526){\makebox(0,0)[lb]{\smash{{\SetFigFontNFSS{12}{14.4}{\rmdefault}{\mddefault}{\updefault}{$\scriptstyle{v_2}$}}}}}\put(5550,-4941){\makebox(0,0)[lb]{\smash{{\SetFigFontNFSS{12}{14.4}{\rmdefault}{\mddefault}{\updefault}{$\scriptstyle{v_1q^2}$}}}}}\put(6010,-5300){\makebox(0,0)[lb]{\smash{{\SetFigFontNFSS{12}{14.4}{\rmdefault}{\mddefault}{\updefault}{$\scriptstyle{v_2q^2}$}}}}}\put(7070,-4350){\makebox(0,0)[lb]{\smash{{\SetFigFontNFSS{12}{14.4}{\rmdefault}{\mddefault}{\updefault}{$\scriptstyle{v_1}$}}}}}\put(6300,-4400){\makebox(0,0)[lb]{\smash{{\SetFigFontNFSS{12}{14.4}{\rmdefault}{\mddefault}{\updefault}{$\scriptstyle{v_2q^{-2}}$}}}}}\put(7900,-5301){\makebox(0,0)[lb]{\smash{{\SetFigFontNFSS{12}{14.4}{\rmdefault}{\mddefault}{\updefault}{$\scriptstyle{v_1}$}}}}}\put(8700,-4700){\makebox(0,0)[lb]{\smash{{\SetFigFontNFSS{12}{14.4}{\rmdefault}{\mddefault}{\updefault}{$\scriptstyle{v_2q^{-2}}$}}}}}\put(9900,-5100){\makebox(0,0)[lb]{\smash{{\SetFigFontNFSS{12}{14.4}{\rmdefault}{\mddefault}{\updefault}{$\scriptstyle{v_2q^2}$}}}}}\put(10600,-5000){\makebox(0,0)[lb]{\smash{{\SetFigFontNFSS{12}{14.4}{\rmdefault}{\mddefault}{\updefault}{$\scriptstyle{v_2q^{-4}}$}}}}}
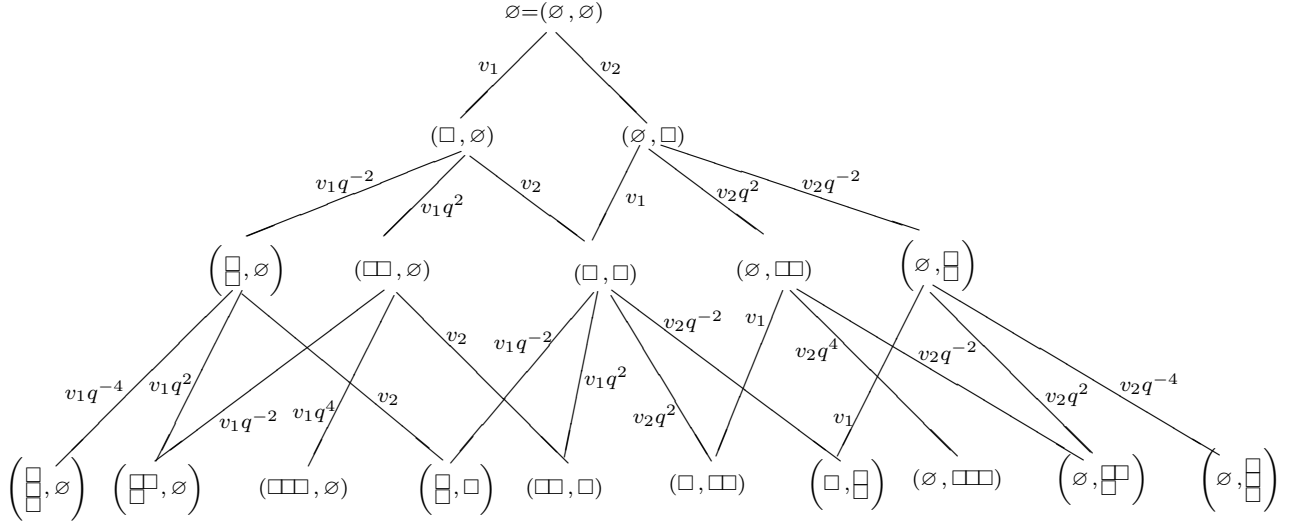
\put(2900,-6800){\begin{pic}\label{Brat-cyclo}Bratteli diagram (four first levels) for $H(m,1,n)$ with $m=2$.
\end{pic}}
\end{picture}

\vskip 2.5cm

\renewcommand{\theequation}{{C}.{\arabic{equation}}}
\setcounter{equation}{0}
\section*{Appendix C.  \hspace{0.2cm} Examples}\addcontentsline{toc}{section}{Appendix C.$\ $ Examples}

Here we illustrate the construction of irreducible representations of the algebras $H(m,1,n)$ on several examples with $m=2$ and small $n$. 
{}For these examples we write down the formulas (\ref{rel-a1})-(\ref{rel-a2}) and (\ref{rep-a1})-(\ref{rep-a2}) from Section \ref{ygcyclo'}. 

\paragraph{1. The representation of $H(2,1,2)$ corresponding to the 2-partition} $\hspace{-.2cm}(\Box\, ,\Box)$.

\vskip .2cm
The dimension of this representation is 2. We choose a basis  
\[{\cal{X}}_1:={\cal{X}}_{\left(\fbox{\scriptsize{$1$}}\,,\,\fbox{\scriptsize{$2$}}\right)}\ ,\quad {\cal{X}}_2:={\cal{X}}_{\left(\fbox{\scriptsize{$2$}}\,,\,\fbox{\scriptsize{$1$}}\right)}\ .\]

The formulas (\ref{rel-a1})-(\ref{rel-a2}) take the form
\[ \textstyle
\left(\sigma_1+\frac{(q-q^{-1})v_2}{v_1-v_2}\right) {\cal{X}}_1={\cal{X}}_2\left(\sigma_1+\frac{(q-q^{-1})v_1}{v_2-v_1}\right)\ ,
\quad \left(\sigma_1+\frac{(q-q^{-1})v_1}{v_2-v_1}\right) {\cal{X}}_2={\cal{X}}_1\left(\sigma_1+\frac{(q-q^{-1})v_2}{v_1-v_2}\right) \]
and  
\[ (\tau-v_1){\cal{X}}_1=0\quad ,
\qquad (\tau-v_2){\cal{X}}_2=0\ .\]

\vskip .2cm 
The matrices corresponding to the action (\ref{rep-a1})-(\ref{rep-a2}) of the generators of $H(2,1,2)$ in the basis above are given by:
\begin{equation}\label{rep1virgule1}\sigma_1\mapsto \left(\begin{array}{cc}-(q-q^{-1})\frac{v_2}{v_1-v_2}&
\frac{q v_1-q^{-1}v_2}{v_1-v_2}\\[1em]
\frac{q v_2-q^{-1}v_1}{v_2-v_1}&-(q-q^{-1})\frac{v_1}{v_2-v_1}\end{array}\right)\ ,\quad
\tau\mapsto \text{diag}(v_1,v_2)
\ .\end{equation}

The Gram matrix $\langle{\cal{X}}_i,{\cal{X}}_j\rangle_{i,j=1,2}\ $, given by (\ref{H-scal-prod1})--(\ref{H-scal-prod2}), reads
$$\text{diag}\left(\frac{q^{-1}v_1-qv_2}{v_1-v_2},\frac{qv_1-q^{-1}v_2}{v_1-v_2}\right)\ .$$

\paragraph{2. The representation of $H(2,1,3)$ corresponding to the 2-partition}
$\hspace{-.2cm}\left(\begin{array}{c}\Box\\[-0.6em] \Box\end{array},\Box\right)$.  

\vskip .2cm
The representation has dimension 3 and we choose a basis 
\[{\cal{X}}_1:={\cal{X}}_{\textrm{\tiny{$\left(\!\!\begin{array}{c}\fbox{\scriptsize{$1$}}\\[-0.10em] 
\fbox{\scriptsize{$2$}}\end{array}\!,\,\fbox{\scriptsize{$3$}}\right)$}}}\ ,\quad 
{\cal{X}}_2:={\cal{X}}_{\textrm{\tiny{$\left(\!\!\begin{array}{c}\fbox{\scriptsize{$1$}}\\[-0.10em] 
\fbox{\scriptsize{$3$}}\end{array}\!,\,\fbox{\scriptsize{$2$}}\right)$}}}\ ,\quad 
{\cal{X}}_3:={\cal{X}}_{\textrm{\tiny{$\left(\!\!\begin{array}{c}\fbox{\scriptsize{$2$}}\\[-0.10em] 
\fbox{\scriptsize{$3$}}\end{array}\!,\,\fbox{\scriptsize{$1$}}\right)$}}}\ .\]

The formulas (\ref{rel-a1})-(\ref{rel-a2}) take the form
\[\begin{array}{ll}
\left(\sigma_1+\frac{(q-q^{-1})v_1q^{-2}}{v_1-v_1q^{-2}}\right){\cal{X}}_1=0\ ,&
\left(\sigma_1+\frac{(q-q^{-1})v_2}{v_1-v_2}\right){\cal{X}}_2={\cal{X}}_3\left(\sigma_1+\frac{(q-q^{-1})v_1}{v_2-v_1}\right)\ ,\\[1em]
\left(\sigma_1+\frac{(q-q^{-1})v_1}{v_2-v_1}\right){\cal{X}}_3={\cal{X}}_2\left(\sigma_1+\frac{(q-q^{-1})v_2}{v_1-v_2}\right)\ ,&
\end{array}\] 

\vskip .1cm
\[\begin{array}{ll}
\left(\sigma_2+\frac{(q-q^{-1})v_2}{v_1q^{-2}-v_2}\right){\cal{X}}_1={\cal{X}}_2\left(\sigma_2+\frac{(q-q^{-1})v_1q^{-2}}{v_2-v_1q^{-2}}\right)\ ,&
\left(\sigma_2+\frac{(q-q^{-1})v_1q^{-2}}{v_2-v_1q^{-2}}\right){\cal{X}}_2={\cal{X}}_1\left(\sigma_2+\frac{(q-q^{-1})v_2}{v_1q^{-2}-v_2}\right)\ ,\\[1em]
\left(\sigma_2+\frac{(q-q^{-1})v_1q^{-2}}{v_1-v_1q^{-2}}\right){\cal{X}}_3=0\ ,&
\end{array}\] 
and
\[(\tau-v_1){\cal{X}}_1=0\ ,
\quad (\tau-v_1){\cal{X}}_2=0\ ,
\quad (\tau-v_2){\cal{X}}_3=0\ .
\]

\vskip .2cm 
The matrices corresponding to the action (\ref{rep-a1})-(\ref{rep-a2}) of the generators of $H(2,1,3)$ in the basis above are given by:
\[\sigma_1\mapsto\!\!\left(\begin{array}{ccc}-q^{-1}&0&0\\[1em]
0&-\frac{(q-q^{-1})v_2}{v_1-v_2}&\frac{qv_1-q^{-1}v_2}{v_1-v_2}\\[1em]
0&\frac{qv_2-q^{-1}v_1}{v_2-v_1}&-\frac{(q-q^{-1})v_1}{v_2-v_1}\end{array}\right)\ ,\quad
\sigma_2\mapsto\!\! \left(\begin{array}{ccc}-\frac{(q-q^{-1})v_2}{v_1q^{-2}-v_2}&\frac{v_1q^{-1}-q^{-1}v_2}{v_1q^{-2}-v_2}&0\\[1em]
\frac{qv_2-v_1q^{-3}}{v_2-v_1q^{-2}}&-\frac{(q-q^{-1})v_1q^{-2}}{v_2-v_1q^{-2}}&0\\[1em]
0&0&-q^{-1}\end{array}\right)\]
and
\[\tau\mapsto \text{diag}(v_1,v_1,v_2)\ .\]

The Gram matrix $\langle{\cal{X}}_i,{\cal{X}}_j\rangle_{i,j=1,2,3}\ $, given by (\ref{H-scal-prod1})--(\ref{H-scal-prod2}), reads
$$\text{diag}\left(\frac{q^{-2}v_1-q^2v_2}{v_1-v_2},1,\frac{qv_1-q^{-1}v_2}{q^{-1} v_1-qv_2}\right)\ .$$

\paragraph{3. The representation of $H(2,1,4)$ labeled by the 2-partition}  
$\hspace{-.2cm}\left(\begin{array}{c}\Box\\[-0.6em] \Box\end{array},\Box\hspace{-0.05cm}\Box\right)$. 

\vskip .2cm
The representation has dimension 6 and we choose a basis 
\[\begin{array}{c}{\cal{X}}_1:={\cal{X}}_{\textrm{\tiny{$\left(\!\!\begin{array}{c}\fbox{\scriptsize{$1$}}\\[-0.10em] \fbox{\scriptsize{$2$}}\end{array}
\!,\,\fbox{\scriptsize{$3$}}\fbox{\scriptsize{$4$}}\right)$}}}\ ,\quad {\cal{X}}_2:={\cal{X}}_{\textrm{\tiny{$\left(\!\!\begin{array}{c}\fbox{\scriptsize{$1$}}\\[-0.10em] 
\fbox{\scriptsize{$3$}}\end{array}\!,\,\fbox{\scriptsize{$2$}}\fbox{\scriptsize{$4$}}\right)$}}}\ ,\quad {\cal{X}}_3:={\cal{X}}_{\textrm{\tiny{$\left(\!\!\begin{array}{c}
\fbox{\scriptsize{$1$}}\\[-0.10em] \fbox{\scriptsize{$4$}}\end{array}\!,\,\fbox{\scriptsize{$2$}}\fbox{\scriptsize{$3$}}\right)$}}}\ ,\quad
{\cal{X}}_4:={\cal{X}}_{\textrm{\tiny{$\left(\!\!\begin{array}{c}\fbox{\scriptsize{$2$}}\\[-0.10em] \fbox{\scriptsize{$3$}}\end{array}\!,\,\fbox{\scriptsize{$1$}}
\fbox{\scriptsize{$4$}}\right)$}}}\ ,\\[3em]
{\cal{X}}_5:={\cal{X}}_{\textrm{\tiny{$\left(\!\!\begin{array}{c}\fbox{\scriptsize{$2$}}\\[-0.10em] \fbox{\scriptsize{$4$}}\end{array}\!,\,\fbox{\scriptsize{$1$}}
\fbox{\scriptsize{$3$}}\right)$}}}\ ,\quad
{\cal{X}}_6:={\cal{X}}_{\textrm{\tiny{$\left(\!\!\begin{array}{c}\fbox{\scriptsize{$3$}}\\[-0.10em] \fbox{\scriptsize{$4$}}\end{array}\!,\,\fbox{\scriptsize{$1$}}
\fbox{\scriptsize{$2$}}\right)$}}}\ .\end{array}\]

The formulas (\ref{rel-a1})-(\ref{rel-a2}) take the form
\[\begin{array}{ll}
\left(\sigma_1+\frac{(q-q^{-1})v_1q^{-2}}{v_1-v_1q^{-2}}\right){\cal{X}}_1=0\, ,&
\left(\sigma_1+\frac{(q-q^{-1})v_2}{v_1-v_2}\right){\cal{X}}_2={\cal{X}}_4\left(\sigma_1+\frac{(q-q^{-1})v_1}{v_2-v_1}\right)\, ,\\[1em]
\left(\sigma_1+\frac{(q-q^{-1})v_2}{v_1-v_2}\right){\cal{X}}_3={\cal{X}}_5\left(\sigma_1+\frac{(q-q^{-1})v_1}{v_2-v_1}\right)\, ,&
\left(\sigma_1+\frac{(q-q^{-1})v_1}{v_2-v_1}\right){\cal{X}}_4={\cal{X}}_2\left(\sigma_1+\frac{v_2}{(q-q^{-1})v_1-v_2}\right)\, ,\\[1em]
\left(\sigma_1+\frac{(q-q^{-1})v_1}{v_2-v_1}\right){\cal{X}}_5={\cal{X}}_3\left(\sigma_1+\frac{(q-q^{-1})v_2}{v_1-v_2}\right)\, ,&
\left(\sigma_1+\frac{(q-q^{-1})v_2q^2}{v_2-v_2q^2}\right){\cal{X}}_6=0\, ,
\end{array}\] 

\vskip .1cm
\[\begin{array}{ll}
\left(\sigma_2+\frac{(q-q^{-1})v_2}{v_1q^{-2}-v_2}\right){\cal{X}}_1={\cal{X}}_2\left(\sigma_2+\frac{(q-q^{-1})v_1q^{-2}}{v_2-v_1q^{-2}}\right)\, ,&
\left(\sigma_2+\frac{(q-q^{-1})v_1q^{-2}}{v_2-v_1q^{-2}}\right){\cal{X}}_2={\cal{X}}_1\left(\sigma_2+\frac{(q-q^{-1})v_2}{v_1q^{-2}-v_2}\right)\, ,\\[1em]
\left(\sigma_2+\frac{(q-q^{-1})v_2q^2}{v_2-v_2q^2}\right){\cal{X}}_3=0\, ,&
\left(\sigma_2+\frac{(q-q^{-1})v_1q^{-2}}{v_1-v_1q^{-2}}\right){\cal{X}}_4=0\, ,\\[1em]
\left(\sigma_2+\frac{(q-q^{-1})v_2q^2}{v_1-v_2q^2}\right){\cal{X}}_5={\cal{X}}_6\left(\sigma_2+\frac{(q-q^{-1})v_1}{v_2q^2-v_1}\right)\, ,&
\left(\sigma_2+\frac{(q-q^{-1})v_1}{v_2q^2-v_1}\right){\cal{X}}_6={\cal{X}}_5\left(\sigma_2+\frac{(q-q^{-1})v_2q^2}{v_1-v_2q^2}\right)\, ,
\end{array}\] 

\vskip .1cm
\[\begin{array}{ll}
\left(\sigma_3+\frac{(q-q^{-1})v_2q^2}{v_2-v_2q^2}\right){\cal{X}}_1=0\, ,&
\left(\sigma_3+\frac{(q-q^{-1})v_2q^2}{v_1q^{-2}-v_2q^2}\right){\cal{X}}_2={\cal{X}}_3\left(\sigma_3+\frac{(q-q^{-1})v_1q^{-2}}{v_2q^2-v_1q^{-2}}\right)\, ,\\[1em]
\left(\sigma_3+\frac{(q-q^{-1})v_1q^{-2}}{v_2q^2-v_1q^{-2}}\right){\cal{X}}_3={\cal{X}}_2\left(\sigma_3+\frac{(q-q^{-1})v_2q^2}{v_1q^{-2}-v_2q^2}\right)\, ,&
\left(\sigma_3+\frac{(q-q^{-1})v_2q^2}{v_1q^{-2}-v_2q^2}\right){\cal{X}}_4={\cal{X}}_5\left(\sigma_3+\frac{(q-q^{-1})v_1q^{-2}}{v_2q^2-v_1q^{-2}}\right)\, ,\\[1em]
\left(\sigma_3+\frac{(q-q^{-1})v_1q^{-2}}{v_2q^2-v_1q^{-2}}\right){\cal{X}}_5={\cal{X}}_6\left(\sigma_3+\frac{(q-q^{-1})v_2q^2}{v_1q^{-2}-v_2q^2}\right)\, ,&
\left(\sigma_3+\frac{(q-q^{-1})v_1q^{-2}}{v_1-v_1q^{-2}}\right){\cal{X}}_6=0\, ,
\end{array}\] 
and
\[\begin{array}{lll}
(\tau-v_1){\cal{X}}_1=0\ ,&
(\tau-v_1){\cal{X}}_2=0\ ,&
(\tau-v_1){\cal{X}}_3=0\ ,\\[1em]
(\tau-v_2){\cal{X}}_4=0\ ,&
(\tau-v_2){\cal{X}}_5=0\ ,&
(\tau-v_2){\cal{X}}_6=0\ .
\end{array}\]

\vskip .2cm
The matrices corresponding to the action (\ref{rep-a1})-(\ref{rep-a2}) of the generators of $H(2,1,4)$ in the basis above are given by:

\[\sigma_1\mapsto\left(\begin{array}{cccccc}-q^{-1}&0&0&0&0&0\\[1em]
0&-\frac{(q-q^{-1})v_2}{v_1-v_2}&0&\frac{qv_1-q^{-1}v_2}{v_1-v_2}&0&0\\[1em]
0&0&-\frac{(q-q^{-1})v_2}{v_1-v_2}&0&\frac{qv_1-q^{-1}v_2}{v_1-v_2}&0\\[1em]
0&\frac{qv_2-q^{-1}v_1}{v_2-v_1}&0&-\frac{(q-q^{-1})v_1}{v_2-v_1}&0&0\\[1em]
0&0&\frac{qv_2-q^{-1}v_1}{v_2-v_1}&0&-\frac{(q-q^{-1})v_1}{v_2-v_1}&0\\[1em]
0&0&0&0&0&q\end{array}\right)\ ,\] 

\vskip 1.2cm
\[\sigma_2\mapsto\left(\begin{array}{cccccc}-\frac{(q-q^{-1})v_2}{v_1q^{-2}-v_2}&\frac{v_1q^{-1}-q^{-1}v_2}{v_1q^{-2}-v_2}&0&0&0&0\\[1em]
\frac{qv_2-v_1q^{-3}}{v_2-v_1q^{-2}}&-\frac{(q-q^{-1})v_1q^{-2}}{v_2-v_1q^{-2}}&0&0&0&0\\[1em]
0&0&q&0&0&0\\[1em]
0&0&0&-q^{-1}&0&0\\[1em]
0&0&0&0&-\frac{(q-q^{-1})v_2q^2}{v_1-v_2q^2}&\frac{qv_1-v_2q}{v_1-v_2q^2}\\[1em]
0&0&0&0&\frac{v_2q^3-q^{-1}v_1}{v_2q^2-v_1}&-\frac{(q-q^{-1})v_1}{v_2q^2-v_1}\end{array}\right)\ ,\] 

\vskip 1.2cm
\[\sigma_3\mapsto\left(\begin{array}{cccccc}q&0&0&0&0&0\\[1em]
0&-\frac{(q-q^{-1})v_2q^2}{v_1q^{-2}-v_2q^2}&\frac{v_1q^{-1}-v_2q}{v_1q^{-2}-v_2q^2}&0&0&0\\[1em]
0&\frac{v_2q^3-v_1q^{-3}}{v_2q^2-v_1q^{-2}}&-\frac{(q-q^{-1})v_1q^{-2}}{v_2q^2-v_1q^{-1}}&0&0&0\\[1em]
0&0&0&-\frac{(q-q^{-1})v_2q^2}{v_1q^{-2}-v_2q^2}&\frac{v_1q^{-1}-v_2q}{v_1q^{-2}-v_2q^2}&0\\[1em]
0&0&0&\frac{v_2q^3-v_1q^{-3}}{v_2q^2-v_1q^{-2}}&-\frac{(q-q^{-1})v_1q^{-2}}{v_2q^2-v_1q^{-2}}&0\\[1em]
0&0&0&0&0&-q^{-1}\end{array}\right)\ ,\] 

\vskip .5cm
\[\tau\mapsto\text{diag}(v_1,v_1,v_1,v_2,v_2,v_2)\ .\]

The Gram matrix $\langle{\cal{X}}_i,{\cal{X}}_j\rangle_{i,j=1,\dots,6}\ $, given by (\ref{H-scal-prod1})--(\ref{H-scal-prod2}), 
is $\text{diag}(z_1,z_2,z_3,z_4,z_5,z_6)$ where
$$\begin{array}{c}z_1=\frac{(q^{-2}v_1-q^2v_2)(q^{-3}v_1-q^3v_2)}{(v_1-v_2)(q^{-1}v_1-qv_2)}\,,\ \ 
z_2=\frac{q^{-3}v_1-q^3v_2}{q^{-1}v_1-qv_2}\,,\ \ z_3=1\,,\ \ 
z_4=\frac{(qv_1-q^{-1}v_2)(q^{-3}v_1-q^3v_2)}{(q^{-1}v_1-qv_2)^2}\,,\\[1.3em]
z_5=\frac{qv_1-q^{-1}v_2}{q^{-1}v_1-qv_2}\,,\ \ 
z_6=\frac{(v_1-v_2)(qv_1-q^{-1}v_2)}{(q^{-1}v_1-qv_2)(q^{-2}v_1-q^2v_2)}\ .
\end{array}$$

\paragraph{Acknowledgements.} We thank S. Khoroshkin, A. Kleshchev and E. Opdam for useful discussions.

\end{document}